\documentclass{amsart}

\usepackage[alphabetic]{amsrefs}
\usepackage{amsmath,amssymb,amsthm,enumerate, amscd}
\usepackage{tikz-cd}
\usepackage{tikz,graphicx}
\usepackage[all,cmtip]{xy}

\makeindex 
\begin{document}

%\externalcitedocument{Lit}

\pagenumbering{arabic}

 \setlength{\parskip}{7pt}

\def\LU{{\mathcal L\mathcal U}}
 \def\dach{\!\widehat{\phantom{G}}}
\def\Aut{\operatorname{Aut}}
\def\Int{\operatorname{Int}}
\def\GL{\operatorname{GL}}  
\def\QQ{{\mathbb{Q}}}
\def\Prim{\operatorname{Prim}}
\def\tr{\operatorname{tr}}
\def\Inn{\operatorname{Inn}} \def\Ad{\operatorname{Ad}}
\def\id{\operatorname{id}} \def\supp{\operatorname{supp}}
\def\csp{\overline{\operatorname{span}}}
\def\Ind{\operatorname{Ind}} \def\nd{\mathrm{es}}
\def\triv{\mathrm{triv}} 
\def\eps{\epsilon} 
\newcommand\veps{\varepsilon}
\def\K{{\mathcal K}}
\def\L{\mathcal L} \def\R{\mathcal R} \def\L{\mathcal L}
\def\C{\mathcal C} \def\E{\mathcal E} \def\Q{\mathcal Q}
\def\F{\mathcal F} \def\B{\mathcal B}
\def\H{\mathcal H}
\def\PU{\mathcal PU}
\def\M{\mathcal M}  \def\V{\mathcal V}

\def\EE{\mathbb E} \def\DD{\mathbb
    D} \def\I{{\mathcal I}} \def\U{{\mathcal U}} \def\UM{{\mathcal
      U}M} \def\ZUM{{\mathcal Z}UM} 
\newcommand{\ZM}{\mathcal{ZM}}
\def\KK{{KK}} \def\RK{\operatorname{RK}}
\def\RKK{\operatorname{RKK}} \def\ind{\operatorname{ind}}
\def\res{\operatorname{res}} \def\inf{\operatorname{inf}}
\def\ker{\operatorname{ker}} \def\infl{\operatorname{inf}}
\def\pt{\operatorname{pt}} 
\newcommand{\comp}{\operatorname{comp}}
\def\infl{\operatorname{infl}}
\def\IND{\operatorname{IND}}
\def\Bott{\operatorname{Bott}}
\def\Mor{\operatorname{Mor}}
\def\Rep{\operatorname{Rep}}
\def\Incr{\operatorname{Incr}}
\def\sp{\operatorname{span}}
\def\Hom{\operatorname{Hom}}
\def\LH{\operatorname{LH}}
\def\k{\operatorname{K}}
\def\EEG{\underline{\underline{\mathcal{E}G}}} 
\def\EG{\underline{EG}}
\def\EH{\underline{EH}}
\def\EGT{\underline{E\tilde{G}}}
\def\EGN{\underline{E G/N}} 
\renewcommand{\top}{\operatorname{top}}
\newcommand{\A}{\mathcal{A}} 
\newcommand{\D}{\mathcal{D}} 
\def\sm{\backslash} \def\Ind{\operatorname{Ind}}
\def\TT{\mathbb T} \def\ZZ{\mathbb Z} \def\CC{\mathbb C}
\def\FF{\mathbb F}
\def\RR{\mathbb R} \def\NN{\mathbb N} \def\om{\omega}
\def\ClmV{C_{\!\!\!{}_{-V}}} \def\ClV{C_{\!{}_{V}}}
\newcommand{\rk}{\rangle}
\newcommand{\lk}{\langle}
\newcommand{\bmtr}{\left(\begin{matrix}}
\newcommand{\emtr}{\end{matrix}\right)}
\newcommand*{\into}{\hookrightarrow}
\newcommand{\op}{\operatorname{op}}
\newcommand*{\Corr}{\mathfrak{Corr}}% category of C*-algebras with homomorphisms
\newcommand{\Cl}{\operatorname{\it{Cl}}}
\newcommand{\sign}{\operatorname{sign}}
\newcommand{\diag}{\operatorname{diag}}
\newcommand{\indx}{\operatorname{index}}
\newcommand{\one}{\mathrm 1}
\newcommand{\Homeo}{\operatorname{Homeo}}
\newcommand{\SL}{\operatorname{SL}}
\newcommand{\SO}{\operatorname{SO}}
\newcommand{\Ort}{\operatorname{O}}
\newcommand{\inte}{\operatorname{int}}
\newcommand{\spn}{\operatorname{span}}
\newcommand{\cspn}{\overline{\operatorname{span}}}
\newcommand{\Spec}{\operatorname{Spec}}
\newcommand{\Om}{\Omega}
\newcommand*{\onto}{\twoheadrightarrow}

\newcommand{\ddS}{\stackrel{\scriptscriptstyle{o}}{S}}
\renewcommand{\oplus}{\bigoplus}

\newtheorem{theorem}{Theorem}[section] \newtheorem{corollary}[theorem]{Corollary}
\newtheorem{lemma}[theorem]{Lemma} \newtheorem{proposition}[theorem]{Proposition}
\newtheorem{lemdef}[theorem]{Lemma and Definition}

\theoremstyle{definition} \newtheorem{definition}[theorem]{Definition}
\theoremstyle{definition} \newtheorem{notation}[theorem]{Notations}

\theoremstyle{remark} \newtheorem{remark}[theorem]{\bf Remark}
\newtheorem{example}[theorem]{\bf Example}
\newtheorem{exercise}[theorem]{\bf Exercise}

\numberwithin{equation}{section} \emergencystretch 25pt

\renewcommand{\theenumi}{\roman{enumi}}
\renewcommand{\labelenumi}{(\theenumi)}

\newpage

\title{Crossed products and the Mackey-Rieffel-Green machine}

\author{Siegfried Echterhoff}
\email{echters@uni-muenster.de}
\address{Mathematisches Institut\\
 Westf\"alische Wilhelms-Universit\"at M\"un\-ster\\
 Einsteinstr.\ 62\\
 48149 M\"unster\\
 Germany}

\begin{abstract} We give an introduction into the ideal structure and representation theory 
 of crossed products by actions of locally compact groups on $C^*$-algebras. In particular, we discuss the Mackey-Rieffel-Green theory of induced representations of crossed products and groups. 
Although we do not give complete proofs of all results, we try at least to explain the main ideas.
For a much more detailed exposition of many of the results presented here we refer to the beautiful book \cite{Wi-crossed} by Dana Williams.
   \end{abstract}
\maketitle

\section{Introduction}\footnote{The content of this note will appear 
as Chapter 2 of the book ``$K$-theory for group C*-algebras and semigroup C*-algebras'' 
which will appear in the Oberwolfach-seminar series of the Birkhäuser publishing company.
The research for this paper has been supported by the DFG through  CRC 878 groups, geometry \& actions.}

If a locally compact group $G$ acts continuously via
$*$-automorphisms on a $C^*$-algebra $A$, one can form
the full and reduced crossed products $A\rtimes G$ and
$A\rtimes_rG$ of $A$ by $G$.
The full crossed product should be thought of as a
skew {\em maximal} tensor product of $A$ with the
full group $C^*$-algebra $C^*(G)$ of $G$ and the
reduced crossed product should be regarded as a skew
{\em minimal} (or spacial) tensor product of
$A$ by the reduced group $C^*$-algebra $C_r^*(G)$ of $G$.

The crossed product construction provides a major 
source of examples in $C^*$-algebra theory, and
it plays  an important r\^ole in many applications
of $C^*$-algebras in other fields of mathematics,
like group representation theory and topology--here in particular in connection with the Baum-Connes conjecture,
which we shall treat in  \cite{KK}.
It is the purpose of this article to present in
a concise way some of the most important constructions
and features of crossed products with emphasis on the 
Mackey-Rieffel-Green machine as a basic technique
to investigate the ideal structure of crossed products. 
The contents of the first six sections of this chapter are also basic 
for the understanding of the contents of  \cite{KK}. 
Detailed proofs of most of the results on crossed products presented in this chapter
(if not given here) can be found in the  monograph \cite{Wi-crossed} by Dana Williams.
Note that the material covered in this article is almost perpendicular
to the material covered in Pedersen's book
\cite{Ped}. Hence we recommend the interested reader to 
also have a look into \cite{Ped} to obtain a more complete and balanced picture
of the theory. Pedersen's book also provides a good introduction 
into the general theory of $C^*$-algebras.
 An incomplete list of other  good references
on the general theory of $C^*$-algebas 
is  \cites{Bla06, Dav, Dix, Mur}. 
The Morita (or correspondence) category has been studied in more detail in 
\cite{EKQR2}.

Some general notation: if $X$ is a locally compact Hausdorff space
and $E$ is a  normed linear space, then we denote be $C_b(X,E)$
the space of bounded continuous $E$-valued functions on $X$ 
and by $C_c(X,E)$ and $C_0(X,E)$ those functions in $C_b(X,E)$ which have 
compact supports or which vanish at infinity. If $E=\CC$, then we simply write
$C_b(X), C_c(X)$ and $C_0(X)$, respectively.
If $E$ and $F$ are two linear spaces, then $E\odot F$ always denotes the 
algebraic tensor product of $E$ and $F$ and we reserve the sign ``$\otimes$''
for certain kinds of topological tensor products.

\section{Some preliminaris}\label{sec-prel}

We shall assume throughout this article that the reader is familiar
with the basic concepts of $C^*$-algebras as can be found in any 
of the standard text books mentioned above. 
However, in order to make this treatment more self-contained 
 we try to recall some basic facts
and notation on $C^*$-algebras which will play an important r\^ole in this article.

\subsection{$C^*$-algebras}\label{subseccstar}
A (complex) $C^*$-algebra is a complex Banach-algebra $A$ together with
an involution $a\mapsto a^*$ such that $\|a^*a\|=\|a\|^2$ for all $a\in A$.
Note that we usually do not assume that $A$ has a unit.
Basic examples are given by the algebras $C_0(X)$ and $C_b(X)$
equipped with the supremum-norm and the involution $f\mapsto \bar{f}$. 
These algebras are clearly commutative and a classical theorem of Gelfand and Naimark
asserts, that {\em all} commutative $C^*$-algebras are isomorphic to some $C_0(X)$
(see paragraph \ref{subsec-comm} below).
Other examples are given by the algebras $\B(H)$ of bounded operators 
on a Hilbert space $H$ with operator norm and involution given by taking 
the adjoint operators, and all closed $*$-subalgebras of $\B(H)$
(like the algebra $\K(H)$ of compact operators on $H$).
Indeed, another classical result by Gelfand and Naimark shows
that every $C^*$-algebra is isomorphic to a closed $*$-subalgebra of some $\B(H)$.
If $S\subseteq A$ is any subset of a $C^*$-algebra $A$, we denote by 
$C^*(S)$ the smallest sub-$C^*$-algebra of $A$ which contains $S$.
A common way to construct $C^*$-algebras is by describing a certain 
set $S\subseteq \B(H)$ and forming the algebra $C^*(S)\subseteq \B(H)$.
If $S=\{a_1, \ldots, a_l\}$ is a finte set of elements of $A$, we shall also write
$C^*(a_1,\ldots, a_l)$ for $C^*(S)$. For example, if $U,V\in \B(H)$ are unitary operators
such that $UV=e^{2\pi i \theta}VU$ for some irrational $\theta\in [0,1]$, then 
$A_{\theta}:=C^*(U,V)$ is the well known {\em irrational rotation algebra} corresponding to $\theta$,
a standard example in $C^*$-algebra theory (in this example one can show that the isomorphism class 
of $C^*(U,V)$ does not depend on the particular choice of $U$ and $V$).

$C^*$-algebras are very rigid objects: If $A$ is a $C^*$-algebra, then every closed 
(two-sided) ideal of $A$ is automatically
selfadjoint and $A/I$, equipped with the obvious operations and the quotient norm
is again a $C^*$-algebra.
If $B$ is any Banach $*$-algebra (i.e., a Banach algebra 
with isometric involution, which does not necessarily satisfy the $C^*$-relation $\|b^*b\|=\|b\|^2$),
and if $A$ is a $C^*$-algebra, then any $*$-homomorphism $\Phi:B\to A$ 
is automatically continuous with $\|\Phi(b)\|\leq\|b\|$ for all $b\in B$.
If $B$ is also a $C^*$-algebra, then $\Phi$ factors through an isometric $*$-homomorphism
$\tilde{\Phi}:B/(\ker\Phi)\to A$. In particular, if $A$ and $B$ are $C^*$-algebras and 
$\Phi:B\to A$ is an injective (resp. bijective) $*$-homomorphism, then $\Phi$ 
is automatically isometric (resp. an isometric isomorphism).

\subsection{Multiplier Algebras} 
The {\em multiplier algebra} $M(A)$ of a $C^*$-algebra $A$
is the largest $C^*$-algebra
which contains $A$ as an essential ideal (an ideal $J$ of a $C^*$-algebra $B$ is called
{\em essential} if for all $b\in B$ we have $bJ=\{0\} \Rightarrow b=0$).
If $A$ is represented faithfully and non-degenerately
on a Hilbert space $H$ (i.e. $A\subseteq \B(H)$ with $AH=H$), 
then $M(A)$ can be realized as
the idealizer
$$M(A)=\{T\in \B(H): TA\cup AT\subseteq A\}$$
of $A$ in $\B(H)$.
In particular we have $M(\K(H))=\B(H)$, where $\K(H)$
denotes the algebra of compact operators on $H$.

The {\em strict topology} on $M(A)$ is the locally convex topology
generated by the semi-norms $m\mapsto \|am\|, \|ma\|$ with $a\in A$.
Notice that $M(A)$ is the strict completion of $A$.
$M(A)$ is always unital and  $M(A)=A$ if (and only if) $A$ is unital.
If $A=C_0(X)$ for some locally compact space 
$X$, then $M(A)\cong C_b(X)\cong C(\beta(X))$, where $\beta(X)$ denotes the 
Stone-\v Cech compactification of $X$. Hence $M(A)$ should be viewed as 
a noncommutative analogue of the Stone-\v Cech compactification.
If $A$ is any $C^*$-algebra, then the algebra $A_1:=C^*(A\cup\{1\})\subseteq M(A)$ 
is called the {\em unitization} of $A$ (notice that $A_1=A$ if $A$ is unital). If
$A=C_0(X)$ for some non-compact $X$, then $A_1\cong C(X_+)$,
where $X_+$ denotes the one-point compactification of $X$.

A $*$-homomorphism
$\pi:A\to M(B)$ is called {\em non-degenerate} if
$\pi(A)B=B$, which by Cohen's factorization theorem is equivalent
to the weaker condition that $\sp\{\pi(a)b:a\in A, b\in B\}$ is
dense in $B$ (e.g. see \cite{RW}*{Proposition 2.33}). 
If $H$ is a Hilbert space, then  $\pi:A\to M(\K(H))=\B(H)$ 
is non-degenerate in the above sense iff $\pi(A)H=H$.
If $\pi:A\to M(B)$ is non-degenerate, then there exists
a unique $*$-homomorphism  $\bar{\pi}:M(A)\to M(B)$
such that $\bar{\pi}|_A=\pi$. We shall usually make no
notational difference between $\pi$ and its extension
$\bar{\pi}$. 

\subsection{Commutative $C^*$-algebras and functional calculus}\label{subsec-comm}
If $A$ is commutative, then we denote by  $\Delta(A)$ the set of all
non-zero algebra homomorphisms $\chi:A\to\CC$ equipped 
 with the weak-$*$ topology.
 Then $\Delta(A)$ is locally compact and it is compact if $A$ is unital. 
 If $a\in A$, then $\widehat{a}:\Delta(A)\to \CC; \widehat{a}(\chi):=\chi(a)$
is an element of $C_0(\Delta(A))$, and the 
Gelfand-Naimark theorem asserts that
$A\to C_0(\Delta(A)): a\mapsto \widehat{a}$
is an (isometric) $*$-isomorphism. 

If $A$ is any $C^*$-algebra, then an element $a\in A$ is called {\em normal} if $a^*a=aa^*$.
If $a\in A$ is normal, then $C^*(a,1)\subseteq A_1$ is a commutative sub-$C^*$-algebra 
of $A_1$. Let $\sigma(a)=\{\lambda\in \CC: a-\lambda 1 \;\text{is not invertible in $A_1$}\}$
denote the {\em spectrum} of $a$, a nonempty compact subset of $\CC$. 
If $\lambda\in \sigma(a)$, then $a-\lambda 1$ generates 
a unique maximal ideal $M_\lambda$ of $C^*(a,1)$ and the quotient map
$C^*(a,1)\to C^*(a,1)/M_\lambda\cong \CC$ determines an element 
$\chi_\lambda\in \Delta(C^*(a,1))$.
One then checks that $\lambda\mapsto \chi_\lambda$ is a homeomorphism between
$\sigma(a)$ and $\Delta(C^*(a,1))$. Thus, the Gelfand-Naimark theorem provides 
a $*$-isomorphism $\Phi:C(\sigma(a))\to  C^*(a,1)$. 
If $p(z)=\sum_{i,j=0}^n \alpha_{ij}z^i\bar{z}^j$ is a polynomial in $z$ and $\bar{z}$ 
(which by the Stone-Weierstra{\ss}  theorem form
a dense subalgebra of $C(\sigma(a))$), then $\Phi(p)=\sum_{i,j=0}^n \alpha_{ij}a^i(a^*)^j$.
In particular, we have $\Phi(1)=1$ and $\Phi(\id_{\sigma(a)})=a$.
In what follows, we always write $f(a)$ for $\Phi(f)$. Note that $\sigma(f(a))=f(\sigma(a))$ 
and if $g\in C\big(\sigma(f(a))\big)$, then $g(f(a))=(g\circ f)(a)$, i.e., the functional calculus
is compatible with composition of functions.
If $A$ is not unital, then $0\in \sigma(a)$ and it is clear that for any polynomial $p$ in $z$ and 
$\bar{z}$ we have $p(a)\in A$ if and only if $p(0)=0$. Approximating functions by polynomials,
it follows that $f(a)\in A$ if and only if $f(0)=0$ and we obtain an isomorphism
$C_0(\sigma(a)\setminus\{0\})\to C^*(a)\subseteq A; f\mapsto f(a)$.

\begin{example}
An element $a\in A$ is called {\em positive} if $a=b^*b$ for some $b\in A$. This is equivalent
to saying that $a$ is selfadjoint (i.e. $a=a^*$) and $\sigma(a)\subseteq[0,\infty)$.  If $a\geq 0$, then the functional calculus
provides  the element $\sqrt{a}\in A$, which is the unique positive element of 
$A$ such that $(\sqrt{a})^2=a$. If $a\in A$ is selfadjoint, then 
$\sigma(a)\subseteq\RR$ and the functional
calculus allows a unique decomposition $a= a_+- a_-$ with $a_+, a_-\geq 0$ such that $a_+\cdot a_-=0$.
Simply take $a_+=f(a)$ with $f(t)=\max\{t,0\}$.
Since we can write any $b\in A$ as a linear combination of two selfadjoint elements
via $b=\frac{1}{2}(a+a^*) +i  \frac{1}{2i}(a-a^*)$, we see that every element of $A$
can be written as a linear combination of four positive elements. 
Since every positive element is a square, it follows that $A=A^2:= \LH\{ab: a,b\in A\}$
(Cohen's factorization theorem even implies that $A=\{ab: a,b\in A\}$).
\end{example}

Every $C^*$-algebra has an {\em approximate unit}, i.e., a net  $(a_i)_{i\in I}$ 
in $A$ such that $\|a_ia-a\|, \|a a_i-a\|\to 0$ for all $a\in A$.
In fact $(a_i)_{i\in I}$ can be chosen so that
$a_i\geq 0$ and $\|a_i\|=1$ for all $i\in I$.
If $A$ is separable ( i.e., $A$ contains a countable dense set), then 
one can find a sequence $(a_n)_{n\in\NN}$ with these properties.

If $A$ is a unital $C^*$-algebra,
then $u\in A$ is called a {\em unitary}, if $uu^*=u^*u=1$. If $u$ is unitary, then 
$\sigma(u)\subseteq\TT=\{z\in \CC: |z|=1\}$ and hence $C^*(u)=C^*(u,1)$ is isomorphic 
to a quotient of $C(\TT)$. Note that if $u,v\in A$ are two unitaries such that
$uv=e^{2\pi i\theta}vu$ for some irrational $\theta\in [0,1]$, then one can show that
$\sigma(u)=\sigma(v)=\TT$, so that $C^*(u)\cong C^*(v)\cong C(\TT)$.
It follows that the irrational rotation algebra $A_\theta=C^*(u,v)$ should be regarded
as (the algebra of functions on) a ``noncommutative product'' of two tori
which results in the expression of a {\em noncommutative $2$-torus}.

\subsection{Representation and ideal spaces of $C^*$-algebras}\label{subsec-spectrum}
If $A$ is a $C^*$-algebra, the {\em spectrum} $\widehat{A}$ is defined as the 
set of all unitary equivalence classes of irreducible representations 
$\pi:A\to \B(H)$ of $A$ on Hilbert space\footnote{A selfadjoint subset $S\subseteq \B(H)$ is called 
{\em irreducible} if there exists no proper nontrivial closed subspace $L\subseteq H$
with $SL\subseteq L$. By Schur's lemma, this is equivalent to saying that the commutator 
of $S$ in $\B(H)$ is equal to $\CC\cdot 1$. A representation $\pi:A\to \B(H)$ is {\em irreducible}
if $\pi(A)$ is irreducible. Two representations $\pi,\rho$ of $A$ on $H_\pi$ and $H_{\rho}$, respectively,
are called unitarily equivalent, if there exists a unitary $V:H_\pi\to H_\rho$ such that
$V\circ \pi(a)=\rho(a)\circ V$ for all $a\in A$.}.
We shall usually make no notational difference between 
an  irreducible representation $\pi$ and its equivalence class $[\pi]\in \widehat{A}$.
The {\em primitive ideals} of $A$ are the kernels of the irreducible representations
of $A$, and we write
$\Prim(A):=\{\ker\pi: \pi\in \widehat{A}\}$ for the set of all
 primitive ideals  of $A$. 
Every closed 
two-sided ideal $I$ of $A$ is an intersection of primitive ideals.
The spaces $\widehat{A}$ and $\Prim(A)$ 
are equipped with the {\em Jacobson topologies}, where the closure
operations are given by $\pi\in \overline{R}:\Leftrightarrow
\ker\pi\supseteq \cap\{\ker\rho:\rho\in R\}$ (resp.
$P\in \overline{R}:\Leftrightarrow
P\supseteq \cap\{Q:Q\in R\}$) for $R\subseteq\widehat{A}$ 
 (resp.  $R\subseteq\Prim(A)$).
In general, the Jacobson topologies are far away 
from being Hausdorff. In fact, while $\Prim(A)$ is at least always a 
T$_0$-space (i.e. for any two different elements in $\Prim(A)$ at least one 
of them has an open neighborhood which does not contain the other), 
this very weak separation property often fails for the space $\widehat{A}$.
If $A$ is commutative, it follows from Schur's lemma that $\widehat{A}=\Delta(A)$
and the Jacobson topology coincides in this case with the weak-$*$ topology.

If $I$ is a closed two-sided ideal of $A$, then $\widehat{A}$
can be identified with the disjoint union of $\widehat{I}$ with 
$\widehat{A/I}$, such that 
$\widehat{I}$ identifies with $\{\pi\in \widehat{A}: \pi(I)\neq \{0\}\}\subseteq
\widehat{A}$ and $\widehat{A/I}$ identifies
with $\{\pi\in \widehat{A}: \pi(I)=\{0\}\}\subseteq\widehat{A}$.
It follows from the definition of the Jacobson topology that 
$\widehat{A/I}$ is closed and $\widehat{I}$ is open in $\widehat{A}$.
The correspondence $I\leftrightarrow \widehat{I}$ (resp
$I\leftrightarrow \widehat{A/I}$) is a one-to-one 
correspondence between the closed two-sided ideals of $A$
and the open (resp. closed) subsets of $\widehat{A}$. 
Similar statements hold for the open or closed subsets of $\Prim(A)$.

A $C^*$-algebra is called {\em simple} if $\{0\}$ is the only proper closed 
two-sided ideal of $A$. Of course, this is equivalent to saying that $\Prim(A)$ has 
only one element (the zero ideal). Simple $C^*$-algebras are thought of as the
basic ``building blocks'' of more general $C^*$-algebras.
Examples of simple algebras are 
the algebras $\K(H)$ of compact operators on a Hilbert space $H$ and 
the irrational rotation algebras $A_\theta$.
Note that while $\widehat{\K(H)}$ has also only one element (the equivalance
class of its embedding into $\B(H)$), one can show that $\widehat{A}_\theta$
is an uncountable infinite set (this can actually be deduced from Proposition \ref{prop-irred} below).

A $C^*$-algebra $A$ is called {\em type I} (or {\em GCR}, or {\em postliminal})
if for every irreducible representation $\pi:A\to \B(H)$ 
we have $\pi(A)\supseteq \K(H)$. We refer to \cite{Dix}*{Chapter 12} for some important
equivalent characterizations of type I algebras.
A $C^*$-algebra $A$ is called {\em CCR} (or {\em liminal}), if $\pi(A)=\K(H)$ for every irreducible 
representation $\pi\in \widehat{A}$.
If $A$ is type I, then the mapping $\widehat{A}\to \Prim(A):\pi\mapsto \ker\pi$
is a homeomorphism, and the converse holds if $A$ is separable (in the 
non-separable case the question whether this converse holds leads to quite
interesting logical implications, e.g. see \cite{AW}).
Furthermore, if $A$ is type I, then $A$ is CCR if and only if $\widehat{A}\cong \Prim(A)$
is a T$_1$-space, i.e., points are closed.

A $C^*$-algebra is said to have 
 {\em continuous trace} if there exists a dense ideal 
$\mathfrak m\subseteq A$ such that for all positve elements $a\in \mathfrak m$ the 
operator $\pi(a)\in \B(H_\pi)$ is trace-class and the resulting map
$\widehat{A}\to [0,\infty); \pi\mapsto \tr(\pi(a))$ is continuous.
Continuous-trace algebras are all CCR with Hausdorff spectrum $\widehat{A}$.
Note that every type I $C^*$-algebra $A$ contains a non-zero closed two-sided ideal $I$ such 
that $I$ is a continuous-trace algebra (see \cite{Dix}*{Chapter 4}).
%
%Since quotients  of type I $$C^*$$-algebras are 
%again type I, it follows that every $C^*$-algebra has a (possibly transfinite) 
%ascending composition
%sequence of ideals $\{I_\alpha\}_{\alpha\leq \alpha_0}$ where the $\alpha$'s  denote
%ordinal numbers, such that $I_{\alpha_0}=A$, $I_{\alpha+1}/I_\alpha$ has continuous trace
%for all $\alpha<\alpha_0$ and $I_{\beta}=\cup_{\alpha<\beta} I_\alpha$ if $\beta\leq \alpha_0$
%is a limit ordinal. Thus the continuous trace algebras can be regarded as the 
%building stones of type I algebras.

\subsection{Tensor products}
The algebraic tensor product $A\odot B$ of two $C^*$-algebras
$A$ and $B$ has a canonical structure as a $*$-algebra.
To make it a $C^*$-algebra, we have to take completions
with respect to suitable cross-norms $\|\cdot \|_{\mu}$ satisfying
$\|a\otimes b\|_{\mu}=\|a\|\|b\|$.
Among the possible choices of such norms there is a maximal
cross-norm $\|\cdot\|_{\max}$ and a minimal cross-norm $\|\cdot\|_{\min}$
giving rise to the {\em maximal tensor product} $A\otimes_{\max}B$
and the {\em minimal tensor product} $A\otimes_{\min}B$
(which we shall always denote by $A\otimes B$).

The maximal tensor product is characterized
by the universal
property that  any commuting pair of $*$-homomorphisms
$\pi:A\to D, \rho:B\to D$ determines a $*$-homomorphism 
$\pi\times \rho:A\otimes_{\max} B\to D$ such that $\pi\times \rho(a\otimes b)=\pi(a)\rho(b)$
for all elementary tensors $a\otimes b\in A\odot B$. 
%$\pi:A\odot B\to \B(H)$
%   extends to the completion $A\otimes_{\max}B$.
The minimal (or spatial) tensor product
$A\otimes B$ is the completion of $A\odot B$ with respect
to
$$\left\|\sum_{i=1}^na_i\otimes b_i\right\|_{\min}=
\left\|\sum_{i=1}^n\rho(a_i)\otimes\sigma(b_i)
\right\|,$$
where $\rho:A\to \B(H_{\rho}), \sigma:B\to \B(H_{\sigma})$ are
faithful representations of $A$ and $B$
and the norm on the right is taken in $\B(H_{\rho}\otimes H_{\sigma})$.
It is a non-trivial fact (due to Takesaki) that
$\|\cdot\|_{\min}$ is the smallest cross-norm on
$A\odot B$ and that it does not depend on the
choice of $\rho$ and $\sigma$ (e.g. see \cite{RW}*{Theorem B.38}).

A $C^*$-algebra $A$ is called {\em nuclear}, if
$A\otimes_{\max}B=A\otimes B$
for all $B$. Every type I $C^*$-algebra is nuclear (e.g. see \cite{RW}*{Corollary B.49})  as well as 
the irrational rotation algebra $A_\theta$ (which will follow from Theorem \ref{thm-nuclear} below).  
In particular,
all commutative $C^*$-algebras are nuclear and
we have $C_0(X)\otimes B\cong C_0(X,B)$ for any locally compact space $X$.
One can show that $\B(H)$ is not nuclear if $H$
is an infinite dimensional Hilbert space. 

If $H$ is an infinite dimensional Hilbert space, then $\K(H)\otimes \K(H)$ 
is isomorphic to $\K(H)$ (which can be deduced from a unitary isomorphism
$H\otimes H\cong H$). A $C^*$-algebra $A$ is called {\em stable} if $A$ is isomorphic 
to $A\otimes \K$, where we write $\K:=\K(l^2(\NN))$. It 
follows from the associativity of taking tensor products that $A\otimes\K$ is always stable
and we call $A\otimes\K$ the {\em stabilisation} of $A$. 
Note that $A\otimes \K$ and $A$ have isomorphic representation and ideal spaces.
For example,  the map $\pi\mapsto \pi\otimes \id_\K$ gives a homeomorphism between
$\widehat{A}\to (A\otimes\K)\dach$. 
Moreover $A$ is type I (or CCR or continuous-trace or nuclear)
if and only if  $A\otimes \K$ is.

\section{Actions and their crossed products}
\label{sec-actions}

\subsection{Haar measure and vector-valued integration on groups}
If $X$ is a locally compact space, we denote by $C_c(X)$ the set of all
continuous functions with compact supports on $X$. 
A {\em positive  integral} on $C_c(X)$ is a linear functional
$\int :C_c(X)\to \CC$ such that $\int_X f(x)\,dx:=\int( f )\geq 0$ if $f\geq 0$.
We refer to \cite{RudRC} for a good treatment of the Riesz representation theorem
which provides a one-to-one  correspondence between integrals on $C_c(X)$ and
positive Radon measures on $X$. If $H$ is a Hilbert space and $f:X\to \B(H)$ is a
weakly continuous function 
(i.e., $x\mapsto \lk f(x)\xi,\eta\rk$ 
is continuous for all $\xi,\eta\in H$)
with compact support, then there exists a unique operator 
$\int_X f(x)\,dx\in \B(H)$ such that 
$$\left\lk \left(\int_X f(x)\,dx\right)\xi,\eta\right\rk=
\int_X\lk f(x)\xi,\eta\rk\,dx\quad\quad\text{for all}\; \xi,\eta\in H.$$
If $A$ is a $C^*$-algebra imbedded faithfully by a  non-degenerate representation into some $\B(H)$ and $f\in C_c(X,A)$
is {\em norm-continuous}, 
then approximating $f$ uniformly with controlled supports by elements in
the algebraic tensor product $C_c(X)\odot A$ shows that $\int_X f(x)\,dx \in A$.
Moreover, if $f:X\to M(A)$ is a strictly continuous function with compact support, then
(via the canonical embedding $M(A)\subseteq \B(H)$) $f$ is weakly continuous as a function into $\B(H)$,
and since $\big(x\mapsto a f(x), f(x) a\big)\in C_c(X,A)$ for all $a\in A$ it follows
that $\int_X f(x)\,dx\in M(A)$.

If $G$ is a locally compact group, then 
there exists a nonzero positive integral 
$\int:C_c(G)\to \CC$, called {\em Haar integral} on $C_c(G)$,
  such that  $\int_G f(gx)\,dx=\int_G f(x)\,dx$ for all $f\in C_c(G)$ and $g\in G$.
The Haar integral is unique up to multiplication with 
a positive number, which implies that  for each $g\in G$ there exists a positive 
number $\Delta(g)$ such that
$\int_G f(x)\, dx=\Delta(g)\int_G f(xg)\, dx$ for all
$f\in C_c(G)$ (since the right hand side of the equation  defines a new Haar integral). 
One can show that $\Delta:G\to (0,\infty)$ is  a continuous group homomorphism.
A group $G$ is called unimodular if $\Delta(g)=1$ for all $g\in G$. All discrete, all compact
and all abelian groups are unimodular, however, the $ax+b$-group, which is the 
semidirect product $\RR\rtimes \RR^*$ via the action of the multiplicative group
$\RR^*:=\RR\setminus\{0\}$ on the additive group 
$\RR$ by dilation, is not unimodular. As a general reference for the Haar integral we refer to
\cite{DE}.
%
% If $A$ is a $C^*$-algebra and
%$f\in C_c(G,A)$, the set of all compactly supported
%$A$-valued functions on $G$,  the integral
%$\int_G f(s)\,ds\in A$ can be most conviniently
%defined via representing $A$ faithfully (and non-degenerately)
%as bounded operators on some Hilbert space $H$ and then
%defining $\int_G f(s)\,ds$ as the unique bounded operator $T$ on $H$
%satisfying $\lk T\xi,\eta\rk=\int_G\lk f(s)\xi,\eta\rk\,ds$.
%Approximating $f$ uniformly by functions in the algebraic tensor product
%$C_c(G)\odot A$ (with controlled supports) shows that $T$ is an 
%element of $A$.

%More general,
%if $f:G\to M(A)$ is a  strictly continuous
%compactly supported
%function, and if $A$ (and hence $M(A)$)  is represented faithfully
%into $\B(H)$, then $f$ becomes a strongly
%continuous compactly supported function into $\B(H)$,
%and, as above, the  operator
%valued integral $T:=\int_G f(s)\, ds\in \B(H)$ is defined.
%Since for all $a\in A$ the functions $s\to f(s)a, af(s)$ are in $C_c(A)$,
%it follows that $\int_G f(s)\, ds\in M(A)$.
%This is basically all integration theory we shall
%need in this article. 

\subsection{$C^*$-dynamical systems and their crossed products}
An \emph{action} of a locally compact group
$G$ on  a $C^*$-algebra $A$ is a
homomorphism  $\alpha:G\to \Aut(A); s\mapsto \alpha_s$  of $G$ into the group
$\Aut(A)$ of $*$-automorphisms of $A$ such that
$s\mapsto \alpha_s(a)$ is continuous for all $a\in A$ (we then say that $\alpha$ 
is {\em strongly continuous}).
The triple $(A,G,\alpha)$ is then called a {\em $C^*$-dynamical system}
(or {\em covariant system}). We also often say that $A$ is a {\em $G$-algebra},
when $A$ is equipped with a given $G$-action $\alpha$. 

\begin{example}[Transformation groups]\label{ex-transf}
If $G\times X\to X; (s,x)\mapsto s\cdot x$ is a continuous action 
of $G$ on a locally compact Hausdorff space $X$, then $G$ acts on 
$C_0(X)$ by
$\big(\alpha_s(f)\big)(x):=f(s^{-1}\cdot x)$, and it is not difficult to see 
that every action on $C_0(X)$ arises in this way.
Thus, general $G$-algebras
are non-commutative analogues
of locally compact $G$-spaces.
\end{example}

If $A$ is a $G$-algebra, then
$C_c(G,A)$ becomes a $*$-algebra with respect to
convolution and involution defined by
\begin{equation}
f*g(s)=\int_G f(t)\alpha_t(g(t^{-1}s))\,dt\quad
\text{and}\quad
f^*(s)=\Delta(s^{-1})\alpha_s(f(s^{-1}))^*.
\end{equation}
%As for algebraic tensor products of two $C^*$-algebras, there are
%two canonical choices for taking $C^*$-completions of the
%$*$-algebra $C_c(G,A)$.
%To describe these completions
%we need the notion of 
A {\em covariant homomorphism}
of $(A,G,\alpha)$ into the multiplier algebra $M(D)$
of some $C^*$-algebra $D$ is a pair
$(\pi, U)$, where
$\pi:A\to M(D)$ is a $*$-homomorphism
and $U:G\to UM(D)$ is a strictly continuous homomorphism
into the group $UM(D)$ of unitaries in $M(D)$
satisfying
\[
\pi(\alpha_s(a))= U_s\pi(a)U_{s^{-1}}\quad
\text{for all}\;s\in G.
\]
We say that $(\pi, U)$ is {\em non-degenerate} if $\pi$ is non-degenerate.
A \emph{covariant representation} of $(A,G,\alpha)$ on a
Hilbert space $H$ is a covariant homomorphism into
$M(\K(H))=\B(H)$.  If $(\pi, U)$ is a covariant homomorphism into
$M(D)$, its \emph{integrated form} $\pi\times U:C_c(G,A)\to
M(D)$ is defined by
\begin{equation}
(\pi\times U)(f):=\int_G \pi(f(s))U_s\,ds\in M(D).
\end{equation}
It is straightforward to check that $\pi\times U$ is a $*$-homomorphism.

Covariant homomorphisms do exist. Indeed, if
$\rho:A\to M(D)$ is any
$*$-homomorphism, then we can construct the 
induced covariant homomorphism
$\Ind\rho:=(\tilde{\rho},1\otimes\lambda)$ of
$(A,G,\alpha)$ into $M\big(D\otimes \K(L^2(G))\big)$ as follows:
Let $\lambda:G\to U(L^2(G))$ denote the {\em left regular
representation} of $G$ given by $(\lambda_s\xi)(t)=\xi(s^{-1}t)$,
and define
$\tilde\rho$ as the composition
$$
\begin{CD}
      A @>\tilde{\alpha} >>  M\big(A\otimes C_0(G)\big)
      @>\rho\otimes M>> M\big(D\otimes \K(L^2(G))\big),
\end{CD}
$$
where the $*$-homomorphism  $\tilde\alpha:A\to C_b(G,A)\subseteq M\big(A\otimes C_0(G)\big)$$\;$\footnote{$C_b(G,A)$ is regarded as a subset of $M\big(A\otimes
C_0(G)\big)$ via the identification $A\otimes C_0(G)\cong C_0(G,A)$
and taking pointwise products of functions.}
is defined by  $\tilde\alpha(a)(s)=\alpha_{s^{-1}}(a)$, and where
$M:C_0(G)\to \B(L^2(G))=M(\K(L^2(G)))$ denotes 
the represention by multiplication
operators. We call $\Ind\rho$ the
{\em covariant homomorphism
induced from $\rho$}, and we shall make no notational difference 
between $\Ind\rho$
and its integrated form $\tilde\rho\times (1\otimes \lambda)$.
$\Ind\rho$ is faithful on $C_c(G,A)$
whenever $\rho$ is faithful
on $A$.
If $\rho=\id_A$, the identity on $A$, then we say that
$$\Lambda_A^G:=\Ind(\id_A): C_c(G,A)\to M\big(A\otimes \K(L^2(G))\big)$$
is {\em the regular representation} of $(A,G,\alpha)$. Notice that 
\begin{equation}\label{eq-regular}
\Ind\rho=(\rho\otimes\id_{\K})\circ \Lambda_A^G
\end{equation}
for all $*$-homomorphisms $\rho:A\to M(D)$.$\;$\footnote{This equation even makes sense if 
$\rho$ is degenerate since $\rho\otimes\id_\K$ is well defined on the image of 
$C_b(G,A)$ in $M(A\otimes\K(L^2(G)))$.}

\begin{remark}\label{rem-reg}
If we start with a representation $\rho:A\to \B(H)=M(\K(H))$ of $A$ 
on a Hilbert space $H$, then $\Ind\rho=(\tilde\rho,1\otimes \lambda)$ is the representation 
of $(A,G,\alpha)$ into $\B(H\otimes L^2(G))$ (which equals $M(\K(H)\otimes \K(L^2(G)))$)
given by the formulas
$$\big(\tilde\rho(a)\xi\big)(t)=\rho(\alpha_{t^{-1}}(a))(\xi(t))\quad
\text{and}\quad \big((1\otimes\lambda)(s)\xi\big)(t)=\xi(s^{-1}t),$$
for $a\in A, s\in G$ and $\xi\in L^2(G,H)\cong H\otimes L^2(G)$. 
Its integrated form is given by the convolution formula
$$f*\xi(t):=\big(\Ind\rho(f)\xi\big)(t)=
\int_G\rho\big(\alpha_{t^{-1}}(f(s))\big)\xi(s^{-1}t)\,ds$$
for $f\in C_c(G,A)$ and $\xi\in L^2(G,H)$. 
\end{remark}

\begin{definition}%[cf Chapter \ref{?}]
\label{def-crossedaction}
Let $(A,G,\alpha)$ be a $C^*$-dynamical system. 
\begin{enumerate}
      \item The {\em full crossed produt} $A\rtimes_{\alpha}G$
      (or just $A\rtimes G$ if $\alpha$ is understood) is the
      completion of $C_c(G,A)$ with respect to
      $$
\|f\|_{\max}:=\sup\{\|(\pi\times U)(f)\|: \text{$(\pi,U)$ is a
covariant representation
of $(A,G,\alpha)$}\}.
$$
\item The {\em reduced crossed product} $A\rtimes_{\alpha,r}G$
(or just $A\rtimes_rG$) is defined as
$$\overline{\Lambda_A^G\big(C_c(G,A)\big)}\subseteq M\big(A\otimes
\K(L^2(G))\big).$$
\end{enumerate}
\end{definition}

\begin{remark}
\label{rem-universal}
{\bf (1)} It follows directly from the above definition that every
integrated form
$\pi\times U:C_c(G,A)\to M(D)$ of a covariant homomorphism $(\pi,U)$
extends to a $*$-homomorphism
of $A\rtimes_{\alpha}G$ into $M(D)$. Conversely, every
{\bf non-degenerate} $*$-homomorphism $\Phi:A\rtimes_{\alpha}G\to M(D)$
is of the form $\Phi=\pi\times U$ for some non-degenerate covariant
homomorphism $(\pi,U)$. To see this consider
the canonical covariant homomorphism
$(i_A,i_G)$ of $(A,G,\alpha)$ into $M(A\rtimes_{\alpha}G)$ given
by the formulas
\begin{align*}
(i_A(a) f)(s) &= a f(s)
&
(i_G(t) f)(s) &= \alpha_t(f(t^{-1}s))
\\
(f i_A(a))(s) &= f(s) \alpha_s(a)
&
(f i_G(t))(s)
&= \Delta(t^{-1}) f(st^{-1}),
\end{align*}
$f\in C_c(G,A)$ (the given 
formulas extend to left and right multiplications of $i_A(a)$ and $i_G(s)$
with elements in $A\rtimes G$). It is then relatively easy to check that
$\Phi=\pi\times U$ with 
$$\pi=\Phi\circ i_A\quad\text{and}\quad U=\Phi\circ i_G.$$
Nondegeneracy of $\Phi$ is needed to have the compositions
$\Phi\circ i_A$ and $\Phi\circ i_G$ well defined.
In the definition of $\|\cdot\|_{\max}$ one could restrict
to non-degenerate or even (topologically)
irreducible representations of $(A,G,\alpha)$ on Hilbert space. However, it is extremely 
useful to consider more general covariant homomorphisms 
into multiplier algebras.
\\
{\bf (2)} The above described
correspondence between non-degenerate representations
of $(A,G,\alpha)$ and $A\rtimes G$ induces a bijection
between the set $(A,G,\alpha)\dach$ of unitary equivalence classes of irreducible 
covariant Hilbert-space representations
of $(A,G,\alpha)$ and $(A\rtimes G)\dach$. 
We topologize $(A,G,\alpha)\dach$ such that this bijection becomes a 
homeomorphism.
\\
{\bf (3)} The reduced crossed product $A\rtimes_{r}G$ does not
enjoy the above described universal properties, and therefore it
is often more difficult to handle. However, it follows from (\ref{eq-regular}) that
whenever $\rho:A\to M(D)$
is a $*$-homomorphism, then 
$\Ind\rho$ factors through a representation of $A\rtimes_rG$ to $M(D\otimes \K(L^2(G)))$
which is faithful iff $\rho$ is faithful.
In particular, if $\rho:A\to \B(H)$ is a faithful representation 
of $A$, then $\Ind\rho$ is a faithful representation
of $A\rtimes_rG$ into $\B(H\otimes L^2(G))$.
\\
{\bf (4)} By construction, the regular representation 
$\Lambda_A^G:C_c(G,A)\to A\rtimes_rG\subseteq M(A\otimes \K(L^2(G)))$ is the integrated
form of the covariant homomorphism $(i_{A,r}, i_{G,r})$ of $(A,G,\alpha)$ into $ M(A\otimes \K(L^2(G)))$ with
$$i_{A,r}=(\id_A\otimes M)\circ \tilde\alpha\quad \text{and}\quad i_{G,r}=1_A\otimes\lambda_G.$$
Since both, $\tilde{\alpha}:A\to M(A\otimes C_0(G))$ and $\id_A\otimes M: A\otimes C_0(G)\to M(A\otimes \K(L^2(G)))$ are faithful, it follows that $i_{A,r}$ is faithful, too.
Since $i_{A,r}=\Lambda_A^G\circ i_A$, where $i_A$ denotes the embedding 
of $A$ into $M(A\rtimes G)$, we see that  $i_A$ is injective, too.
\\
{\bf (5)} If $G$ is discrete, then $A$ embeds into $A\rtimes_{(r)}G$ via 
$a\mapsto \delta_e\otimes a\in C_c(G,A)\subseteq A\rtimes_{(r)}G$.
If, in addition, $A$ is unital, then $G$ also embeds into $A\rtimes_{(r)}G$
via $g\mapsto \delta_g\otimes 1$. If we identify $a\in A$ and $g\in G$ with their images
in $A\rtimes_{(r)}G$, we obtain the relations $ga=\alpha_g(a)g$ for all $a\in A$ 
and $g\in G$. The full crossed product is then the universal C*-algebra generated 
by $A$ and $G$ (viewed as a group of unitaries) subject to the relation $ga=\alpha_g(a)g$.
\\
{\bf (6)} In case $A=\CC$, the maximal crossed product $C^*(G):=\CC\rtimes G$ 
is called the {\em full  group
$C^*$-algebra} of $G$ (note that $\CC$ has only the trivial $*$-automorphism). 
The universal
properties of $C^*(G)$
translate into a one-to-one correspondence between
the  unitary representations of
$G$ and the non-degenerate $*$-representations of $C^*(G)$ which induces 
a bijection between the set $\widehat{G}$ of equivalence classes of irreducible unitary
Hilbert-space representations of $G$ and $\widehat{C^*(G)}$.
Again, we topologize $\widehat{G}$ so that this bijection becomes a homeomorphism.
\newline
The {\em reduced group $C^*$-algebra} $C_r^*(G):=\CC\rtimes_rG$
is realized as the closure $\overline{\lambda\big(C_c(G)\big)}\subseteq \B(L^2(G))$,
where $\lambda$ denotes the regular representation of $G$. 
\\
{\bf (7)}
If $G$ is compact, then every irreducible representation of $G$ is finite dimensional
and the Jacobson topology on $\widehat{G}=\widehat{C^*(G)}$ is the discrete
topology. Moreover,  it follows from the Peter-Weyl theorem (e.g. see \cites{DE, Fol}) that 
$C^*(G)$ and $C_r^*(G)$ are isomorphic to the 
$C^*$-direct sum $\bigoplus_{U\in \widehat{G}}M_{\dim{U}}(\CC)$.
In particular, we have $C^*(G)=C_r^*(G)$ if $G$ is compact.
\\
{\bf (8)} The convolution algebra $C_c(G)$, and hence also its completion $C^*(G)$,
is commutative if and only if $G$ is abelian.
In that case $\widehat{G}$ coincides with
the set of continuous homomorphisms from
$G$ to the circle group $\TT$, called {\em characters of $G$},
equipped with the compact-open topology.
The Gelfand-Naimark theorem for commutative $C^*$-algebras then implies that
 $C^*(G)\cong C_0(\widehat{G})$ (which also coincides with $C_r^*(G)$ in this case).
 Note that $\widehat{G}$, equipped with the pointwise multiplication of characters, is again a locally compact abelian group and the Pontrjagin duality theorem asserts that 
 $\widehat{\widehat{G}}$ is isomorphic to $G$
via $g\mapsto \widehat{g}\in \widehat{\widehat{G}}$ defined by 
$\widehat{g}(\chi)=\chi(g)$. 
Notice that the Gelfand isomorphism $C^*(G)\cong C_0(\widehat{G})$  extends 
the {\em Fourier transform} 
$$\mathcal F: C_c(G)\to C_0(\widehat{G}); \mathcal F(f)(\chi)=\chi(f)=\int_G f(x)\chi(x)\,dx.$$
For the circle group $\TT$ we have $\ZZ\cong \widehat{\TT}$ via $n\mapsto \chi_n$ with
$\chi_n(z)=z^n$, and one checks that the above Fourier transform coincides with 
the classical Fourier transform on $C(\TT)$. Similarly, if $G=\RR$, then 
$\RR\cong \widehat{\RR}$ via $s\mapsto \chi_s$ with $\chi_s(t)=e^{2\pi i st}$
and we recover the classical Fourier transform on $\RR$. We refer to \cite{DE}*{Chapter 3} for a detailed 
treatment of Pontrjagin duality and its connection to the Gelfand isomorphism\end{remark}

\begin{example}[Transformation group algebras]\label{ex-crossed}
      If $(X,G)$ is a topological dynamical system, then  we
      can form the crossed products $C_0(X)\rtimes G$ and 
      $C_0(X)\rtimes_{r}G$ with respect to the corresponding action
      of $G$ on $C_0(X)$. These algebras are often called the
      (full and reduced) {\em transformation group algebras}
      of the dynamical system $(X,G)$. Many important $C^*$-algebras
      are of this type. For instance if $X=\TT$ is the circle group
      and $\ZZ$ acts on $\TT$ via $n\cdot z=e^{i 2\pi \theta n}z$,
      $\theta\in [0,1]$,
      then $A_{\theta}=C(\TT)\rtimes\ZZ$ is the (rational or irrational)
      rotation algebra corresponding to $\theta$ (compare with \S \ref{subseccstar} above).
      Indeed, since $\ZZ$ is discrete and $C(\TT)$ is unital, we have canonical embeddings 
      of $C(\TT)$ and $\ZZ$ into $C(\TT)\rtimes \ZZ$. If we denote by $v$ the image
      of $\id_\TT\in C(\TT)$ and by $u$ the image of $1\in \ZZ$ under these embeddings,
      then the relations given in part (5) of the above remark show that $u,v$ are unitaries
      which satisfy the basic commutation relation  $uv=e^{2\pi i \theta }vu$.
      It is this realization as a crossed product of $A_{\theta}$ which motivates the
      notion ``rotation algebra''.

      There is quite some interesting and deep work on crossed products by
      actions of $\ZZ$ (or $\ZZ^d$) on compact spaces, which we cannot cover 
      in this article. We refer the interested reader to the article \cite{GPS} for a survey
      and for further references to this work.
  \end{example}

\begin{example}[Decomposition action]\label{ex-deco}
 Assume that $G=N\rtimes H$ is the semi-direct product ot two
      locally compact groups. If $A$ is a $G$-algebra, then
$H$ acts canonically on $A\rtimes N$ (resp. $A\rtimes_rN$)
via the extension of the action $\gamma$
of $H$ on $C_c(N,A)$ given by
$$\big(\gamma_h(f)\big)(n)=
      \delta(h)\alpha_h\big(f(h^{-1}\cdot n)\big),$$
  where $\delta:H\to \RR^+$ is determined by the equation
  $\int_N f(h\cdot n)\, dn=\delta(h)\int_N f(n)\,dn$ for all $f\in C_c(N)$.
      The  inclusion $C_c(N,A)\subseteq
      A\rtimes_{(r)}N$ determines an inclusion
      $C_c(N\times H,A)\subseteq C_c(H,A\rtimes_{(r)}N)$ which extends
      to isomorphisms $A\rtimes (N\rtimes H)\cong (A\rtimes N)\rtimes H$
and $A\rtimes_r(N\rtimes H)\cong (A\rtimes_rN)\rtimes_rH$.
In particular, if $A=\CC$, we obtain canonical isomorphisms
$C^*(N\rtimes H)\cong C^*(N)\rtimes H$ and $C_r^*(N\rtimes H)
\cong C_r^*(N)\rtimes_rH$.
\end{example}

We shall later extend the notion of crossed products to allow also the
decompostion of crossed products by group extensions which are
not topologically split.

\begin{remark}\label{rem-approx} When working with crossed products, it is often
useful to use the following concrete realization of an approximate unit
in $A\rtimes G$ (resp. $A\rtimes_rG$) in terms of a given 
approximate unit $(a_i)_{i\in I}$ in $A$: Let $\mathcal U$ be any
neighborhood basis of the identity $e$ in $G$, and for each 
$U\in \U$ let $\varphi_U\in C_c(G)^+$ with $\supp \varphi_U\subseteq U$, $\varphi_U(s)=\varphi_U(s^{-1})$ for all $s\in G$,
and such that $\int_G\varphi_U(t)\,dt=1$.
Let $\Lambda=I\times\U$ with $(i_1, U_1)\geq (i_2, U_2)$ if 
$i_1\geq i_2$ and $U_1\subseteq U_2$. Then a straightforward computation 
in the dense subalgebra $C_c(G,A)$ 
shows that $(\varphi_U\otimes a_i)_{(i,U)\in \Lambda}$ is an approximate
unit of $A\rtimes G$ (resp. $A\rtimes_rG$), where we write
$\varphi\otimes a$ for the function $(t\mapsto \varphi(t)a)\in C_c(G,A)$ if 
$\varphi\in C_c(G)$ and $a\in A$.
\end{remark}

\section{Crossed products versus tensor products}\label{sec-cvt}
The following lemma indicates the conceptual similarity
of full crossed products with maximal tensor products
and of reduced crossed products with minimal tensor products of
$C^*$-algebras.

\begin{lemma}
\label{lem-tensor}
   Let $(A,G,\alpha)$ be a $C^*$-dynamical system  and let
$B$ be a $C^*$-algebra. Let $\id\otimes_{\max}\alpha:G\to
\Aut(B\otimes_{\max}A)$
be the diagonal action of $G$ on $B\otimes_{\max}A$ (i.e., $G$ acts
trivially on $B$), and let
$\id\otimes \alpha:G\to \Aut(B\otimes A)$ denote the
diagonal action on $B\otimes A$.
Then the obvious map $B\odot C_c(G,A)\to C_c(G,B\odot A)$ induces
isomorphisms
$$B\otimes_{\max}(A\rtimes_{\alpha}G)\cong
(B\otimes_{\max}A)\rtimes_{\id\otimes\alpha}G
\quad\text{and}
\quad
B\otimes(A\rtimes_{\alpha,r}G)\cong
(B\otimes A)\rtimes_{\id\otimes\alpha,r}G.
$$
\end{lemma}

\begin{proof}[Sketch of proof] For the full crossed products
check that both sides have the same non-degenerate representations
and use the universal properties of the
full crossed products and the maximal tensor product.
For the reduced crossed products observe that the
map $B\odot C_c(G,A)\to C_c(G, B\odot A)$ identifies
$\id_B\otimes\Lambda_A^G$ with
$\Lambda_{B\otimes A}^G$.
\end{proof}

\begin{remark}\label{rem-fullgroup}
As a special case of the above lemma (with $A=\CC$) we see in
particular that
$$B\rtimes_{\id}G\cong B\otimes_{\max}C^*(G)\quad\text{and}
\quad B\rtimes_{\id,r}G\cong B\otimes C_r^*(G).$$
\end{remark}

We now want to study an important condition on $G$ which
implies that full and reduced crossed products by $G$
always coincide.

\begin{definition}\label{def-amenable}
      Let $1_G:G\to\{1\}\subseteq\CC$ denote the trivial representation of $G$.
      Then $G$ is called {\em amenable} if $\ker 1_G\supseteq
      \ker\lambda$ in $C^*(G)$, i.e., if the integrated form
      of $1_G$ factors through a homomomorphism $1_G^r:C_r^*(G)\to
      \CC$.{\ }\footnote{In particular, it follows that
      $1_G^r(\lambda_s)=1_G(s)=1$ for all $s\in G$!}
\end{definition}

\begin{remark} The above definition is not the standard
definition of amenability of groups, but it is
  one of the many equivalent formulations
for amenability (e.g. see \cites{Dix, Pat}),
and it is best suited for our purposes.
It is not hard to check
(even using the above $C^*$-theoretic definition)
that abelian groups and compact groups are amenable.
Moreover, extensions, quotients, and closed subgroups of amenable
groups are again amenable. In particular, all solvable groups
are amenable.

On the other side, one can show that the non-abelian free group
$F_2$ on two generators, and hence any group which contains
$F_2$ as a closed subgroup, is not amenable. This shows that 
non-compact semi-simple Lie groups are never amenable.
For extensive studies of amenability of groups (and groupoids) we refer
to \cites{Pat, ADR}.
\end{remark}

If $(\pi,U)$ is a covariant representation of $(A,G,\alpha)$
on some Hilbert space $H$,
then the covariant representation $(\pi\otimes 1, U\otimes \lambda)$
of $(A,G,\alpha)$ on $H\otimes L^2(G)\cong L^2(G,H)$ is unitarily
equivalent to $\Ind\pi$ via the unitary
$W\in U(L^2(G,H))$
defined by $(W\xi)(s)=U_s\xi(s)$ (this simple fact is known as {\em Fell's trick}).
Thus, if $\pi$ is faithful on $A$, then
$(\pi\otimes 1)\times (U\otimes\lambda)$
factors through a faithful representation of $A\rtimes_rG$.
As an important application we get

\begin{proposition}
\label{prop-amenable}
If $G$ is amenable, then
$\Lambda_A^G:A\rtimes_{\alpha}G\to A\rtimes_{\alpha,r}G$
is an isomorphism.
\end{proposition}
\begin{proof}
Choose any faithful representation $\pi\times U$
of $A\rtimes_{\alpha}G$ on some Hilbert space $H$.
Regarding $(\pi\otimes 1, U\otimes\lambda)$ as a representation
of $(A,G,\alpha)$ into $M(\K(H)\otimes C_r^*(G))$, we obtain
the equation
\[
(\id\otimes 1_G^r)\circ \bigl((\pi\otimes 1)\times(U\otimes \lambda)\bigr)
=\pi\times U.
\]
Since $\pi$ is faithful, it follows that
$$\ker\Lambda_A^G=\ker(\Ind\pi)=\ker\bigl((\pi\otimes 1)\times(U\otimes
\lambda)\bigr)
\subseteq\ker(\pi\times U)=\{0\}.$$
\end{proof}

The special case $A=\CC$ gives

\begin{corollary}\label{cor-amenable}
      $G$ is amenable if and only if $\lambda:C^*(G)\to C_r^*(G)$ is an
      isomorphism.
\end{corollary}

A combination of Lemma \ref{lem-tensor} with
Proposition \ref{prop-amenable} gives the following important result:

\begin{theorem}\label{thm-nuclear}
Let $A$ be a nuclear $G$-algebra with $G$ amenable.
The $A\rtimes_{\alpha}G$ is nuclear.
\end{theorem}
\begin{proof} Using Lemma \ref{lem-tensor} and
Proposition \ref{prop-amenable} we get
\begin{align*}
B\otimes_{\max}(A\rtimes_{\alpha}G)
&\cong
(B\otimes_{\max}A)\times_{\id\otimes\alpha}G
\cong
(B\otimes A)\times_{\id\otimes\alpha}G
\\
&\cong (B\otimes A)\times_{\id\otimes\alpha,r}G
\cong B\otimes (A\rtimes_{\alpha,r}G)
\cong B\otimes (A\rtimes_{\alpha}G).
\end{align*}
\end{proof}

If $(A,G,\alpha)$ and $(B,G,\beta)$ are two systems, then a
$G$-equivariant homomorphism
$\phi:A\to M(B)$ \footnote{where we uniquely extend
$\beta$ to an action of $M(B)$, which may fail to be strongly continuous}
induces a $*$-homomorphism
$$
\phi \rtimes G:=(i_B\circ \phi)\times i_G :A\rtimes_{\alpha}G\to
M(B\times_{\beta}G)
$$
where $(i_B, i_G)$ denote the canonical embeddings of $(B,G)$ into
$M(B\rtimes_{\beta}G)$, and a similar $*$-homomorphism
$$
\phi\rtimes_rG:= \Ind\phi:A\rtimes_{\alpha,r}G\to M(B\rtimes_{\beta,r}G)
\subseteq M\big(B\otimes\K(L^2(G))\big).
$$
Both maps are given on the level of functions by
$$\phi\rtimes_{(r)}G(f)(s)=\phi(f(s)),\;\; f\in C_c(G,A).$$
If $\phi(A)\subseteq B$, then
$\phi\rtimes G(A\rtimes_{\alpha}G)\subseteq B\rtimes_{\beta}G$ and
similarly for the reduced crossed products.
Moreover, $\phi\rtimes_rG=\Ind\phi$ is faithful if and only if $\phi$ is --- a
result which does {\bf not} hold in general for $\phi\rtimes G$!

On the other hand, the following proposition shows
that taking full crossed products gives an exact functor
between the category of $G$-$C^*$-algebras and the category
of $C^*$-algebras, which is not always true for the
reduced crossed-product functor!

\begin{proposition}\label{prop-exactfull}
      Assume that $\alpha:G\to \Aut(A)$ is an action and $I$ is a
      $G$-invariant closed ideal in $A$. Let
      $j:I\to A$ denote the inclusion and let $q:A\to A/I$
      denote the quotient map. Then the sequence
      $$
      \begin{CD}
	0\to I\rtimes_{\alpha}G @>{j\rtimes G}>> A\rtimes_{\alpha}G
      @>{q\rtimes G}>>
      (A/I)\rtimes_{\alpha}G\to 0
      \end{CD}
      $$
      is exact.
\end{proposition}
\begin{proof}
      If $(\pi,U)$ is a non-degenerate representation of
      $(I,G,\alpha)$ into $M(D)$, then $(\pi, U)$ has a canonical
      extension to a covariant homomorphism of $(A,G,\alpha)$
      by defining $\pi(a)(\pi(b)d)=\pi(ab)d$ for $a\in A, b\in I$ and
      $d\in D$. By the definition of $\|\cdot\|_{\max}$,
 this  implies that the inclusion $I\rtimes_{\alpha}G\to
      A\rtimes_{\alpha}G$ is isometric.

      Assume now that $p:A\rtimes_{\alpha}G\to
      (A\rtimes_{\alpha}G)/(I\rtimes_{\alpha}G)$
      is the quotient map. Then $p=\rho\times V$ for some
      covariant homomorphism $(\rho,V)$ of $(A,G,\alpha)$ into
      $M\big((A\rtimes G)/(I\rtimes G)\big)$. 
      Let $i_A:A\to M(A\rtimes G)$ denote the embedding. Then we have
      $i_A(I)C_c(G,A)= C_c(G,I)\subseteq I\rtimes G$ from which it follows
      that  
      $$\rho(I)(\rho\times V(C_c(G,A)))= \rho\times V\big(i_A(I)(A\rtimes G)\big)
      \subseteq \rho\times V(I\rtimes G)=\{0\}.$$
      Since $\rho\times V(C_c(G,A))$ is dense in $A/I\rtimes G$, it follows that 
      $\rho(I)=\{0\}$.
 Thus $\rho$ factors through a representation of $A/I$ and
      $p=\rho\times V$ factors through $A/I\rtimes_{\alpha}G$.
      This shows that the crossed product sequence is exact in the
      middle term. Since $C_c(G,A)$ clearly maps onto a dense
      subset in $A/I\rtimes_{\alpha}G$, $q\rtimes G$ is surjective
      and the result follows.
\end{proof}

For quite some time it was an open question whether the analogue
of Proposition \ref{prop-exactfull} also holds
for the reduced crossed products. This problem lead to

\begin{definition}[Kirchberg -- S. Wassermann]\label{def-exact}
A locally compact group $G$ is called {\em $C^*$-exact} (or simply exact)
if for any
system $(A,G,\alpha)$ and any $G$-invariant ideal $I\subseteq A$
the sequence
$$
\begin{CD}
0\to I\rtimes_{\alpha,r}G @>j\rtimes_r G>> A\rtimes_{\alpha,r}G
@> {q\rtimes_r G}>>
A/I\rtimes_{\alpha,r}G\to 0
\end{CD}
$$
is exact.
\end{definition}

Let us remark that only exactness in the middle term is the problem,
since $q\rtimes_rG$ is clearly surjective, and $j\rtimes_rG=\Ind j$
is injective since $j$ is. We shall later report on
Kirchberg's and S. Wassermann's permanence results on
exact groups, which imply that the class of exact groups
is indeed very large. However, a construction based on ideas of  Gromov 
(see \cites{Gr1, Ghys, Osajda}) implies that there do exist
  finitely generated discrete groups which are not exact!

\section{The correspondence categories}\label{sec-mor}

In this section we want to give some theoretical
background for  the discussion of  imprimitivity 
theorems for crossed products and for the theory of induced 
representations on one side, and for the construction of Kasparov's 
bivariant $K$-theory groups on the other side.
The basic notion for this is the notion of the correspondence category
in which the objects are $C^*$-algebras and  the morphisms are unitary
equivalence classes of Hilbert bimodules. Having this at hand, 
the theory of induced representations
will reduce to taking  compositions of morphisms in the correspondence category.
All this is based on the fundamental idea of Rieffel (see \cite{Rie1}) who first 
made a systematic approach to the theory of induced representations of $C^*$-algebras
in terms of (pre-)  Hilbert modules, and who showed how the theory of induced
group representations can be seen as part of this more general theory.
However, it seems that a systematic categorical treatment of this theory was 
first given in \cite{EKQR1} and, in parallel work,  by Landsman in \cite{Lan}.
The standard reference for Hilbert modules is \cite{Lance}.

\subsection{Hilbert modules} \label{subsec-Hilbert}
If $B$ is a $C^*$-algebra, then a 
(right) Hilbert $B$-module is a complex Banach 
space $E$ equipped with a right $B$-module structure and
a positive definit $B$-valued inner 
product (with respect to positivity in $B$)
$\lk\cdot, \cdot\rk_B:E\times E\to B$, which is linear in the 
second 
and antilinear in the first variable and satisfies
$$(\lk 
\xi, \eta\rk_B)^*=\lk \eta,\xi\rk_B,\quad \lk \xi, \eta\rk_Bb=\lk 
\xi,\eta\cdot b\rk_B,
\quad\text{and}\quad
\|\xi\|^2=\|\lk \xi,\xi\rk_B\|$$ 
for all $\xi, \eta\in E$ and $b\in B$.
With the obvious modifications we can also define
{\em left}-Hilbert $B$-modules. The Hilbert $\CC$-modules 
are precisely the Hilbert spaces. Moreover, every $C^*$-algebra 
$B$ becomes a Hilbert $B$-module by defining $\lk b,c\rk_B:= b^*c$.
We say that $E$ is a {\em full} Hilbert $B$-module, if 
$$B=\lk E,E\rk_B:=\csp\{\lk \xi,\eta\rk_B:\xi,\eta\in E\}.$$
In general $\lk E,E\rk_B$ is a closed two-sided ideal of $B$.

If $E$ and $F$ are Hilbert $B$-modules, then a linear map $T:E\to F$ is called {\em adjointable}
if there exists a map $T^*:F\to E$ such that $\langle T\xi,\eta\rangle_B=\langle \xi,T^*\eta\rangle_B$ 
for all $\xi\in E, \eta\in F$.\footnote{Note that, different from the operators on Hilbert space, a bounded $B$-linear 
operator $T:E\to F$ is {\bf not} automatically adjointable}
Every adjointable operator from $E$ to $F$ is automatically 
bounded and $B$-linear. We write $\mathcal L_B(E,F)$  for the set of adjointable operators from $E$ to $F$.
Then
$$\mathcal L_B(E):=\mathcal L_B(E,E)$$
becomes a $C^*$-algebra with respect to the usual operator norm. Every pair $\xi,\eta$ with $\xi\in F$, $\eta\in E$ determines 
an element $\Theta_{\xi,\eta}\in \mathcal L_B(E,F)$ given by
\begin{equation}\label{eq-theta} \Theta_{\xi,\eta}(\zeta)=\xi\langle \eta,\zeta\rangle_B\end{equation}
with adjoint $\Theta_{\xi,\eta}^*=\Theta_{\eta,\xi}$. The closed linear span of all such operators  forms 
the set of compact operators $\mathcal K_B(E,F)$ in $\mathcal L_B(E,F)$. If $E=F$, then $\mathcal K_B(E):=\mathcal K_B(E,E)$
is a closed ideal in $\mathcal L_B(E)$. Note that there is an obvious $*$-isomorphism between the multiplier algebra $M(\mathcal K_B(E))$ 
and $\mathcal L_B(E)$, which is given by extending the action of $\mathcal K_B(E)$ on $E$ to all of $M(\mathcal K_B(E))$ in the canonical way.

\begin{example}\label{ex-Hilbert}
{\bf (1)} If $B=\CC$ and  $H$ is a Hilbert space, then $\L_\CC(H)=\B(H)$ and  $\K_\CC(H)=\K(H)$.\\
{\bf (2)} If  a  $C^*$-algebra $B$ is viewed as a Hilbert $B$-module
with respect to the inner product $\lk b,c\rk_B=b^*c$ and the obvious right module operation
then $ \K_B(B)=B$, where we let $B$ act on itself via left multiplication,
and  we have $\L_B(B)=M(B)$. 

It is important to notice that, in case $B\neq\CC$,
the notion of compact operators as given above does {\bf not} coincide
with the standard notion of compact operators on a Banach space (i.e., that the image of the unit
ball has compact closure). For example, if $B$ is unital, then $\L_B(B)=\K_B(B)=B$
and we see that the identity operator on $B$ is a compact operator 
in the sense of the above definition. But if $B$ is not finite dimensional, the 
identity operator is not a compact operator in the usual sense of Banach-space operators.
\end{example}

There is a one-to-one correspondence between right and left Hibert $B$-modules
given by the operation $E\mapsto E^*:=\{\xi^*:\xi\in E\}$, 
with left action of $B$ on $E^*$ given by $b\cdot\xi^*:=(\xi\cdot b^*)^*$ and 
with inner product $_B\lk\xi^*,\eta^*\rk:=\lk\xi,\eta\rk_B$ (notice that the inner product of a  
left Hilbert $B$-module is linear in the first and antilinear in the second variable).
We call $E^*$ the {\em adjoint module} of $E$. Of course, if $F$ is a left 
Hilbert $B$-module, a similar construction yields an adjoint $F^*$ -- a right Hilbert $B$-module.
Clearly, the notions of adjointable and compact operators also have their left analogues
(thought of as acting on the right),
and we have $\L_B(E) = \L_B(E^*)$ (resp. $\K_B(E) = \K_B(E^*)$)
via $\xi^*T:=(T^*\xi)^*$.

There are several important
 operations on Hilbert modules (like taking the direct sum 
$E_1\oplus E_2$ of two Hilbert $B$-modules $E_1$ and $E_2$ in the obvious way).
 But for our considerations the construction of the interior tensor products is most
important. For this assume that $E$ is a (right) Hilbert $A$-module, $F$ is a
(right) Hilbert $B$-module, and $\Psi:A\to \L_B(F)$ is a $*$-homomorphism.
Then the {\em interior tensor product}
$E\otimes_AF$ is defined as the 
Hausdorff completion of $E\odot F$ with respect to the 
$B$-valued inner 
product
$$\lk \xi\otimes \eta, \xi'\otimes \eta'\rk_B=\lk  \eta, \Psi(\lk 
\xi,\xi'\rk_A)\cdot\eta'\rk_B,$$
where $\xi,\xi'\in E$ and 
$\eta,\eta'\in F$. 
With this inner product, $E\otimes_AF$ becomes a 
Hilbert $B$-module. Moreover, if $C$ is a third $C^*$-algebra 
and if $\Phi:C\to\L_A(E)$ is a $*$-representation of $C$ on $\L_A(E)$,
then $\Phi\otimes  1: C\to \L_B(E\otimes_AF)$ with 
$\Phi\otimes 1(c)(\xi\otimes \eta)=\Phi(c)\xi\otimes\eta$ becomes a $*$-representation 
of $C$ on $E\otimes_AF$ (we refer to \cites{Lance, RW} for more details). 
The construction of this representation is absolutely crucial
in what follows below.

\subsection{Morita equivalences}
The notion of Morita equivalent $C^*$-algebras, which goes back to Rieffel \cite{Rie1}
 is one of the most important tools in the study of crossed products.

\begin{definition}[Rieffel]\label{def-imp}
Let $A$ and $B$ be 
$C^*$-algebras. An {\em $A$-$B$ imprimitivity bimodule}\footnote{often called 
an $A$-$B$ {\em equivalence} bimodule in the literature}
$X$ is a 
Banach space $X$ which carries the structure of both, a 
right Hilbert 
$B$-module and  a  left Hilbert $A$-module with commuting actions of $A$ and $B$
 such
that 
\begin{enumerate}
\item $_A\lk X,X\rk=A$ and 
$\lk X,X \rk_B=B$  (i.e., both inner products on $X$ are full);
 \item $_A\lk 
\xi,\eta\rk \cdot\zeta=\xi\cdot\lk\eta,\zeta\rk_B$ for all 
$\xi,\eta,\zeta\in X$.
\end{enumerate}
$A$ and $B$ are called {\em Morita equivalent} 
if such 
$A$-$B$ bimodule 
$X$ exists.
\end{definition}

\begin{remark}\label{rem-Mor}
{\bf (1)} It follows from the above definition together with
(\ref{eq-theta}) that,
if $X$ is an 
$A$-$B$ imprimitivity bimodule, then $A$ canonically identifies with 
$\K_B(X)$ and $B$ 
canonically identifies with $\K_A(X)$.
Conversely, if $E$ is any Hilbert $B$-module,
then $_{\K(E)}\lk \xi,\eta\rk:=\Theta_{\xi,\eta}$
(see (\ref{eq-theta})) defines a full $\K_B(E)$-valued inner product on 
$E$, and $E$ becomes a $\K_B(E)$-$\lk E,E\rk_B$ imprimitivity bimodule.
In particular, if $E$ is a full Hilbert $B$-module (i.e., $\lk E,E\rk_B=B$), then $B$ is 
Morita equivalent to $\K_B(E)$. 
\\
{\bf (2)} As a very special case of (1) we see that $\CC$ is Morita equivalent to $\K(H)$
for every Hilbert space $H$.
\\
{\bf (3)} 
It is easily checked that Morita equivalence 
is an equivalence relation:
If $A$ is any $C^*$-algebra, then 
$A$ becomes an $A$-$A$ imprimitivity bimodule 
with respect to 
$_A\lk a,b\rk= ab^*$ and $ \lk a,b\rk_A=a^*b$ for $a,b\in 
A$.
If $X$ is an $A$-$B$ imprimitivity bimodule and $Y$ is a $B$-$C$ 
imprimitivity
bimodule, then $X\otimes_BY$ is an $A$-$C$ imprimitivity 
bimodule. Finally, if $X$ is an $A$-$B$-imprimitivity  bimodule, then
the adjoint module $X^*$
is a $B$-$A$ imprimitivity bimodule.
 \\
{\bf (4)} 
Recall that a $C^*$-algebra $A$ is a {\em full corner} of the 
$C^*$-algebra $C$,
if there exists a full projection  $p\in M(C)$ 
(i.e., $\overline{CpC}=C$) 
such that  $A=pCp$. Then $pC$ equipped 
with the canonical  
 inner products and actions coming from 
multiplication and involution on $C$
becomes an $A$-$C$ imprimitivity   
bimodule. Thus, if $A$ and $B$ can be 
represented as full corners of 
a $C^*$-algebra $C$, they are Morita equivalent.
Conversely,  let 
$X$ be an 
$A$-$B$ imprimitivity bimodule. Let 
$L(X)=\left(\begin{smallmatrix} A&X\\ X^*&B\end{smallmatrix}
\right)$ with 
multiplication and involution defined by
\begin{align*}
\left(\begin{matrix} a_1& 
\xi_1\\ \eta_1^*&b_1\end{matrix}\right)
\left(\begin{matrix} a_2& 
\xi_2\\ \eta_2^*&b_2\end{matrix}\right)
&=\left(\begin{matrix} a_1a_2+ 
{_A\lk \xi_1,\eta_2\rk}& a_1\cdot \xi_2+\xi_1\cdot b_2\\ 
\eta_1^*\cdot a_2+b_1\cdot \eta_2^*& \lk \eta_1, 
\xi_2\rk_B+b_1b_2\end{matrix}\right)
\;\text{and}\\
\left(\begin{matrix} a& \xi \\ 
\eta^*&b\end{matrix}\right)^*&=
\left(\begin{matrix} a^*& \eta\\ 
\xi^*&b^*\end{matrix}\right).
\end{align*}
Then $L(X)$ has a canonical embedding as a 
closed subalgebra of
 the adjointable operators on the Hilbert 
$B$-module
$X\oplus B$ via
$$\left(\begin{matrix} a& \xi \\ 
\eta^*&b\end{matrix}\right)\left(\begin{matrix} \zeta\\ d\end{matrix}\right)=
\left(\begin{matrix}a \zeta + \xi d\\ \lk \eta,\zeta\rk_B+ bd\end{matrix}\right)$$
which makes $L(X)$ a $C^*$-algebra. 
If 
$p=\left(\begin{smallmatrix} 1& 0\\ 0&0\end{smallmatrix}\right)\in 
M(L(X))$, then
$p$ and $q:=1-p$ are full projections such that 
$A=pL(X)p$, $B=qL(X)q$ and 
$X=pL(X)q$. 
The algebra $L(X)$ is 
called the {\em linking algebra} of $X$. It often serves as
a valuable 
tool for the study of imprimitivity bimodules.
\\
{\bf (5)} It 
follows from {\bf (4)} that  $A$ is Morita equivalent to $A\otimes \K(H)$ 
for any 
Hilbert space $H$ (since $A$ is a full corner of $A\otimes \K(H)$). 
Indeed, a deep theorem of Brown, Green and Rieffel (see \cite{BGR})
shows 
that if $A$ and $B$ are $\sigma$-unital\footnote{A $C^*$-algebra is called {\em $\sigma$-unital},
if it has a countable approximate unit. In particular, all separable and all unital $C^*$-algebras are 
$\sigma$-unital}, then $A$ and $B$ are Morita equivalent
 if and only if they are {\em stably isomorphic},
i.e., there exists an isomorphism between $A\otimes \K(H)$ and $B\otimes \K(H)$
with $H=l^2(\NN)$. A similar result does not hold if the $\sigma$-unitality 
assumption is dropped (see \cite{BGR}).
\\
{\bf (6)} The above results indicate that many important
properties of $C^*$-algebras are preserved by Morita equivalences.
Indeed, among these properties are: nuclearity, exactness, simplicity,
the property of being a type I algebra (and many more). 
Moreover, Morita equivalent $C^*$-algebras have homeomorphic 
primitive ideal spaces and isomorphic $K$-groups. Most of  these 
properties will be discussed later in more detail
(e.g., see Propositions
\ref{prop-rief-cor}, \ref{prop-induce} and \ref{prop-typeI} below).
The $K$-theoretic implications are discussed in  \cite{KK}.
\end{remark}

A very important tool when working with imprimitivity bimodules is 
the Rieffel correspondence. To explain this correspondence suppose that 
$X$ is an $A$-$B$ imprimitivity bimodule and that
$I$ is a closed ideal of $B$. Then $X\cdot I$ is a closed $A$-$B$ submodule 
of $X$ and $\Ind^XI:={_A\lk X\cdot I, X\cdot I\rk}$ (taking the closed span)
is a closed ideal of $A$.  The following proposition implies that Morita equivalent 
$C^*$-algebras have equivalent ideal structures:

\begin{proposition}[Rieffel correspondence]\label{prop-rief-cor}
Assume notation as above. Then
\begin{enumerate}
\item The assignments $I\mapsto X\cdot I$,  $I\mapsto 
\Ind^XI$ and $I\mapsto J_I:=\left(\begin{smallmatrix} \Ind^XI & X\cdot I\\ I\cdot X^*& I
\end{smallmatrix}\right)$
provide  inclusion preserving bijective correspondences 
between
the closed two-sided ideals of $B$, the closed 
$A$-$B$-submodules of $X$, 
the closed two-sided ideals of $A$, and the closed two-sided ideals of the linking algebra $L(X)$,
respectively.
\item $X\cdot I$ is an 
$\Ind^XI$-$I$ imprimitivity 
bimodule and $X/(X\cdot I)$, equipped with the obvious inner products 
and bimodule actions, becomes an
$A/(\Ind^XI)$-$B/I$ imprimitivity bimodule.
Moreover, we have $J_I=L(X\cdot I)$ and $L(X)/J_I\cong L(X/X\cdot I)$.
\end{enumerate}
\end{proposition}

\begin{remark}\label{rem-imphom}
Assume that $X$ is an $A$-$B$ imprimitivity bimodule and $Y$ is a $C$-$D$ imprimitivity
bimodule. An {\em imprimitivity bimodule homomorphism} from $X$ to $Y$
is then a triple $(\phi_A, \phi_X,\phi_B)$ such that $\phi_A:A\to C$ and $\phi_B:B\to D$
are $*$-homomorphisms and $\phi_X:X\to Y$ is a linear map such that the triple
$(\phi_A, \phi_X,\phi_B)$ satisfies the obvious compatibility conditions with respect 
to the inner products and module actions on $X$ and $Y$ (e.g. 
$\lk\phi_X(\xi),\phi_X(\eta)\rk_D=\phi_B(\lk \xi,\eta\rk_B)$, $\phi_X(\xi b)=\phi_X(\xi)\phi_B(b)$,
etc.). 

If $(\phi_A, \phi_X,\phi_B)$ is such an imprimitivity bimodule homomorphism,
then one can check that $\ker\phi_A$, $\ker\phi_X$ and $\ker\phi_B$ 
all correspond to each other under the Rieffel correspondence for $X$ (e.g., see \cite{EKQR2}*{Chapter 1}).
\end{remark}

As a simple application of the Rieffel correspondence and the above remark 
we now show 

\begin{proposition}\label{prop-nuclear} 
Suppose that $A$ and $B$ are Morita equivalent $C^*$-algebras. Then $A$ is nuclear if and only 
if $B$ is nuclear.
\end{proposition}
\begin{proof}[Sketch of proof]
Let $X$ be an $A$-$B$ imprimitivity bimodule.
If $C$ is any other 
$C^*$-algebra, we can equip $X\odot C$ with $A\odot C$- and $B\odot C$-valued
inner products and an $A\odot C$-$B\odot C$ module structure in the obvious way.
Then one can check that $X\odot C$ completes to an $A\otimes_{\max} C$-$B\otimes_{\max}C$
imprimitivity bimodule $X\otimes_{\max} C$ as well as to an $A\otimes C$-$B\otimes C$
imprimitivity bimodule $X\otimes C$. The identity map on $X\odot C$ then extends
to a quotient map $X\otimes_{\max} C\to X\otimes C$ 
which together with the quotient maps $A\otimes_{\max}C\to A\otimes C$ 
and $B\otimes_{\max}C\to B\otimes C$ is an imprimitivity bimodule homomorphism.
But then it follows from the above remark and the Rieffel correspondence 
that injectivity of any one of these quotient maps implies injectivity of all three of them.
\end{proof}

\subsection{The correspondence categories.} \label{subsec-corr}
We now come to the definition of the correspondence categories.
Suppose that $A$ and $B$ are $C^*$-algebras. A {\em (right) Hilbert $A$-$B$ bimodule}
is a pair $(E,\Phi)$ in which $E$ is a Hilbert $B$-module and $\Phi:A\to \L_B(E)$
is a $*$-representation of $A$ on $E$.
We say that $(E,\Phi)$ is non-degenerate, if $\Phi(A)E=E$ (this is equivalent
to $\Phi:A\to M(\K_B(E))=\L_B(E)$ being non-degenerate in the usual sense).
Two Hilbert $A$-$B$ bimodules $(E_i,\Phi_i)$, $i=1,2$ are called {\em unitarily equivalent} if
there exists an isomorphism $U:E_1\to E_2$ preserving the $B$-valued inner products 
 such that $U\Phi_1(a)=\Phi_2(a)U$ for all $a\in A$.
 Note that for any Hilbert $A$-$B$ bimodule $(E,\Phi)$ the module $(\Phi(A)E,\Phi)$ 
 is a non-degenerate $A$-$B$ sub-bimodule of $(E,\Phi)$. Note that 
 $\Phi(A)E=\{\Phi(a)\xi: a\in A, \xi\in E\}$ equals $\cspn(\Phi(A)E)$ by  Cohen's 
 factorisation theorem.
 %We denote by $[E,\Phi]$ the unitary equivalence class of $(E,\Phi)$.

\begin{definition}[{cf. \cites{BEW, EKQR1, EKQR2,  Lan}}]\label{defn-Morcat}
The {\em correspondence category} (also called the {\em Morita category}) 
$\Corr$ is the category whose 
objects are $C^*$-algebras and where the morphisms  from $A$ to $B$
are given by equivalence
classes $[E,\Phi]$ of  Hilbert $A$-$B$ bimodules $(E,\Phi)$ under the equivalence 
relation
$$(E_1,\Phi_i)\sim (E_2,\Phi_2)\Leftrightarrow \Phi_1(A)E_1\cong \Phi_2(A)E_2,$$
where $\cong$ denotes unitary equivalence.
The identity morphism from $A$ to $A$ is represented by the trivial 
$A$-$A$ bimodule $(A,\id)$
and composition of 
two morphisms  $[E,\Phi]\in \Mor(A,B)$ and $[F,\Psi]\in \Mor(B,C)$
is given by taking the interior tensor product $[E\otimes_BF,\Phi\otimes 1]$.

The {\em compact correspondence category}  $\Corr_c$ is the subcategory of $\Corr$ in which we 
additionally require
$\Phi(A)\subseteq \K_B(E)$ for a morphism $[E,\Phi]\in \Mor_c(A,B)$.
\end{definition}

\begin{remark}\label{rem-set}
{\bf (1)} We should note that the correspondence category is not a category in the strong sense, since 
the morphisms $\Mor(A,B)$ from $A$ to $B$ do not form a set. This problem can be overcome
by restricting the size of the objects and the underlying modules for the morphisms by 
assuming that they contain dense subsets of a certain maximal cardinality.
But  for most practical  aspects this does not cause any problems.\\
{\bf (2)} We should also note that in most places of the literature (e.g., in \cite{EKQR2}) 
the correspondence category is defined 
as the category with objects the $C^*$-algebras and with morphism sets 
$\Mor(A,B)$ given by unitary equivalence classes of {\em non-degenerate} 
$A$-$B$ bimodules. But the correspondence category $\Corr$ defined above is equivalent 
to the one of \cite{EKQR2} where the equivalence is given by 
the identity map on objects and by assigning a morphism  $[E,\Phi]\in \Mor(A,B)$ 
to the unitary equivalence class $[\Phi(A)E,\Phi]$ in the morphism set 
as in \cite{EKQR2}. 
\\
{\bf (3)}  Note that every  $*$-homomorphism
$\Phi:A\to M(B)$ determines a morphism $[E,\Phi]\in \Mor(A,B)$ in $\Corr$ with $E=B$,
and $[E,\Phi]$ is a morphism in  $\Corr_c$ if and only if $\Phi(A)\subseteq B$.
\\
{\bf (4)} Taking direct sums of bimodules allows to define sums of morphisms 
in the correspondence categories (and hence a semi-group structure with 
neutral element given by the zero-module). It is easy to check that this operation is 
commutative and satisfies the distributive  law with respect to composition.
\end{remark}

If $X$ is an $A$-$B$ imprimitivity bimodule, then the adjoint module $X^*$ satisfies
$X\otimes_BX^*\cong A$ as $A$-$A$ bimodule (the isomorphism given on elementary tensors 
by ${x\otimes y^*}$ $\mapsto {_A\lk x,y\rk}$) and $X^*\otimes_AX\cong B$ as $B$-$B$ bimodule,
so $X^*$ is an inverse of $X$ in the correspondence categories. Indeed we have

\begin{proposition}[cf \cites{Lan, EKQR2}]\label{prop-isomor}
The isomorphisms in the categories $\Corr$ and $\Corr_c$ are precisely the Morita equivalences.
\end{proposition}

\subsection{The equivariant correspondence categories}\label{subsec-equivmor}
If $G$ is a locally compact group, then the {\em $G$-equivariant correspondence category}
$\Corr(G)$ is the category in which the objects are systems $(A,G,\alpha)$
and  morphisms from $(A,G,\alpha)$ to $(B,G,\beta)$ are the
equivalence classes (as in the non-equivariant case) of equivariant 
$A$-$B$ Hilbert bimodules $(E,\Phi,u)$, i.e., $E$ is equipped with a 
strongly continuous homomorphism
 $u:G\to \Aut(E)$ such that
\begin{equation}\label{eq-mor}
\begin{split}
 \lk u_s(\xi),u_s(\eta)&\rk_B=\beta_s(\lk \xi,\eta\rk_B), \quad 
 u_s(\xi\cdot b)= u_s(\xi)\beta_s(b)\\
 &\text{and}\quad
 u_s(\Phi(a)\xi)=\Phi(\alpha_s(a)) u_s(\xi).
 \end{split}
\end{equation}
Of course we require that a unitary equivalence $U:E_1\to E_2$ between two 
$G$-equivariant Hilbert bimodules also intertwines that actions of $G$ on $E_1, E_2$.
Again, composition of morphisms is given by taking interior tensor products
equipped with the diagonal actions, and the isomorphisms in this category 
are just the equivariant Morita equivalences.
%
%More generally, if $N$ is a fixed normal subgroup of $G$, one can consider the 
%{\em twisted equivariant correspondence category} $\Corr(G,N)$ in which the objects are 
%the twisted systems $(A,G,N,\alpha,\tau)$ and in which a morphism from 
%$(A,G,N,\alpha,\tau)$ to $(B,G,N,\beta,\sigma)$ is given by morphism 
%$[E,\Phi, u]$ from $(A,G,\alpha)$ to $(B,G,\beta)$ in $\Corr(G)$ 
%which preserves the twists in the sense that
%\begin{equation}\label{eq-mortwist}
%\Phi(\tau_n)\xi= u_n(\xi)\sigma_n\quad\text{for all $n\in N$.}
%\end{equation}

Notice that the  crossed product constructions
 $A\rtimes G$  and  $A\rtimes_rG$
extend to descent functors
$$\rtimes_{(r)}:\Corr(G)\to \Corr.$$
In particular,  Morita equivalent systems have 
Morita equivalent full (resp. reduced) crossed products.
If $[E,\phi, u]$ is a morphism from $(A,G,\alpha)$ to $(B,G,\beta)$,
then the crossed product $[E\rtimes_{(r)}G, \Phi\rtimes_{(r)}G]\in \Mor(A\rtimes_{(r)}G,B\rtimes_{(r)}G)$ is given as the completion of 
$C_c(G,E)$ with respect to the $B\rtimes_{(r)}G$-valued inner product
$$\lk\xi,\eta\rk_{B\rtimes_{(r)}G}(t)=\int_G\beta_{s^{-1}}(\lk \xi(s), \eta(st)) \rk_B)\,ds$$
(taking values in $C_c(G,B)\subseteq B\rtimes_{(r)}G$) and with left 
action of $C_c(G,A)\subseteq A\rtimes_{(r)}G$ on $E\rtimes_{(r)}G$ 
given by
$$\big(\Phi\rtimes_{(r)}G(f)\xi\big)(t)=\int_G \Phi(f(s)) u_s(\xi(s^{-1}t))\,ds.$$
%In the twisted case the crossed products $E\rtimes_{(r)}(G,N)$ can be defined as 
%the quotients $(E\rtimes G)/\big((E\rtimes G)\cdot I_{(r)}\big)$, with 
%$I_{(r)}:=\ker\big(B\rtimes G\to B\rtimes_{(r)}(G,N)\big)$.  Of course, we can also consider
%the compact equivariant correspondence categories $\Corr_c(G)$ and $\Corr_c(G,N)$ by 
%requiring $\Phi:A\to \L_B(E)$ taking image in $\K_B(E)$. The above defined descent 
%functors then restrict to functors between $\Corr_c(G)$ (resp. $\Corr_c(G,N)$) to $\Corr_c$.
%
%The correspondence categories (equivariant for actions and coactions of groups) have been studied 
%extensively in \cite{EKQR1, EKQR2}, but the existence of such categories has been 
%indicated by several authors before (e.g. see the introduction of \cite{EKQR2}).
%However, as far as I know, this is the first place where they are called correspondence categories.
The crossed product constructions for equivariant bimodules first appeared (to my knowledge) in 
Kasparov's famous Conspectus \cite{K2}, which circulated as a preprint 
from the early eighties. A more detailed study in case of imprimitivity bimodules
has been given in \cite{Combes}.  A very extensive study of the equivariant correspondence categories
for actions and coactions of groups together with their relations to duality theory
are given in \cite{EKQR2}.

\subsection{Induced representations and ideals}\label{subsec-ind-rep}
If $B$ is a $C^*$-algebra we denote by $\Rep(B)$ the 
collection of all unitary equivalence classes of non-degenerate 
$*$-representations of $B$ on Hilbert space. 
In terms of the correspondence category, $\Rep(B)$ coincides with the 
collection $\Mor(B,\CC)$ of morphisms from  $B$ to $\CC$
in $\Corr$ (every morphism can be represented by a non-degenerate
$*$-representation which is unique up to unitary equivalence). 
Thus, if $A$ is any other $C^*$-algebra and if $[E,\Phi]\in \Mor(A,B)$,
then composition with $[E,\Phi]$ determines a map
$$\Ind^{(E,\Phi)}:\Rep(B)\to \Rep(A); [H,\pi]\mapsto  [H,\pi]\circ [E,\Phi]=
[E\otimes_BH,\Phi\otimes 1].$$
If confusion seems unlikely, we will simply  write $\pi$ for the representation $(H,\pi)$
and for its class $[H,\pi]\in \Rep(A)$
and we write $\Ind^E\pi$ for the representation $\Phi\otimes 1$ of $A$ on
$\Ind^EH:=E\otimes_BH$.
We call $\Ind^E\pi$  the {\em
representation of $A$ induced from $\pi$ via $E$}.

Note that in the above, we did not require the action of $A$ on $E$ to be non-degenerate.
If it fails to be non-degenerate, the representation $\Phi\otimes 1$ of $A$ on $E\otimes_BH$ may also fail
to be non-degenerate. We then  pass to the restriction of $\Phi\otimes 1$ to
$\Phi\otimes 1(A)(E\otimes_BH)\subseteq E\otimes_BH$ to obtain a 
non-degenerate representative of $[E\otimes_BH,\Phi\otimes 1]\in \Mor(A,\CC)=\Rep(A)$.

\begin{remark}\label{rem-ind}
{\bf (1)} A special case of the above procedure is given in case when 
$\Phi:A\to M(B)$ is a non-degenerate $*$-homomorphism and  $[B,\Phi]\in \Mor(A,B)$ is the corresponding morphism in $\Corr$. 
Then the induction map $\Ind^B:\Rep(B)\to \Rep(A)$ coincides with 
the obvious map 
$$\Phi^*:\Rep(B)\to \Rep(A); \pi\mapsto \Phi^*(\pi):=\pi\circ \Phi.$$
{\bf (2) Induction in steps.} If 
$[H,\pi]\in \Rep(B)$, $[E,\Phi]\in \Mor(A,B)$ and $[F,\Psi]\in \Mor(D,A)$
for some $C^*$-algebra $D$, then it follows directly from the
associativity of composition in $\Corr$ that (up to equivalence)
$$\Ind^F(\Ind^E\pi)=\Ind^{F\otimes_AE}\pi.$$
{\bf (3)} If $X$ is an $A$-$B$ imprimitivity bimodule, then 
$\Ind^X:\Rep(B)\to \Rep(A)$ gets  inverted by $\Ind^{X^*}:\Rep(A)\to \Rep(B)$,
where $X^*$ denotes the adjoint of $X$ (i.e., the inverse of $[X]$ in $\Corr$). 
Since composition  of morphisms 
in $\Corr$ preserves direct sums, it follows from this 
that induction via $X$ maps 
irreducible representations of $B$ to irreducible representations of $A$ and hence
induces a bijection $\Ind^X:\widehat{B}\to \widehat{A}$ between the spectra.
\end{remark}

It is useful to consider a similar induction map on the set $\I(B)$
of closed two sided ideals of the $C^*$-algebra $B$. If $(E,\Phi)$ is any 
Hilbert  $A$-$B$ bimodule,   we define
\begin{equation}\label{eq-indideal}
\Ind^E:\I(B)\to \I(A);\; \Ind^EI:=\{a\in A: \lk \Phi(a)\xi,\eta\rk_B\in I\;\text{for all}\; \xi,\eta\in E\}.\;\footnote{If $X$ is an $A$-$B$ imprimitivity bimodule,
the induced ideal $\Ind^XI$ defined 
here  coincides with the induced ideal $\Ind^XI={_A\lk X\cdot I,X\cdot I\rk}$ 
of the Rieffel correspondence (see Proposition \ref{prop-rief-cor}).}
\end{equation}
It is clear that induction preserves inclusion of ideals
and with a little more work one can check that 
\begin{equation}\label{eq-kerindpi}
\Ind^E(\ker\pi)=\ker(\Ind^E\pi) \quad\text{for all $\pi\in \Rep(B)$.}
\end{equation}
Hence it follows from part (3) of Remark \ref{rem-ind} that,
if $X$ is an $A$-$B$ imprimitivity bimodule, then induction of ideals 
via $X$ restricts to give a bijection $\Ind^X:\Prim(B)\to \Prim(A)$
between the primitive ideal spaces of $B$ and $A$.
Since induction preserves inclusion 
of ideals, the next proposition follows directly from the description of the closure operations in 
$\widehat{A}$ and $\Prim(A)$ (see \S \ref{subsec-spectrum}).

%

%Recall that the closure operations in the Jacobson topologies
%of  $\widehat{B}$ and $\Prim(B)$ are given by
%$$\pi\in \bar{R}\Leftrightarrow \ker\pi\supseteq\cap\{\ker\rho:\rho\in R\}
%\quad\text{and}\quad P\in \bar{S}\Leftrightarrow P\supseteq\cap\{Q:Q\in S\},$$
%if $R\subseteq \widehat{B}$ and $S\subseteq \Prim(B)$, respectively.

%

\begin{proposition}[Rieffel]\label{prop-induce}
Let $X$ be an $A$-$B$ imprimitivity bimodule. Then the bijections 
$$\Ind^X:\widehat{B}\to \widehat{A}  \quad\text{and}\quad \Ind^X:\Prim(B)\to \Prim(A)$$
are homeomorphisms.
\end{proposition}

Notice that these homeomorphisms are compatible with the Rieffel-correspondence
(see Proposition \ref{prop-rief-cor}): If $I$ is any closed ideal of 
$B$ and if we identify $\widehat{B}$ with the disjoint union 
$\widehat{I}\cup \widehat{B/I}$ in the canonical way (see \S \ref{subsec-spectrum}),
then induction via $X$ ``decomposes'' into induction via $Y:=X\cdot I$ 
from $\widehat{I}$ to $({\Ind^XI})\dach$ and induction via 
$X/Y$ from $\widehat{B/I}$ to $({A/\Ind^XI})\dach$. This helps to prove

\begin{proposition}\label{prop-typeI}
Suppose that $A$ and $B$ are Morita equivalent $C^*$-algebras. Then 
\begin{enumerate}
\item $A$ is type I
if and only if $B$ is type I.
\item $A$ is CCR if and only if $B$ is CCR.
\item $A$ has continuous trace if and only if $B$ has continuous trace.
%\footnote{By \ref{} we know that
% the property of being a 
%continuous-trace $C^*$-algebra  is also invariant under Morita equivalence.}
\end{enumerate}
\end{proposition}
\begin{proof} Recall from \S \ref{subsec-spectrum} that a $C^*$-algebra 
 $B$ is type I if and only if for each $\pi\in \widehat{B}$ the image $\pi(B)\subseteq\B(H_{\pi})$ contains $\K(H_{\pi})$. Furthermore, $B$ is CCR
if and only if $B$ is type I and points are closed in $\widehat{B}$.

If $X$ is an $A$-$B$ imprimitivity bimodule and $\pi\in \widehat{B}$, 
we may  pass to $B/\ker\pi$ and $A/\ker(\Ind^X\pi)$ via the Rieffel correspondence
to assume that $\pi$ and $\Ind^X\pi$ are injective, and hence
that $B\subseteq \B(H_{\pi})$ and $A\subseteq \B(X\otimes_BH_{\pi})$.
If $B$ is type I, it follows that  $\K:=\K(H_{\pi})$ is an ideal of $B$.
Let $Z:=X\cdot \K$. Then $Z$ is an 
$\Ind^X\K-\K$ imprimitivity bimodule and $Z\otimes_\K H_{\pi}$,
the composition of $Z$ with the $\K-\CC$ imprimitivity bimodule $H_{\pi}$, is  
an $\Ind^X\K-\CC$ imprimitivity bimodule.
It follows that $\Ind^X\K\cong \K(Z\otimes_KH_{\pi})$.
Since $Z\otimes_\K H_{\pi}\cong X\otimes_B H_{\pi}$
via the identity map on both factors, we conclude that $\Ind^X\pi(A)$ contains 
the compact operators $\K(X\otimes_BH_{\pi})$. This proves (i). Now  (ii) follows from (i) since $\widehat{B}$ is homeomorphic to $\widehat{A}$.
The proof of (iii) needs a bit more room and we refer the interested reader to
\cite{Wi-crossed}.
\end{proof}

Of course, similar induction procedures as described above can be 
defined in the equivariant settings: If $(A,G,\alpha)$ is a system, then the 
morphisms from $(A,G,\alpha)$ to $(\CC,G,\id)$ in $\Corr(G)$ are just the  
unitary equivalence classes of non-degenerate
covariant representations of $(A,G,\alpha)$ on Hilbert space, which we shall denote
by $\Rep(A,G)$ (surpressing the given action $\alpha$ in our notation).
Composition with a fixed equivariant morphism $[E,\Phi, u]$ between 
two systems $(A,G,\alpha)$ and $(B,G,\beta)$ 
gives an induction map
$$\Ind^{E}:\Rep(B,G)\to \Rep(A,G); [H, (\pi,U)]\mapsto [E,\Phi,u]\circ [H,\pi,U].$$
As above, we shall write $\Ind^EH:=E\otimes_BH,\;\Ind^E\pi:=\Phi\otimes 1$, and
$\Ind^EU:=u\otimes U$, so that the composition $[E,\Phi,u]\circ [H,\pi,U]$ becomes the 
triple $[\Ind^EH,\Ind^E\pi,\Ind^EU]$.
Taking integrated forms allows to identify
$\Rep(A,G)$ with $\Rep(A\rtimes G)$.
A more or less straight-forward computation gives:

\begin{proposition}\label{prop-cross-ind}
Assume that $[E,\Phi,u]$ is a morphism from $(A,G,\alpha)$ to $(B,G,\beta)$ 
in $\Corr(G)$ and let $[E\rtimes G, \Phi\rtimes G]\in \Mor(A\rtimes G, B\rtimes G)$ 
denote its crossed product. 
Then, for each $[H, (\pi,U)]\in \Rep(B,G)$ we have
$$[\Ind^EH, \Ind^E\pi\times \Ind^EU]=
[\Ind^{E\rtimes G}H, \Ind^{E\rtimes G}(\pi\times U)]\quad\text{in}\quad \Rep(A\rtimes G).$$
Hence induction from $\Rep(B,G)$ to $\Rep(A,G)$ via $[E,\Phi,u]$ 
is equivalent to induction from $\Rep(B\rtimes G)$ to $\Rep(A\rtimes G)$
via $[E\rtimes G, \Phi\rtimes G]$ under the canonical identifications
$\Rep(A,G)\cong \Rep(A\rtimes G)$ and $\Rep(B,G)\cong \Rep(B\rtimes G)$.
\end{proposition}
\begin{proof} Simply check that the map 
$$W: C_c(G,E)\odot H\to E\otimes_B H; \quad
W(\xi\otimes v)=\int_G\xi(s)\otimes U_sv\,ds$$
extends to a unitary from $(E\rtimes G)\otimes_{B\rtimes G}H$ to $E\otimes_BH$
which intertwines both representations (see \cite{E-mor} or \cite{EKQR2} for more details).
\end{proof}
%
%More generally, all
% this extends to the case of  twisted systems $(A,G,N,\alpha,\tau)$, for which 
% the collection $\Rep(A,G,N)$ of all unitary equivalence classes of non-degenerate 
% twisted covariant representations can be identified with the 
% morphisms in $\Corr(G,N)$ from $(A,G,N,\alpha,\tau)$ to $(\CC, G/N, \id)$,
% where the latter is viewed as a twisted system via inflation.
% As above, any morphism $[E,\Phi,u]$ between two twisted systems
% $(A,G,N,\alpha,\tau)$ and $(B,G,N,\beta,\sigma)$ induces an induction map
% from $\Rep(B,G,N)$ to $\Rep(A,G,N)$ via composition of morphisms,
% and Proposition \ref{prop-cross-ind} remains to be true in this 
% more general situation.

We close this section with a brief discussion of corners: 
If $A$ is a C*-algebra and $p\in M(A)$ is a projection, 
then $Ap$ is a Hilbert $pAp$-module with inner product given by 
$\lk ap, bp\rk_{pAp}=pa^*bp$, 
and multiplication from the left turns $Ap$ into 
an $A$-$pAp$ correspondence $[Ap, \phi]$. We then have $\K(Ap)\cong \overline{ApA}$, the 
ideal of $A$ generated by $p$. In a similar way, we may regard $pA$ as an
$pAp$-$A$ correspondence, with inner product given by $\lk pa, pb\rk_{A}=a^*pb$.
Note that $pA$ is then isomorphic to the adjoint module $(Ap)^*$ with 
isomorphism given by $ap\mapsto pa^*$. 

Recall that $p$ is called full, iff $\overline{ApA}=A$. In this case $Ap$ is an $A$-$pAp$ imprimitivity bimodule
and induction from $\Rep(pAp)$ to 
$\Rep(A)$ via $Ap$ gives a bijection between $\Rep(pAp)$ and $\Rep(A)$ with 
inverse given by induction via the adjoint module $(Ap)^*=pA$.

In general, the induction map $\Ind^{Ap}:\Rep(pAp)\to \Rep(A)$ is split injective with 
converse given via {\em compression} by $p$: If $\pi:A\to \B(H_\pi)$ is a non-degenerate representation,
we define $H_{\comp(\pi)}:=\pi(p)H_\pi$ and 
$$\comp(\pi): pAp\to \B(H_{\comp(\pi)}) \quad\text{by}\quad 
\comp(\pi)(pap)=\pi(pap).$$ 
Note that in general, $\comp(\pi)$ could be the zero representation, which 
happens precisely if $\pi(p)=0$. Since $\pi$ is non-degenerate, this is equivalent 
to $\pi(ApA)=0$. 

\begin{proposition}\label{prop-fullp}
Let $p\in M(A)$ be as above. 
Then the following are true:
\begin{enumerate}
\item The compression map $\comp:\Rep(A)\to\Rep(pAp)$ coincides 
with the induction map $\Ind^{pA}:\Rep(A)\to\Rep(pAp)$.
\item For all $\rho\in \Rep(pAp)$ we have $\comp(\Ind^{Ap}\rho)\cong \rho$.
\item $p$ is full if and only if $\comp$ is an inverse for $\Ind^{Ap}$.
\item $p$ is full if and only if $\pi(p)\neq 0$ for all $\pi\in \Rep(A)$.
\end{enumerate}
\end{proposition}
\begin{proof} For (i) just check that for every non-degenerate representation 
$\pi:A\to \B(H_\pi)$ the map
$$pA\otimes_{A}H_\pi\to \pi(p)H_{\pi}; pa\otimes \xi\mapsto \pi(pa)\xi$$
is an isomorphism which intertwines $\Ind^{pA}\pi$ with $\comp(\pi)$.

For the proof of (ii) observe that $pA\otimes_AA_p\cong pAp$ as a $pAp$ bimodule, 
hence by (i) we get $\comp\circ \Ind^{Ap}=\Ind^{pA}\circ \Ind^{Ap}=\Ind^{pAp}=\id_{\Rep(pAp)}$.
%
%
%let $\rho:pAp\to \B(H_\rho)$ be any non-degenerate representation. Then 
%$$H_{\comp(\Ind^{Ap}\!\rho)}=\Ind^{Ap}\!\rho(p)(Ap\otimes_{pAp}H_\rho)=pAp\otimes_{pAp}H_{\rho}\cong H_\rho$$
%via $pap\otimes \xi\mapsto \rho(pap)\xi$, and it is trivial to check that this isomorphism 
%of Hilbert spaces intertwines $\comp(\Ind^{Ap}\!\rho)$ with $\rho$. 

For (iii) we first observe that if $p$ is full, then $Ap$ is an equivalence bimodule and 
induction via $(Ap)^*\cong pA$ is inverse to $\Ind^{Ap}$. Together with (i) this shows 
that $\comp$ is an inverse to $\Ind^{Ap}$. 
 Conversely, if $\overline{ApA}$ is a proper ideal of $A$, there exist 
non-zero, non-degenerate representations $\pi$ of $A$ which vanish on $ApA$, and hence on $p$.
It is then clear that $\comp(\pi)$ is the zero representation, and then 
$\Ind^{Ap}(\comp(\pi))$ is the zero representation as well. Hence $\pi\not\cong\Ind^{Ap}(\comp(\pi))$.

The proof of (iv) is left as an exercise to the reader.
\end{proof}

%\begin{exercise}\label{ex-corner-ideals}
%Let $p\in M(A)$ be a projection. Then show the following statements:
%\begin{enumerate}
%\item If $\rho\in \Rep(pAp)$, then $\ker\Ind^{Ap}\!\rho= \overline{A(\ker\rho)A}$, the ideal in 
%$A$ generated by $\ker\rho\subseteq pAp\subseteq A$.
%As a consequence, If $I$ is an ideal in $pAp$, then $\Ind^{Ap}I= \overline{AIA}\subseteq A$.
%\item If $\pi\in \Rep(A)$, then $\ker(\comp(\pi))=\ker\pi\cap pAp$. As a consequence, if 
%$J$ is an ideal in $A$, then $\Ind^{pA}J=J\cap pAp$.
%\end{enumerate}
%\end{exercise}
%

\subsection{The Fell-topologies and weak containment.}\label{subsec-fell}
For later use and for completeness
 it is necessary to discuss some more topological 
notions on the spaces $\Rep(B)$ and $\I(B)$:
For $I\in \I(B)$ let $U(I):=\{J\in \I(B): J\smallsetminus I\neq \emptyset\}$.
Then  $\{U(I):  I\in \I(B)\}$ is a sub-basis for the 
{\em Fell topology} on $\I(B)$. The Fell topology on $\Rep(B)$ 
is then defined as the inverse image topology with respect to the map
$\ker:\Rep(B)\to \I(B); \pi\mapsto\ker\pi$. \footnote{Recall that $\Rep(B)$ is a set only if we restrict
the cardinality of the Hilbert spaces.} 
The Fell topologies restrict to the
Jacobson topologies on $\Prim(B)$ and $\widehat{B}$, respectively.
Convergence of nets in $\Rep(B)$ (and hence also in $\I(B)$) can conveniently 
be described in terms of {\em weak containment}:
If $\pi\in \Rep(B)$ and $R$ is a subset of $\Rep(B)$, then 
$\pi$ is said to be {\em weakly contained} in $R$ (denoted $\pi\prec R$)
if $\ker\pi\supseteq\cap\{\ker\rho:\rho\in R\}$.
Two subsets $S,R$ of $\Rep(A)$ are said to be {\em weakly equivalent}
($S\sim R$) if $\sigma\prec R$ for all $\sigma\in S$ and $\rho\prec S$ for
all $\rho\in R$. 

\begin{lemma}[Fell]\label{lem-Felltop}
Let $(\pi_j)_{j\in J}$ be a net in $\Rep(B)$ and let $\pi,\rho\in \Rep(B)$. Then
\begin{enumerate}
\item $\pi_j\to\pi$ if and only if $\pi$ is weakly contained in every subnet
of $(\pi_j)_{j\in J}$.
\item If $\pi_j\to \pi$ and if $\rho\prec \pi$, then $\pi_j\to\rho$.
\end{enumerate}
\end{lemma}
For the proof see \cite{Fell2}*{Propositions 1.2 and 1.3}.  As a direct 
consequence 
of this and the fact that induction via bimodules preserves  
inclusion of ideals
we get

\begin{proposition}\label{prop-cont}
Let $[E,\Phi]\in \Mor(A,B)$. Then induction via $E$ preserves weak containment 
and the maps
$$\Ind^E:\Rep(B)\to \Rep(A)\quad\text{and}\quad\Ind^E:\I(B)\to\I(A)$$
are continuous with respect to the Fell topologies. Both maps are 
homeomorphisms if $E$ is an imprimitivity bimodule.
\end{proposition}

Another important observation  is the fact that tensoring representations and ideals
of $C^*$-algebras is continuous:

\begin{proposition}\label{prop-tensor}
Suppose that $A$ and $B$ are $C^*$-algebras. For $\pi\in \Rep(A)$ and $\rho\in \Rep(B)$
let $\pi\otimes \rho\in \Rep(A\otimes B)$ denote the tensor product representation on the 
minimal tensor product $A\otimes B$. Moreover, if $I\in \I(A)$ and $J\in \I(B)$, define
$I\circ  J$ as the closed two-sided ideal of $A\otimes B$ generated by $I\otimes B+A\otimes J$.
Then the maps
\begin{align*}
\Rep(A)&\times \Rep(B)\to \Rep(A\otimes B);\; (\pi,\rho)\mapsto \pi\otimes \rho\\
\text{and} \quad\quad
\I(A)&\times \I(B)\to \I(A\otimes B);\;(I,J)\mapsto I\circ J
\end{align*}
are continuous with respect to the Fell-topologies.
\end{proposition}
\begin{proof} Notice first that if $I=\ker\pi$ and $J=\ker \rho$, then $I\circ  J=\ker(\pi\otimes \rho)$.
Since tensoring ideals clearly preserves inclusion of ideals, the map $(\pi,\rho)\mapsto \pi\otimes \rho$
preserves weak containment in both variables. Hence 
the result follows from Lemma \ref{lem-Felltop}.
\end{proof}

It follows from deep work of Fell 
(e.g. see \cites{Fell1, Dix}) that weak containment (and hence the topologies
on $\I(B)$ and $\Rep(B)$)
can be described completely in terms of matrix coefficients
of representations.  In particular, if $G$ is a locally compact group
and if we identify the collection $\Rep(G)$ of equivalence classes
of unitary representations of $G$ 
with $\Rep(C^*(G))$ via integration, then it is shown in 
\cites{Fell1, Dix} that weak containment for representations of $G$ can be described 
 in terms of convergence of positive definite functions
on $G$ associated to the given representations. 

%Finally, using the canonical identification between $\Rep(A,G)$ and $\Rep(A\rtimes G)$
% (or between $\Rep(A,G,N)$ with $\Rep\big(A\rtimes (G,N)\big)$ in the twisted case) via taking
% integrated forms, it follows from Proposition \ref{prop-cross-ind} that all continuity results
% remain valid for induction via morphisms in $\Corr(G)$ (resp. $\Corr(G,N)$).

% 

\section{Green's imprimitivity theorem and applications}
\label{sec-imp}
\subsection{The imprimitivity theorem}\label{subsec-imp}
We are now presenting (a slight extension of)
Phil Green's imprimitivity
theorem as presented in \cite{Green1}. For this we start with the construction 
of an induction functor
$$\Ind_H^G: \Corr(H)\to\Corr(G); (A,H,\alpha)\mapsto \big(\Ind_H^G(A,\alpha), G, \Ind\alpha\big),$$
 if $H$ is a closed subgroup of $G$ and $\alpha:H\to \Aut(A)$ an action of $H$ on the $C^*$-algebra $A$.
The induced $C^*$-algebra $\Ind_H^G(A,\alpha)$ (or just $\Ind A$ if all data are understood)
is defined as
$$\Ind_H^G(A,\alpha):=
\left\{ f\in C^b(G,A): \begin{matrix} f(sh)=\alpha_{h^{-1}}(f(s)) \;\text{
for all $s\in G, h\in H$}\\ \text{and $(sH\mapsto\|f(s)\|)\in C_0(G/H)$}\end{matrix}\right\},$$
equipped with the pointwise operations and the supremum norm.
The 
induced action $\Ind\alpha:G\to\Aut(\Ind A)$ is
given by 
$$\big(\Ind\alpha_s(f)\big)(t):=f(s^{-1}t)\quad\quad\text{for all}\quad s,t\in G.$$
A similar construction works for morphisms in $\Corr(H)$, i.e., if 
$[E,\Phi,u]$ is a morphism from $(A,H,\alpha)$ to $(B,H,\beta)$,
then  a fairly obvious extension of the above construction yields  the induced morphism $[\Ind_H^G(E,u),\Ind\Phi,\Ind u]$ from $(\Ind_H^G(A,\alpha),G,\Ind\alpha)$ 
to $(\Ind_H^G(B,\beta),G,\Ind\beta)$. One then checks that 
induction preserves composition of morphisms, and hence gives a functor
from $\Corr(H)$ to $\Corr(G)$ (see \cite{EKQR1} for more details).

\begin{remark}\label{rem-Ind}
{\bf (1)} If we start with an action $\alpha:G\to\Aut(A)$ and restrict
this action to the closed subgroup $H$ of $G$, then
$\Ind_H^G(A,\alpha)$ is canonically $G$-isomorphic to 
$C_0(G/H,A)\cong C_0(G/H)\otimes A$ equipped with the diagonal action 
$l\otimes \alpha$, where $l$ denotes the left-translation action 
of $G$ on $G/H$. The isomorphism is given by
$$\Phi: \Ind_H^G(A,\alpha)\to C_0(G/H,A); \quad
\Phi(f)(sH)=\alpha_s(f(s)).$$
{\bf (2)} The construction of the induced algebra $\Ind_H^G(A,\alpha)$
is the $C^*$-analogue of the 
usual construction of the induced $G$-space $G\times_HY$ of 
a topological $H$-space $Y$, which is defined as the quotient 
of $G\times Y$ by the  $H$-action $h(g,y)=(gh^{-1}, hy)$ and which is
equipped with the obvious $G$-action given by the left-translation action on the first factor. 
Indeed, if $Y$ is locally compact, 
then $\Ind_H^GC_0(Y)\cong C_0(G\times_HY)$.
\end{remark}

A useful characterization of induced systems is given 
by the following result:

\begin{theorem}[{cf \cite{E-ind}*{Theorem}}]\label{thm-ind}
Let $(B,G,\beta)$ be a system and let $H$ be a 
closed subgroup of $G$. Then $(B,G,\beta)$ is isomorphic to 
an induced system $(\Ind_H^G(A,\alpha), G,\Ind\alpha)$ 
if and only if there exists a continuous $G$-equivariant map
$\varphi:\Prim(B)\to G/H$,
where $G$ acts on $\Prim(B)$ via $s\cdot P:=\beta_s(P)$.
\end{theorem}

Indeed, we can always define a continuous $G$-map 
$\varphi:\Prim(\Ind A)\to G/H$ by sending a primitive ideal
$P$ to $sH$ iff $P$ contains the ideal
$I_s:=\{f\in\Ind A: f(s)=0\}$. 
Conversely, if $\varphi:\Prim(B)\to G/H$
is given, define $A:=B/I_e$ with
$$I_e:=\cap\{P\in \Prim(B): \varphi(P)=eH\}.$$
Since
$I_e$ is $H$-invariant, the action $\beta|_H$ induces an action 
$\alpha$ of $H$ on $A$
and $(B,G,\beta)$ is isomorphic to $(\Ind_H^G A, G, \Ind\alpha)$
via $b\mapsto f_b\in \Ind_H^G A$; $f_b(s):=
\beta_{s^{-1}}(b)+I_e$. 
We should remark at this point that a much more general result 
has been shown  by 
Le Gall in \cite{LeG}
in the setting of Morita equivalent groupoids.
Applying Theorem \ref{thm-ind}  to commutative $G$-algebras, one gets: 

\begin{corollary}\label{cor-ind}
Let $X$ be a locally compact $G$-space and let $H$ be a closed subgroup
of $G$. Then $X$ is $G$-homeomorphic to $G\times_HY$ for some locally compact $H$-space
$Y$ if and only if there exists a continuous $G$-map $\varphi:X\to G/H$.
If such a map is given, then $Y$ can be chosen as $Y=\varphi^{-1}(\{eH\})$
and the homeomorphism $G\times_HY\cong X$ is given by $[g,y]\mapsto gy$.
\end{corollary}

In what follows let $B_0=C_c(H,A)$ and $D_0=C_c(G,\Ind A)$, viewed as 
dense subalgebras of the full (resp. reduced) 
crossed products $A\rtimes_{(r)} H$ and 
$\Ind A\rtimes_{(r)} G$, respectively.
Let $X_0(A)=C_c(G,A)$. 
We
define left and right module actions of $D_0$ and $B_0$ on $X_0(A)$,
and $D_0$- and $B_0$-valued inner products on $X_0(A)$ by the formulas
\begin{equation}\label{eq-products}
\begin{split}
e \cdot x(s)
&
= \int_G e(t,s) x(t^{-1}s) \Delta_G(t)^{1/2} dt
\\
x\cdot b(s)
&
= \int_H \alpha_h\bigl( x(sh) b(h^{-1}) \bigr) \Delta_H(h)^{-1/2} dh
\\
_{D_0}\lk x,y\rk(s,t)
&
= \Delta_G(s)^{-1/2}
\int_H \alpha_h\bigl( x(th) y(s^{-1}th)^* \bigr) \,dh
\\
\lk x,y\rk_{B_0}(h)
&
= \Delta_H(h)^{-1/2}
\int_G x(t^{-1})^* \alpha_h(y(t^{-1}h)) \,dt,
\end{split}
\end{equation}
for $e\in D_0, x,y\in X_0(A)$, and $b\in B_0$.
The $C_c(H,A)$-valued inner product on $X_0(A)$ 
provides $X_0(A)$  with two different norms: 
$\|\xi\|_{\max}^2:=\|\lk \xi,\xi\rk_{B_0}\|_{\max}$ and 
$\|\xi\|_{r}^2:=\|\lk \xi,\xi\rk_{B_0}\|_{r}$,
where $\|\cdot\|_{\max}$ and $\|\cdot\|_{r}$ denote 
the maximal and reduced norms on $C_c(G,A)$.
Then Green's imprimitivity theorem reads as follows:

\begin{theorem}[{Green}]\label{thm-Green}
The actions and inner products on $X_0(A)$ extend to the completion 
$X_H^G(A)$ of $X_0(A)$ with respect to $\|\cdot\|_{\max}$ 
such that $X_H^G(A)$
becomes an
  $(\Ind_H^GA\rtimes G)$-$(A\rtimes H)$
imprimitivity bimodule.  

Similarly, the completion
$X_H^G(A)_r$ of $X_0(A)$ with respect to $\|\cdot\|_{r}$
becomes an $(\Ind_H^GA\rtimes_rG)$-$(A\rtimes_rH)$
imprimitivity bimodule.
\end{theorem}

\begin{remark}\label{rem-Green}
{\bf (1)} Although the statement of Green's theorem looks quite straightforward,
the proof requires a fair amount of work. The main 
problem is to show positivity of the inner products and 
continuiuty of the left and right actions of $D_0$ and $B_0$ 
on $X_0$ with respect to the appropriate norms.
In \cite{Green1} Green only considered full crossed products.
The reduced versions were obtained later by Kasparov 
(\cite{K3}), by Quigg and Spielberg (\cite{QS})
and by Kirchberg and Wassermann (\cite{KW2}).

The reduced module $X_H^G(A)_r$ can also be realized as the quotient
of $X_H^G(A)$ by the  submodule $Y:=X_H^G(A)\cdot I$ with 
$I:=\ker\big(A\rtimes H\to A\rtimes_rH\big)$.
This follows from the fact
that the ideal $I$ corresponds to the ideal 
$J:=\ker\big(\Ind A\rtimes G\to \Ind A\rtimes_rG\big)$ in $\Ind A\rtimes G$ via
the Rieffel correspondence (see Proposition
\ref{prop-rief-cor}). We shall give an argument for this fact
 in Remark \ref{rem-reducedimp} below.

{\bf (2)} In his original work \cite{Green1}, Green first considered the special 
case where the action of $H$ on $A$ restricts from an action of 
$G$ on $A$. In this case one obtains a 
Morita equivalence between $A\rtimes_{(r)}H$ and 
$C_0(G/H,A)\rtimes_{(r)}G$ (compare with Remark \ref{rem-Ind} above).
Green then deduced from this 
a more general result (see \cite{Green1}*{Theorem 17}),
which by Theorem \ref{thm-ind}
is equivalent to the above formulation for full crossed products.

{\bf (3)} In \cite{EKQR1} it is shown that the construction of the equivalence 
bimodule $X(A)$, viewed as an isomorphism in the correspondence category $\Corr$,
provides a natural equivalence between the descent functor
$\rtimes:\Corr(H)\to \Corr; (A,H,\alpha)\mapsto A\rtimes H$ with the composition 
$\rtimes\circ \Ind_H^G:\Corr(H)\to \Corr; (A,H,\alpha)\mapsto \Ind A\rtimes G$
(and similarly for the reduced descent functors $\rtimes_r$). This shows 
that the assignment $(A,H,\alpha)\mapsto X_H^G(A)$ is, in a very strong sense,
natural in $A$.
\end{remark}

Let us now present
some basic examples:

\begin{example}\label{ex-Green}
{\bf (1)} Let $H$ be a closed subgroup of $G$. Consider the 
trivial action of $H$ on $\CC$. Then $\Ind_H^G\CC=C_0(G/H)$
and Green's theorem provides a Morita equivalence 
between $C^*(H)$ and $C_0(G/H)\rtimes G$, and similarly
between $C_r^*(H)$ and $C_0(G/H)\rtimes_rG$. 
It follows then from Proposition \ref{prop-cont} that induction via 
$X_H^G(\CC)$  identifies the representation 
spaces $\Rep(H)$ and $\Rep(C_0(G/H),G)$. 
This is a very strong version of Mackey's 
original imprimitivity theorem  for groups 
(e.g., see \cites{Ma1, Ma2, Blat1}). 
%Notice that it follows from the Rieffel correspondence together with
%part (2) of Remark \ref{rem-Green} and Corollary \ref{cor-amenable}
%that $C_0(G/H)\rtimes G=C_0(G/H)\rtimes_rG$ if and only if  $H$
%is amenable.

{\bf (2)} If $H=\{e\}$ is the trivial subgroup of $G$, we obtain a Morita equivalence between 
$A$ and $C_0(G,A)\rtimes G$, where $G$ acts on itself by left translation. 
Indeed, in this case we obtain a unitary isomorphism between Green's bimodule 
$X_{\{e\}}^G(A)$ and the Hilbert $A$-module $L^2(G,A)\cong A\otimes L^2(G)$
via the transformation 
$$U:X_{\{e\}}^G(A)\to L^2(G,A); \big(U(x)\big)(s)=\Delta(s)^{-\frac{1}{2}}x(s).$$
It follows from this that $C_0(G,A)\rtimes G$ is isomorphic to 
$\K(A\otimes L^2(G))\cong A\otimes \K(L^2(G))$.
In particular, it follows that $C_0(G)\rtimes G$ is isomorphic to $\K(L^2(G))$
if $G$ acts on itself by translation.

Since full and reduced crossed products by the trivial group coincide, it follows 
from part (2) of Remark \ref{rem-Green}  
that $C_0(G,A)\rtimes_rG\cong C_0(G,A)\rtimes G$, and hence that 
$C_0(G,A)\rtimes_rG\cong A\otimes \K(L^2(G))$, too.

%{\bf (2)} If $H=\{e\}$ is the trivial subgroup of $G$ in the above example, we obtain a Morita equivalence between 
%$\CC=C^*(\{e\})$ and $C_0(G)\rtimes G$, where $G$ acts on itself by left translation. In particular, it follows that $C_0(G)\rtimes G$ 
%has a unique equivalence class of irreducible representations.
%A representative is given by the pair 
%$(M,\lambda)$, where 
%$M:C_0(G)\to\B(L^2(G))$ denotes the representation by multiplication 
%operators and $\lambda:G\to U(L^2(G))$ is the regular representation 
%of $G$. If $f\in C_c(G\times G)\subseteq C_c(G,C_0(G))$, then 
%$$\big(M\times\lambda(f)\xi\big)(s)=\int_G f(s, t^{-1}s)\xi(t^{-1}s)\,dt$$
%is an integral operator whose kernel has compact support in $G\times G$.
%Since such operators 
%form a dense subset of $\K(L^2(G))$, it follows  that
%$C_0(G)\rtimes G\cong M\times\lambda(C_0(G)\rtimes G)=\K(L^2(G))$. 
%\\
{\bf (3)} Let $H_3$ denote the three-dimensional real Heisenberg group,
i.e., $H_3=\RR^2\rtimes \RR$ with action of $\RR$ on $\RR^2$ 
given by $x\cdot(y,z)=(y, z+xy)$. We want to use Green's theorem to 
analyse the structure of $C^*(H_3)\cong C^*(\RR^2)\rtimes \RR$. 
We first identify $C^*(\RR^2)$ with $C_0(\RR^2)$ via Fourier
transform. The transformed action of $\RR$ 
on $\RR^2$ is then given by $x\cdot(\eta,\zeta)=(\eta- x\zeta,\zeta)$.
The short exact sequence
$$0\to C_0(\RR\times\RR^*)
\to C_0(\RR^2)\to C_0(\RR\times\{0\})\to 0$$
determines a short exact sequence
$$0\to C_0\big(\RR\times\RR^*)\rtimes\RR
\to C^*(H_3)\to C_0(\RR\times\{0\})\rtimes\RR\to 0.$$
Since the action of $\RR$ on the quotient $C_0(\RR)\cong C_0(\RR\times\{0\})$
is trivial, we see that 
$C_0(\RR\times\{0\})\rtimes\RR\cong C_0(\RR)\otimes C^*(\RR)\cong C_0(\RR^2)$.
The homeomorphism $h:\RR\times\RR^*\to \RR\times \RR^*; h(\eta,\zeta)=
(-\frac{\eta}{\zeta},\zeta)$ transforms the action of $\RR$ on 
$C_0(\RR\times \RR^*)\cong C_0(\RR)\otimes C_0(\RR^*)$ to the diagonal action 
$l\otimes\id$, where $l$ denotes left translation.
Thus, it follows from (2) and Lemma \ref{lem-tensor} that 
$C_0(\RR\times\RR^*)\rtimes\RR\cong 
C_0(\RR^*)\otimes \K(L^2(\RR))$
and  we obtain a short exact sequence
$$0\to C_0(\RR^*)\otimes\K(L^2(\RR))\to C^*(H_3)\to C_0(\RR^2)\to 0$$
for $C^*(H_3)$.

{\bf (4)} Let $\RR$ act on the two-torus $\TT^2$ by an irrational flow,
i.e. there exists an irrational number $\theta\in (0,1)$ such that
$t\cdot(z_1, z_2)=(e^{2\pi it}z_1, e^{ 2\pi i\theta t}z_2)$.
Then $\TT^2$ is $\RR$-homeomorphic to the 
induced space $\RR\times_{\ZZ}\TT$, where $\ZZ$ acts on $\TT$ 
by irrational rotation given by $\theta$ (compare with 
Example  \ref{ex-crossed}). Hence, it follows from 
Green's theorem that $C(\TT^2)\rtimes_{\theta}\RR$ is Morita equivalent
to the irrational rotation algebra $A_{\theta}=C(\TT)\rtimes_{\theta}\ZZ$.
\end{example}

In \cite{Green2} Phil Green shows that for second countable $G$ we always have 
a decomposition $\Ind_H^GA\rtimes_{\Ind\alpha}G\cong (A\rtimes_{\alpha}H)\otimes \K(L^2(G/H))$,
where the $L^2$-space is taken with respect to some quasi invariant measure on $G/H$.
In case where $H=\{e\}$ is the trivial group, this is part (2) of Example \ref{ex-Green} above, but the 
general proof is more difficult because of some measure theoretic technicalities.
But if $H$ is open in $G$, the proof of Green's structure theorem becomes quite easy:

\begin{proposition}\label{prop-Green-structure}
Let $H$ be an open subgroup of $G$ and let $\alpha:H\to\Aut(A)$ be an action.
Let $c:G/H\to G$ be a cross-section for the quotient map $q:G\to G/H$ such that $c(eH)=e$.
Then there is an isomorphism 
$\Psi: X_H^G(A)\stackrel{\cong}{\to} (A\rtimes_\alpha H)\otimes \ell^2(G/H)$ of Hilbert $A\rtimes_\alpha H$ modules
given on the dense subspace $X_0(A)=C_c(G,A)$ by
$$\Psi(x)=\sum_{G/H} x_{sH}\otimes \delta_{sH},$$
where $x_{sH}\in C_c(H,A)$ is defined by $x_{sH}(h)=\Delta(c(s)h)^{-1/2}\alpha_{h}(x(c(s)h))$.
A similar decomposition holds for the reduced module: $X_H^G(A)_r\cong (A\rtimes_rH)\otimes \ell^2(G/H)$.
As a consequence, we get isomorphisms 
\begin{align*}
&\Ind_H^GA\rtimes G\cong (A\rtimes_{\alpha}H)\otimes \K(\ell^2(G/H))\quad\text{and}\; \\
&\Ind_H^GA\rtimes_r G\cong (A\rtimes_{\alpha,r} H)\otimes\K(\ell^2(G/H)).
\end{align*}
\end{proposition}
\begin{proof} It is easy to check that the mapping $x\mapsto \Psi(x)$ is a bijection between 
$C_c(G,A)$ and the algebraic tensor product $C_c(H,A)\odot C_c(G/H)$, which is  dense 
in $(A\rtimes_{\alpha}H)\otimes \ell^2(G/H)$. It therefore suffices 
to check that 
$\lk \Psi(x),\Psi(y)\rk_{A\rtimes H}=\lk x,y\rk_{A\rtimes H}$ for all $x,y\in C_c(G,A)$ and that 
$\Psi$ intertwines the right action of $A\rtimes_\alpha H$ on both modules. We show the first and 
leave the second as an exercise for the reader:
Note that the inner products on both dense subspaces take values in $C_c(H,A)$.  Since
 $H$ is open in $G$ we have $\Delta_G=\Delta_H$ and the 
formula $\int_G f(t)\,dt=\sum_{tH\in G/H} \int_H f(tl)\,dl$ for all $f\in C_c(G,A)$.  
We then compute 
for $x,y\in C_c(G,A)$ and $h\in H$:
\begin{align*}
\lk \Psi(x), \Psi(y)\rk(h)&= \sum_{G/H} (x_{sH}^**y_{sH})(h)=\sum_{G/H}\int_H\alpha_l( x_{sH}(l)^*y_{sH}(lh))\,dl\\
&=\sum_{G/H} \Delta(h)^{-1/2}\int_H\Delta(c(s)l)^{-1}x(c(s)l)^*\alpha_h(y(c(s)lh))\,dl\\
&\stackrel{l\mapsto c(s)^{-1}sl}{=}\sum_{G/H} \Delta(h)^{-1/2}\int_H\Delta(sl)^{-1}x(sl)^*\alpha_h(y(slh))\,dl\\
&=\int_G\Delta(h)^{-1/2}\Delta(s)^{-1}x(s)^*\alpha_h(y(sh))\,ds\\
&=\Delta(h)^{-1/2}\int_G x(s^{-1})^*\alpha_h(y(s^{-1}h))\,ds=\lk x,y\rk(h).
\end{align*}
The decomposition of $\Ind_H^GA\rtimes G$ now follows from 
$$\Ind_H^GA\rtimes G=\K(X_H^G(A))\cong \K((A\rtimes H)\otimes \ell^2(G/H))=(A\rtimes H)\otimes \K(\ell^2(G)),$$
and similarly for $\Ind_H^GA\rtimes_rG$.
\end{proof}

We close this section with a proof of Green's imprimitivity theorem in the case 
where $H$ is open in $G$. 
For this let $\alpha:H\to \Aut(A)$ be an action of $H$ on a C*-algebra $A$. 
We shall see that in this case there exists a canonical full projection $p\in M(\Ind_H^GA\rtimes G)$ 
such that the crossed product 
$A\rtimes_{\alpha}H$ is isomorphic to the corner $p(\Ind_H^GA\rtimes_{\Ind\alpha}G)p$.
Of course, if we believe in the validity of Green's theorem, this follows basically 
from Proposition \ref{prop-Green-structure}.

Consider the canonical embedding $\Psi:C_0(G/H)\to M(\Ind_H^GA)$ given by the (central) action
$$(\Psi(\varphi)F)(t)=\varphi(tH)F(t),\quad \varphi\in C_0(G/H), F\in \Ind_H^GA, t\in G.$$
In what follows we shall often write $\varphi\cdot F$ for $\Psi(\varphi)F$.
Let $\tilde{p}=\varphi(\delta_{eH})\in M(\Ind_H^GA)$ and let $p$ be the image of $\tilde{p}$
under the extension to $M(\Ind_H^GA)$ of the embedding $i_{\Ind_H^GA}:\Ind_H^GA\to M(\Ind_H^GA\rtimes_{\Ind\alpha}G)$.

\begin{proposition}\label{prop-corner}
Let $p\in M(\Ind_H^GA\rtimes_{\Ind\alpha} G)$ be as above. Then $p$ is a full projection 
in $M(\Ind_H^GA\rtimes_{\Ind\alpha} G)$ such that there is a canonical isomorphism 
$$A\rtimes_{\alpha}H\cong p(\Ind_H^GA\rtimes_{\Ind\alpha}G)p.$$
Moreover, the resulting $\Ind_H^GA\rtimes_{\Ind\alpha}G$-$A\rtimes_{\alpha}H$ equivalence bimodule \linebreak
$(\Ind_H^GA\rtimes_{\Ind\alpha}G)p$ is  isomorphic to Green's equivalence bimodule 
 $X_H^G(A)$ of Theorem \ref{thm-imp}.

A similar result holds for the reduced crossed products. 
\end{proposition}

\begin{proof}
We first observe that $A$ can be identified with the direct summand $\tilde{p}\Ind_H^GA$ of  $\Ind_H^GA$ via 
the embedding $\psi:A\to \tilde{p}\Ind_H^GA; a\mapsto F_a$ with 
\begin{equation}\label{eq-summand}
F_a(s)=\left\{\begin{matrix} \alpha_{s^{-1}}(a)& \text{if $s\in H$}\\ 0&\text{else}\end{matrix}\right\}.
\end{equation}

If $F\in C_c(G,\Ind_H^GA)$, then using the formulas in Remark \ref{rem-universal} we compute
\begin{equation}\label{eq-incl}
(pFp)(s)=\tilde{p}F(s)\Ind\alpha_s(\tilde{p})=\delta_{eH}F(s)\delta_{sH}=\left\{\begin{matrix}\delta_{eH}F(s)&\text{if $s\in H$}\\ 0&\text{else}\end{matrix}\right\}.
\end{equation}
Using (\ref{eq-summand}) and (\ref{eq-incl}) one easily checks that there is a canonical $*$-isomorphism 
$$\Phi: C_c(H,A)\stackrel{\cong}{\to} pC_c(G,\Ind_H^GA)p\subseteq p(\Ind_H^GA\rtimes_{\Ind\alpha}G)p$$
which maps a function $f\in C_c(H,A)$ to the function $F\in C_c(G,\Ind_H^GA)$ given by
\begin{equation}\label{eq-corner} 
F(s,t)=\left\{\begin{matrix} \alpha_{t^{-1}}(f(s))& \text{if $s,t\in H$}\\ 0&\text{else}\end{matrix}\right\}.
\end{equation}
Indeed, by a straightforward but lengthy computation one checks that
$\Phi$ coincides with the integrated form of the covariant homomorphism
$(i_{\Ind A}\circ \psi, i_G|_H)$ of $(A,H,\alpha)$ into $M(\Ind_H^GA\rtimes_{\Ind\alpha} G)$, where $\psi$ is as above and 
$(i_{\Ind A}, i_G)$ is the canonical covariant homomorphism of $(\Ind_H^GA, G,\Ind\alpha)$ into $M(\Ind_H^GA\rtimes G)$.
It follows from this that $\Phi$ extends to a surjective $*$-homomorphism $\Phi:A\rtimes_{\alpha}H\onto p(\Ind_H^GA\rtimes_{\Ind\alpha}G)p$. 

To see that $\Phi$ is injective let $(\rho,V)$ be any covariant representation of $(A,H,\alpha)$. Then we construct an induced representation 
$(\Ind\rho, \Ind V)$ of $(\Ind_H^GA, G,\Ind\alpha)$ as follows: We define 
$$H_{\Ind V}=\left\{\xi:G\to H_\rho: \begin{matrix} \text{$\xi(th)=V_{h^{-1}}\xi(t)$ for all $t\in G, h\in H$}\\
\text{ and $\sum_{tH\in G/H} \|\xi(t)\|^2<\infty$.}\end{matrix}\right\}$$
equipped with the inner product 
$$\lk\!\lk \xi, \eta\rk\!\rk=\sum_{tH\in G/H} \lk \xi(t), \eta(t)\rk.$$
Note that this sum is well defined since $\xi(th)=V_{h^{-1}}\xi(t)$ for all $t\in G, h\in H$.
We then define $(\Ind\rho, \Ind V)$ by
$$(\Ind\rho(F)\xi)(t)=F(t)\xi(t)\quad\text{and}\quad (\Ind V_s\xi)(t)=\xi(s^{-1}t),$$
for $F\in \Ind_H^GA$ and $s\in G$. It is then straightforward to check the following items:
\begin{itemize}
\item $(\Ind\rho, \Ind V)$ is a covariant representation of $(\Ind_H^GA, G,\Ind\alpha)$.
\item The composition of the compression $\comp(\Ind\rho\rtimes \Ind V)$ of $\Ind\rho\rtimes \Ind V$ to the corner
$p(\Ind_H^GA\rtimes G)p$ with $\Phi: A\rtimes H\to p(\Ind_H^GA\rtimes G)p$ is equivalent to $\rho\rtimes V$.
\end{itemize}
Hence, if we choose $\rho\rtimes V$ to be faithful on $A\rtimes_{\alpha}H$, we see that $\Phi$ must be faithful as well.

To check that $p$ is a full projection it suffices to show that no non-zero representation $\pi\rtimes U$ of 
$\Ind_H^GA\rtimes G$ vanishes on $p$. By definition of $p$ we have $\pi\rtimes U(p)=\pi(\tilde p)=\pi(\delta_{eH})$.
So assume to the contrary that $\pi(\delta_{eH})=0$, where we regard $C_0(G/H)$ as a subalgebra of $M(\Ind_H^GA)$ 
as described in the discussion preceding the proposition. Then for all $t\in G$ we have
$$\pi(\delta_{tH})=\pi(\Ind\alpha_t\delta_{eH})=U_t\pi(\delta_{eH})U_t^*=0$$
as well, and since $\sum_{tH\in G/H}\delta_{tH}$ converges strictly to $1$ in $M(\Ind_H^GA)$, it follows from 
this that $\pi$ is the zero-representation. But then $\pi\rtimes U$ is zero as well, which contradicts our assumption.

We now have seen that $A\rtimes_\alpha H$ is isomorphic to the full corner $p(\Ind_H^GA\rtimes G)p$.
We want to compare Green's module $X_H^G(A)$ with the module $(\Ind_H^GA\rtimes G)p$.
For this we first compute for $F\in C_c(G, \Ind_H^GA)$ that
$$(Fp)(s,t)=\big(F(s)\cdot\Ind\alpha_s(\delta_{eH})\big)(t)= \big(F(s)\cdot \delta_{sH}\big)(t)=\left\{\begin{matrix} F(s,t)&\text{if $t\in sH$}\\ 0& \text{else}
\end{matrix}\right\}.$$
Therefore, because of the condition $F(s,th)=\alpha_{h^{-1}}(F(s,t))$, it follows that $Fp$ is completely determined 
by the values of $F$ on the diagonal $\Delta_G=\{(s,s): s\in G\}$. 
Recall that $X_H^G(A)$ is the completion of $X_0(A)=C_c(G,A)$ with respect to the $C_c(H,A)$ valued 
inner product
$$\lk x, y\rk(h)=\Delta_H(h)^{-1/2}\int_G x(t^{-1})^*\alpha_h(y(t^{-1}h))\,dt.$$
We then obtain a well defined bijective map
$$\Theta: C_c(G,\Ind_H^GA)p\to C_c(G,A); Fp\mapsto (s\mapsto \sqrt{\Delta_G(s)}F(s,s)).$$
To show that $\Theta$ preserves the $A\rtimes_{\alpha}H$-valued inner products, it follows 
from (\ref{eq-corner}) that we need to check that for all $F_1,F_2\in C_c(G,\Ind_H^GA)$
and all $(h,l)\in H$ we have
$$\alpha_{l^{-1}}(\lk \Theta(F_1),\Theta(F_2)\rk(h))=F_1^**F_2(h,l)=\int_G\Delta_G(t^{-1}) F_1(t^{-1}, t^{-1}l)^*F_2(t^{-1}h, t^{-1}l)\, dt.$$
But this follows from a straightforward calculation using that $\Delta_G|_H=\Delta_H$, since $H$ is open in $G$.
One also easily checks that $\Theta$ intertwines the left action of $\Ind_H^GA\rtimes G$ on both modules.

In order to prove the analogue for the reduced case one checks that compression of a regular representation 
of $\Ind_H^GA\rtimes G$ gives a regular representation of $A\rtimes H$. We leave this as an exercise 
for the reader (or see Remark \ref{rem-reducedimp} below).
\end{proof}

\begin{remark}\label{rem-induced} 
{\bf (1)} If $A$ is a $G$-algebra and $H$ is a closed subgroup of $G$, then we saw in Remark \ref{rem-Ind} that 
$\Ind_H^GA$ is isomorphic to $C_0(G/H,A)$. Thus, if $H$ is also open in $G$, it follows  from 
Proposition \ref{prop-corner} that $A\rtimes_{\alpha}H$ is a full corner in $C_0(G/H,A)\rtimes_{\tau\otimes \alpha} G$,
and similarly for the reduced crossed products. 

{\bf (2)} In  \cite{KK}*{Section 3.5.3}, we need to 
investigate the structure of crossed products of the form $C_0(I, A)\rtimes G$ in which 
$I$ is a {\em discrete} $G$-space, $A$ is a $G$-algebra, and $G$ act diagonally on 
$C_0(I,A)\cong C_0(I)\otimes A$. In this case we can decompose 
$I$ as a disjoint union of $G$-orbits $Gi=\{si:s\in G\}$ which induces a direct sum decomposition 
$$C_0(I, A)\rtimes G\cong \oplus_{G\backslash I} C_0(Gi, A)\rtimes G.$$
If $G_i=\{s\in G: si=i\}$ denotes the stabiliser of $i\in I$ for the action of $G$ (which is open in $G$ since $I$ is discrete), we get 
$G$-equivariant bijections $G/G_i\cong Gi; sG_i\mapsto si$, and then the above decomposition becomes
$$C_0(I, A)\rtimes G\cong \oplus_{G\backslash I} C_0(G/G_i, A)\rtimes G.$$
Now by Green's imprimitivity theorem (or by Proposition \ref{prop-corner}) each summand
$C_0(G/G_i, A)\rtimes G$ is Morita equivalent to $A\rtimes_{\alpha}G_i$, and hence we 
see that $C_0(I, A)\rtimes G$ is Morita equivalent to $\oplus_{G\backslash I} A\rtimes_{\alpha}G_i$.
Indeed, if $p_i\in M(C_0(G/G_i,A)\rtimes G)$ is the full projection as in Proposition \ref{prop-corner},
we observe that the sum $\sum_{G\backslash I} p_i$ converges strictly in $M(C_0(I,A)\rtimes G)$ to a projection $p$ and then
$\oplus_{G\backslash I} A\rtimes_{\alpha}G_i$ is isomorphic to the full corner $p(C_0(I,A)\rtimes G)p$ in 
$C_0(I,A)\rtimes G$. All this goes through without change for the reduced crossed products.
\end{remark}

\subsection{The Takesaki-Takai duality theorem.}\label{subsec-Tak}
From part (2) of Example \ref{ex-Green} it is fairly easy to obtain the 
Takesaki-Takai duality theorem for crossed products by abelian groups.
For this assume that $(A,G,\alpha)$ is a system with $G$ abelian. The
{\em dual action} $\widehat{\alpha}:\widehat{G}\to \Aut(A\rtimes G)$ of the dual group
$\widehat{G}$
on the crossed product $A\rtimes G$ is defined by
$$\widehat{\alpha}_{\chi}(f)(s):=\chi(s)f(s)\quad\quad \text{for}\;\chi\in \widehat{G}\;\text{and}\;f\in C_c(G,A)\subseteq A\rtimes G.$$
With a similar action of $\widehat{G}$ on crossed products $E\rtimes G$ for an equivariant bimodule $(E,\Phi,u)$ we obtain from this a descent functor
$$\rtimes:\Corr(G)\to \Corr(\widehat{G}).$$
The  double dual crossed product 
$A\rtimes G\rtimes\widehat{G}$ is isomorphic to
$C_0(G,A)\rtimes G$ with respect to the diagonal action $l\otimes \alpha$
of $G$ on $C_0(G,A)\cong C_0(G)\otimes A$.
Indeed, we have canonical (covariant) representations 
$(k_A,k_G, k_{\widehat{G}})$ of the triple 
$(A,G,\widehat{G})$ into $M\big(C_0(G,A)\rtimes G\big)$
given by the formulas
\begin{align*}
\big(k_A(a)\cdot f\big)(s,t)=&a\big(f(s,t)\big),\quad \big(k_G(r)\cdot f\big)(s,t)=
\alpha_r\big(f(r^{-1}s,r^{-1}t)\big),\quad\text{and}\\
&\big(k_{\widehat{G}}(\chi)\cdot f\big)(s,t)=\chi(t)f(s,t),
\end{align*}
for $f$ in the dense subalgebra $C_c(G, C_0(G,A))$ of $C_0(G,A)\rtimes G$.
Making extensive use of the universal properties, one checks 
 that the integrated form
$$(k_A\times k_G)\times k_{\widehat{G}}:(A\rtimes G)\rtimes \widehat{G}\to 
M(C_0(G,A)\rtimes G)$$ 
gives the desired isomorphism $A\rtimes G\rtimes\widehat{G}\cong C_0(G,A)\rtimes G$.
Using the isomorphism $C_0(G,A)\rtimes G\cong A\otimes\K(L^2(G))$ of 
Example \ref{ex-Green} (2) and checking what this isomorphism does on the 
double-dual action $\widehat{\widehat{\alpha}}$ we arrive at

\begin{theorem}[Takesaki-Takai]\label{thm-taktak}
Suppose that $(A,G,\alpha)$ is a system with $G$ abelian. Then the double dual 
system $(A\rtimes G\rtimes\widehat{G}, G,\widehat{\widehat{\alpha}})$
is equivariantly isomorphic to the system $(A\otimes \K(L^2(G)), G, \alpha\otimes\Ad\rho)$,
where $\rho:G\to U(L^2(G))$ denotes the right regular representation 
of $G$ on $L^2(G)$.
\end{theorem}

Recall that the {\em right regular representation} $\rho:G\to U(L^2(G))$ is defined by 
$(\rho_t\xi)(s)=\sqrt{\Delta(t)}\xi(st)$ for $\xi\in L^2(G)$ (but if $G$ is abelian, the modular function $\Delta$ dissapears).
Note that the system $(A\otimes \K(L^2(G)), G, \alpha\otimes\Ad\rho)$ in the Takesaki-Takai
theorem is Morita equivalent
to the original system $(A,G,\alpha)$ via the equivariant  imprimitivity bimodule  
$(A\otimes L^2(G), \alpha\otimes\rho)$. In fact, the assignment $(A,G,\alpha)\mapsto
(A\otimes L^2(G),\alpha\otimes\rho)$ is easily seen to give a natural equivalence between 
the identity functor on $\Corr(G)$ and the composition
$$
\begin{CD}
\Corr(G)  @>\rtimes >> \Corr(\widehat{G})  @>\rtimes >> \Corr(G).
\end{CD}
$$

In general, if $G$ is not abelian, one can obtain similar duality theorems 
 by replacing the dual action of $\widehat{G}$ 
by a dual coaction of the group algebra $C^*(G)$ on $A\rtimes G$. A fairly complete
account of that theory in the group case is given in the appendix of 
\cite{EKQR2} --- however a much more  general duality theory for Hopf-$C^*$-algebras
was developed by Baaj and Skandalis in \cite{BS2} and Kustermans and Vaes in \cite{KV}.

\subsection{Permanence properties of exact groups}\label{subsec-permanence}
As a further  application of Green's imprimitivity theorem we now 
present  some of Kirchberg's and Wassermann's permanence results for $C^*$-exact groups.
Recall from Definition \ref{def-exact} that a group $G$ is called $C^*$-exact
(or just exact)
if for every system $(A,G,\alpha)$ and for every $G$-invariant ideal $I\subseteq A$
the sequence
$$0\to I\rtimes_rG\to A\rtimes_rG\to (A/I)\rtimes_rG\to 0$$
is exact (which is equivalent to exactness of the sequence in the middle term).
Recall from Proposition \ref{prop-exactfull} that the corresponding sequence
of full crossed products is always exact. Using Proposition \ref{prop-amenable},
this implies that all amenable groups are exact.

In what follows we want to relate exactness of $G$ with exactnesss of a closed 
subgroup $H$ of $G$. For this we start with a system $(A,H,\alpha)$ and a 
closed $H$-invariant ideal $I$ of $A$. Recall that Green's 
$\Ind A\rtimes_rG-A\rtimes_rH$ imprimitivity bimodule $X_H^G(A)_r$ 
is a completion of $C_c(G,A)$. Using the formulas for the actions 
and inner products as given in (\ref{eq-products}) one observes
that $X_H^G(I)_r$ can be identified
with the closure of $C_c(G,I)\subseteq C_c(G,A)$ in $X_H^G(A)_r$.
It follows that  the ideals $\Ind I\rtimes_r G$ and $I\rtimes_rH$ 
are linked via the Rieffel correspondence with respect to
$X_H^G(A)_r$ (see Proposition \ref{prop-rief-cor}).
Similarly, the imprimitivity bimodule $X_H^G(A/I)_r$ is isomorphic 
to the quotient $X_H^G(A)_r/Y$ with 
$Y:=X_H^G(A)_r\cdot\ker\big(A\rtimes_rH\to A/I\rtimes_rH\big)$,
which implies that the ideals
$$\ker\big(A\rtimes_rH\to A/I\rtimes_rH\big)\quad\text{and}\quad
\ker\big(\Ind A\rtimes_rG\to \Ind(A/I)\rtimes_rG\big)$$
are also linked via the Rieffel correspondence.
Since the Rieffel correspondence is one-to-one, we obtain
\begin{equation}\label{eq-correspond}
\begin{split}
I\rtimes_rH&=\ker\big(A\rtimes_rH\to A/I\rtimes_rH\big)\\
\Longleftrightarrow\;\;\Ind I\rtimes_rG&=\ker\big(\Ind A\rtimes_rG\to \Ind(A/I)\rtimes_rG\big).
\end{split}
\end{equation}
Using this, we now  give proofs of two of the main results of \cite{KW2}.

\begin{theorem}[Kirchberg and Wassermann]\label{thm-KW1}
Let $G$ be a locally compact group. Then the following are true:
\begin{enumerate}
\item If $G$ is exact and $H$ is a closed subgroup of $G$, then $H$ is exact.
\item Let $H$ be a closed subgroup of $G$ such that $G/H$ is compact. Then $H$ 
exact implies $G$ exact.
\end{enumerate}
\end{theorem}
\begin{proof} Suppose that $I$ is an $H$-invariant ideal of the $H$-algebra 
$A$. If $G$ is exact, then $\Ind I\rtimes_rG=\ker\big(\Ind A\rtimes_rG\to \Ind(A/I)\rtimes_rG\big)$ 
and hence $I\rtimes_r H=\ker\big(A\rtimes_rH\to A/I\rtimes_rH\big)$ by  (\ref{eq-correspond}). 
This proves (i).

To see (ii) we start with an arbitrary $G$-algebra $A$ and a $G$-invariant ideal 
$I$ of $A$. Since $A$, $I$, and $A/I$ are $G$-algebras and $G/H$ is compact,
we have $\Ind_H^G A\cong C(G/H,A)$ and similar statements hold  for $I$ and $A/I$.
 Since $H$ is exact we see that  the lower row of the commutative diagram
{\small $$
\begin{CD}
0 @>>> I\rtimes_rG @>>> A\rtimes_rG @>>> (A/I)\rtimes_rG @>>> 0\\
@.  @VVV  @VVV @VVV \\
0 @>>> C(G/H,I)\rtimes_rG @>>> C(G/H,A)\rtimes_rG @>>> C(G,A/I)\rtimes_rG @>>> 0,
\end{CD}
$$}
 is exact, where the vertical maps are induced by the canonical inclusions of 
$I$, $A$, and $A/I$ into $C(G/H,I), C(G/H,A)$ and $C(G/H,A/I)$, respectively.
Since these inclusions are injective, all vertical maps are injective, too
(see the remarks preceeding Proposition \ref{prop-exactfull}). This and 
the exactness of the lower horizontal row imply
that
$$\ker\big(A\rtimes_rG\to (A/I)\rtimes_rG\big)=:J=\big(A\rtimes_rG\big)\cap\big(C(G/H,I)\rtimes_rG\big).$$
Let $(x_i)_i$ be a bounded approximate unit of $I$ and 
let $(\varphi_j)_j$ be an approximate unit of $C_c(G)$ 
(compare with Remark \ref{rem-approx}).
Then $z_{i,j}:=\varphi_j\otimes x_i\in C_c(G,I)$ serves as an approximate unit of 
$I\rtimes_rG$ and of $J:=\big(A\rtimes_rG\big)\cap\big(C(G/H,I)\rtimes_rG\big)$.
Thus if $y\in J$, then $z_{i,j}\cdot y\in I\rtimes_rG$ and  $z_{i,j}\cdot y$ converges to $y$.
Hence $J\subseteq I\rtimes_rG$. 
\end{proof}

\begin{corollary}\label{cor-subgroup}
Every closed subgroup of an almost connected group is exact (in particular, 
every free group in countably many generators is exact). Also, every closed subgroup 
of $\GL(n,\QQ_p)$, where $\QQ_p$ denotes the field of $p$-adic rational numbers equipped 
with the Hausdorff topology is exact.
\end{corollary}
\begin{proof} Recall first that a locally compact group is called {\em almost connected}
if the component $G_0$ of the identity in $G$ is cocompact.
By part (i) of Theorem \ref{thm-KW1} it is enough to show that
every almost connected group $G$ is exact and that $\GL(n,\QQ_p)$ is exact for all $n\in \NN$.
But structure theory for those groups implies that in both cases one can find 
an amenable cocompact subgroup. Since amenable groups are exact (by Propositions 
\ref{prop-amenable} and \ref{prop-exactfull}),
 the result then follows from part (ii) of the theorem.
\end{proof}

\begin{remark}\label{rem-exact} 
We should mention that Kirchberg and Wassermann proved some further 
permanence results: If $H$ is a closed subgroup of $G$ such that $G/H$ 
carries a finite invariant  measure, then $H$ exact implies $G$ exact.
Another important result is the extension result: If $N$ is a closed normal subgroup of
$G$ such that $N$ and $G/N$ are exact, then $G$ is exact. The proof of this result
needs the notion of twisted actions and twisted crossed products. We shall present 
that theory and the proof of the extension result for exact groups in \S \ref{sec-twistmor} below.
We should also mention that the proof of part (ii) of Theorem \ref{thm-KW1}, and 
hence of Corollary \ref{cor-subgroup} followed
some ideas of Skandalis (see also the discussion at the end of \cite{KW1}).

By work of Ozawa and others (e.g. see \cite{Oz} for a general discussion),
the class of discrete exact groups is known to be 
 identical to the class of all discrete groups
which can act amenably on some compact Hausdorff space $X$ (we refer to \cite{ADR}
for a quite complete exposition of amenable actions). An analogous result for 
general second countable locally compact groups has been announced very recently 
by Brodzki, Cave and Li in \cite{BCL}. This implies a new proof that exactness 
passes to closed subgroups, since the restriction of an amenable action to a closed 
subgroup is amenable.
If we apply the exactness condition of a group $G$ to trivial actions, it follows 
from Remark \ref{rem-fullgroup} that $C_r^*(G)$ is an exact $C^*$-algebra 
if $G$ is exact --- the converse is known for discrete groups by \cite{KW1} but 
is still open in the general case. As mentioned at the end of \S\ref{sec-cvt}, it is now
known that there exist non-exact finitely generated discrete groups.
\end{remark}

\section{Induced representations and the ideal structure of crossed products}\label{sec-Mackey}
In this section we  use Green's imprimitivity theorem as a basis for computing the representation 
theory and/or the ideal structure of C*-group algebras and crossed products. First ideas towards this theory appeared 
in the work of Frobenius and Schur on representations of finite groups. In the 1940's George Mackey  
introduced the theory of induced representations of second countable locally compact groups 
together with a procedure (now known as the Mackey machine)
 to compute the irreducible representations of a second countable locally compact group
$G$  in terms of representations 
of a (nice!) normal subgroup $N$ and projective representations of the stabilisers for the action of $G$ on $\widehat{N}$
(see \cites{Ma1,Ma2,Ma3, Ma4, Ma5}.)
For most of the theory, the separability assumption has been eliminated by Blattner in \cites{Blat1, Blat2}.
An extension of this theory to crossed products has first been worked out by Takesaki in \cite{T1}.
In the 1970's Marc Rieffel first showed that the theory of induced representations of groups can be embedded 
into his more algebraic theory of induced representations of C*-algebras as introduced in Section \ref{subsec-ind-rep} (\cites{Rie1,Rie2}).
The full power of this theory became evident with  the fundamental work of Phil Green on twisted crossed products (\cite{Green1}).
In what follows we will try to explain the basics of Green's theory by first restricting to ordinary crossed products. 
The twisted crossed products will be studied later in Section \ref{sec-twists}.
We will also report on the important work of Sauvageot, Gootmann, and Rosenberg (\cites{S1, S2, GR}) 
on the generalised Effros-Hahn conjecture, i.e., 
on the ideal structure of (twisted) crossed products $A\rtimes_\alpha G$ in which the action of $G$ on $\Prim(A)$ 
does not have very good properties.

Many of the results explained in this section also carry over to groupoids and to crossed products 
by (twisted) actions of groupoids on C*-algebras (e.g., see \cites{Ren1, Ren2, IW}), but we shall stick to 
(twisted) crossed products by group actions in these notes.

\subsection{Induced representations of groups and crossed products}\label{subsec-indrep}
If $(A,G,\alpha)$ is a system and $H$ is a closed subgroup of $G$, 
then
Green's imprimitivity theorem provides an imprimitivity bimodule 
$X_H^G(A)$ between $C_0(G/H,A)\rtimes G\cong \Ind_H^GA\rtimes G$ 
and $A\rtimes H$. In particular, $C_0(G/H,A)\rtimes G$ identifies with 
the compact operators $ \K(X_H^G(A))$ on $X_H^G(A)$. There is a  canonical covariant homomorphism
$$(k_A, k_G):(A,G)\to M(C_0(G/H,A)\rtimes G)\cong \L(X_H^G(A)),$$
where  $k_A=i_{C_0(G/H,A)}\circ j_A$ denotes the composition of the  inclusion 
$j_A:A\to M\big(C_0(G/H,A)\big)$ with the inclusion 
$i_{C_0(G/H,A)}:C_0(G/H,A)\to 
M\big(C_0(G/H,A)\rtimes G\big)$ and $k_G$ denotes the canonical inclusion of $G$
into $M\big(C_0(G/H,A)\rtimes G\big)$.
The  integrated 
form 
$$k_A\times k_G:A\rtimes G\to M\big(C_0(G/H,A)\rtimes G\big)\cong \L(X_H^G(A))$$
determines a left action of $A\rtimes G$ on $X_H^G(A)$ and 
we obtain a canonical element 
$[X_H^G(A), k_A\times k_G]\in \Mor(A\rtimes G, A\rtimes H)$ --- a morphism from 
$A\rtimes G$ to $A\rtimes H$ in the correspondence category. Using the techniques of 
\S \ref{subsec-ind-rep}, we can define induced representations of $A\rtimes G$ as follows:

\begin{definition}\label{def-indrep}
For $\rho\times V\in \Rep(A\rtimes H)$ we define the induced representation 
$\ind_H^G(\rho\times V)\in \Rep(A\rtimes G)$ as the representation induced 
from $\rho\times V$ via
$[X_H^G(A), k_A\times k_G]\in \Mor(A\rtimes G, A\rtimes H)$.

Similarly, for $J\in \I(A\rtimes H)$, we define the induced ideal $\ind_H^GJ\in \I(A\rtimes G)$
as the ideal induced from $J$ via  $[X_H^G(A), k_A\times k_G]$.
\end{definition}

On the other hand, if we restrict the canonical embedding $i_G:G\to M(A\rtimes G)$ to $H$,
we obtain a non-degenerate homomorphism $i_A\times i_G|_H:A\rtimes H\to M(A\rtimes G)$
which induces a morphism $[A\rtimes G, i_A\times i_G|_H]\in \Mor(A\rtimes H, A\rtimes G)$.
This leads to

\begin{definition}\label{def-resrep}
For $\pi\times U\in \Rep(A\rtimes G)$ we define the {\em restriction} 
$\res_H^G(\pi\times U)\in \Rep(A\rtimes H)$ as the representation induced 
from $\pi\times U$ via
$[A\rtimes G, i_A\times i_G|_H]\in \Mor(A\rtimes H, A\rtimes G)$.

Similarly, for $I\in \I(A\rtimes G)$, we define the restricted ideal $\res_H^GI\in \I(A\rtimes H)$
as the ideal induced from $I$ via  $[A\rtimes G,i _A\times i_G|_H]$.
\end{definition}

\begin{remark}\label{rem-res} It is a good exercise to show that for any 
$\pi\times U\in \Rep(A\rtimes G)$ we have $\res_H^G(\pi\times U)=\pi\times U|_H$ --- 
the integrated form of the restriction of $(\pi,U)$ to $(A,H,\alpha)$.
\end{remark}

As a consequence of Definitions \ref{def-indrep} and \ref{def-resrep}
 and Proposition \ref{prop-cont} we get

\begin{proposition}\label{prop-contHG}
The maps $\ind_H^G:\Rep(A\rtimes H)\to \Rep(A\rtimes G)$ 
and $\ind_H^G:\I(A\rtimes H)\to\I(A\rtimes G)$
as well as the maps 
$\res_H^G:\Rep(A\rtimes G)\to\Rep(A\rtimes H)$ and 
$\res_H^G:\I(A\rtimes G)\to \I(A\rtimes H)$  are continuous with respect to the 
Fell topologies.
\end{proposition}

\begin{remark}\label{rem-induce}
{\bf (1)} Note that the left action of $A\rtimes G$ on $X_H^G(A)$ can be 
described conveniently on the level of $C_c(G,A)$ via convolution: If $f\in C_c(G,A)\subseteq A\rtimes G$ and 
$\xi\in C_c(G,A)\subseteq X_H^G(A)$, then 
$k_A\times k_G(f)\xi=f*\xi$. 

{\bf (2)} For $A=\CC$ we obtain, after identifying unitary representations 
of $G$ (resp. $H$) with $*$-representations of $C^*(G)$ (resp. $C^*(H)$)
an induction map $\ind_H^G:\Rep(H)\to \Rep(G)$. 
With a bit of work one can check that $\ind_H^GU$ for $U\in \Rep(H)$
coincides (up to equivalence) with the induced representation
defined by Mackey in \cite{Ma1} or Blattner in \cite{Blat1}.
Similarly, the induced representations for $C^*$-dynamical systems
as defined above coincide up to equivalence 
with the induced representations as  constructed  by
Takesaki in \cite{T1}. We will present some more details on these facts in 
Proposition \ref{prop-blattind} and Corollary \ref{cor-indblatt} below.

{\bf (3)} If $\rho$ is a non-degenerate representation of $A$ on a Hilbert space
$H_{\rho}$, then 
$\ind_{\{e\}}^G\rho$ is equivalent to the regular representation $\Ind\rho$
of $A\rtimes G$ on $L^2(G,H_{\rho})$ (see Remark \ref{rem-reg}). The intertwining unitary $V:X_H^G(A)\otimes_A H_{\rho}\to 
L^2(G,H_{\rho})$ is given by 
$$\big(V(\xi\otimes v)\big)(s)=\rho\big(\alpha_{s^{-1}}(\xi(s))\big)v$$
for $\xi\in C_c(G,A)\subseteq X_H^G(A)$ and $v\in H_{\rho}$.
%
%{\bf (4)} Using twisted systems $(A,G,N,\alpha,\tau)$
% and their twisted crossed products everywhere above
%provides an induction map $\ind_H^G:\Rep(A\rtimes (H,N))\to \Rep(A\rtimes(G,N))$
%whenever we start with a closed subgroup $H$ of $G$ which contains the 
%domain $N$ of the twist $\tau$. Indeed, one easily checks that 
%in this case the covariant homomorphism 
%$(k_A,k_G): (A,G)\to M\big(C_0(G/H,A)\rtimes (G,N)\big)$ preserves the twist $\tau$ and hence
%integrates to a $*$-homomorphism of the twisted crossed product $A\rtimes (G,N)$
%into $M\big(C_0(G/H,A)\rtimes (G,N)\big)\cong\L(Y_H^G(A))$, where $Y_H^G(A)$
%is the $C_0(G/H,A)\rtimes (G,N)-A\rtimes (H,N)$ imprimitivity bimodule as described
%in \S\ref{sec-packrae}.

%{\bf (5)} All of the 
% general properties of induced representations stated below have complete 
%analogues in the twisted case. In most cases one can deduce these 
% facts from the ordinary case by using 
%Theorem \ref{thm-twistedequiv} together with compatibility 
%of induction with Morita equivalence
%(e.g.\ see \cite{E-mor, EKQR1}).
\end{remark}

The construction of $[X_H^G(A), k_A\times k_G]$
shows that we have a decomposition 
$$[X_H^G(A), k_A\times k_G]=[C_0(G/H,A)\rtimes G, k_A\times k_G]\circ [X_H^G(A)]
$$
as morphisms in the correspondence category.
Hence the induction map $\ind_H^G:\Rep(A\rtimes H)\to \Rep(A\rtimes G)$ factors
as the composition 
$$
\begin{CD}
\Rep(A\rtimes H)  @>\Ind^{X_H^G(A)}>\cong> \Rep\big(C_0(G/H,A)\rtimes G\big)
@>(k_A\times k_G)^*>> \Rep(A\rtimes G)
\end{CD}
$$
(see Remark \ref{rem-ind} for the meaning of $(k_A\times k_G)^*$).
The representations of $C_0(G/H,A)\rtimes G$ are of the form 
$(P\otimes \pi)\times U$, where $P$ and $\pi$ are commuting representations 
of $C_0(G/H)$ and $A$, respectively  (we use the identification $C_0(G/H,A)\cong C_0(G/H)\otimes A$).
The covariance condition for $(P\otimes \pi,U)$ is equivalent  
to $(\pi,U)$ and $(P,U)$ being covariant representations 
of $(A,G,\alpha)$ and $(C_0(G/H),G,l)$, respectively (where $l:G\to \Aut(C_0(G/H))$ is the left translation action).
One then checks that 
$$\big((P\otimes \pi)\times U\big)\circ (k_A\times k_G)=\pi\times U.$$
Since induction from $\Rep(A\rtimes H)$ to $\Rep\big(C_0(G/H,A)\rtimes G\big)$
via $X_H^G(A)$ is a bijection, we obtain the following 
general version of Mackey's classical imprimitivity theorem for group
representations (see \cite{Ma1} and \cite{T1}):

\begin{theorem}[Mackey-Takesaki-Rieffel-Green]\label{thm-imp}
Suppose that $(A,G,\alpha)$ is a system and let $H$ be a closed 
subgroup of $G$. Then:
\begin{enumerate}
\item A representation $\pi\times U\in \Rep(A\rtimes G)$ on a Hilbert space
$H_{\pi}$ 
is induced from a representation $\sigma\times V\in \Rep(A\rtimes H)$
if and only if there exists a non-degenerate representation 
$P:C_0(G/H)\to \B(H_{\pi})$ which commutes with $\pi$  
such that $(P,U)$ is a covariant representation of 
$(C_0(G/H),G,l)$.
\item If $\pi\times U\in \Rep(A\rtimes G)$ is induced from the irreducible 
representation $\sigma\times V\in \Rep(A\rtimes H)$, and if 
$P:C_0(G/H)\to \B(H_{\pi})$ is the corresponding representation 
such that $(P\otimes \pi)\times U\cong \Ind^{X_H^G(A)}(\sigma\times V)$,
then $\pi\times U$ is irreducible if and only if every
$W\in \B(H_{\pi})$ which intertwines with $\pi$ and $U$ (and hence 
with $\pi\times U$) also intertwines with $P$.
\end{enumerate} 
\end{theorem}

\begin{proof}
The first assertion follows directly from the above discussions. The second 
statement follows from Schur's irreducibilty criterion (a representation is irreducible
iff every intertwiner is a multiple of the identity) together with the fact that 
induction via imprimitivity bimodules preserves irreducibility of representations 
in both directions (see Proposition \ref{prop-induce}).
\end{proof}

In many situations it is convenient to have a more concrete realization 
of the induced representations. The following construction follows 
Blattner's construction of induced group representations of groups (see \cites{Blat1, Fol}). 
We start with the situation of an induced system:
 Assume that $H$ is a closed subgroup of $G$ and that $\alpha:H\to \Aut(A)$ is an 
$H$-action. 
%For $h\in H$ let  $\mu_H(h):=\big(\Delta_H(h)/\Delta_G(h)\big)^{\frac{1}{2}}$.
If $\rho\times V\in \Rep(A\rtimes H)$ is a representation on the Hilbert space 
$H_\rho$ we put 
$$\F_{\rho\times V}:=\left\{\xi:G\to H_{\rho}: \begin{matrix}
\xi(sh)=\sqrt{\Delta_H(h)/\Delta_G(h)}V_{h^{-1}}\xi(s)\;\text{for all $s\in G, h\in H$}\\
\text{and $\xi$ is continuous  with compact support modulo $H$}\end{matrix}\right\}.
$$
Let $c:G\to [0,\infty)$ be a Bruhat section for $H$, i.e., $c$ is continuous with
$\supp c\cap C\cdot H$ compact for all compact $C\subseteq G$ and such that
$\int_H c(sh)\, dh=1\quad\text{for all}\; s\in G$ (for the existence of such $c$ see
\cite{BourT}). Then
%

%
% the surjective linear map given by 
%$$\tau_H(\varphi)(sH):=\int_H \varphi(sh)\,dh.$$
%For $\xi,\eta\in \F_{\rho\times V}$ we can use 
%Urysohn's theorem to find a function $\varphi_{\xi,\eta}\in C_c(G)$ such that 
%$\tau_H(\varphi_{\xi,\eta})\equiv 1$ on $\big(\supp(\xi)\cup\supp(\eta)\big)/H$.

$$\lk \xi,\eta\rk:=\int_G c(s)\lk \xi(s),\eta(s)\rk\,ds$$
determines a well defined  inner product on $\F_{\rho\times V}$
and we let $H_{\ind(\rho\times V)}$ denote its Hilbert space
completion. We can now define representations 
$\sigma$ and $U$ of $\Ind_H^GA$ and $G$ 
on $H_{\ind(\rho\times V)}$, respectively,  by
\begin{equation}\label{eq-ind1}
(\sigma(f)\xi)(s):=\rho(f(s))\xi(s)\quad \text{and}\quad
(U_t\xi)(s):=\xi(t^{-1}s).
\end{equation}
Then $\sigma\times U$ is a representation of $\Ind_H^GA\rtimes G$
on $H_{\ind(\rho\times V)}$ and a straightforward but lengthy computation
gives:

\begin{proposition}\label{prop-blattind}
Let $X:=X_H^G(A)$ denote  Green's $\Ind_H^GA\rtimes G-A\rtimes H$ imprimitivity bimodule
and let $\rho\times V$ be a representation of $A\rtimes H$ on $H_{\rho}$. Then 
there is a unitary $W:X\otimes_{A\rtimes H} H_{\rho}\to 
H_{\ind(\rho\times V)}$, given 
on elementary tensors $x\otimes v\in X\odot H_\rho$ by
$$W(x\otimes v)(s)=\Delta_G(s)^{-\frac{1}{2}}\int_H\Delta_H(h)^{-\frac{1}{2}}V_h\rho(x(sh))v\,dh,$$
which implements a unitary equivalence between $\Ind^X(\rho\times V)$ 
and the representation $\sigma\times U$ defined above.
\end{proposition}

Observe that in case where $H$ is open in $G$, the representation $(\sigma, U)$ constructed above 
coincides with the representation $(\Ind\rho,\Ind V)$ as constructed in the proof of 
Proposition \ref{prop-corner}.

In the special case where $A$ is a $G$-algebra we identify 
$\Ind_H^GA$ with $C_0(G/H,A)$ via the isomorphism
$\Phi$ of Remark \ref{rem-Ind} (1). It is then easy to check that 
the representation $\sigma$ defined above corresponds to the 
representation $P\otimes \pi$ of $C_0(G/H,A)\cong C_0(G/H)\otimes A$ on 
$H_{\ind(\rho\times V)}$ given by the formula
\begin{equation}\label{eq-ind2}
(P(\varphi)\xi)(s):=\varphi(sH)\xi(s)\quad\text{and}\quad (\pi(a)\xi)(s)=\rho(\alpha_{s^{-1}}(a))\xi(s).
\end{equation}
Hence, as a direct corollary of the above proposition we get:

\begin{corollary}\label{cor-indblatt}
Let $(A,G,\alpha)$ be a system and let $\rho\times V\in \Rep(A\rtimes H)$ for
some closed subgroup $H$ of $G$. Then $\ind_H^G(\rho\times V)$, as defined in Definition \ref{def-indrep},
is unitarily equivalent to the representation $\pi\times U$ of $A\rtimes G$ 
on $H_{\ind(\rho\times V)}$
with $\pi$ and $U$ as in Equations (\ref{eq-ind2}) and (\ref{eq-ind1}), respectively.
\end{corollary}

Another corollary which we can easily obtain from Blattner's realisation 
is the following useful observation: Assume that $H$ is a closed subgroup of $G$ 
and that $A$ is an $H$-algebra. Let $\eps_e:\Ind_H^GA\to A$ be the 
$H$-equivariant surjection defined by evaluation of  functions $f\in \Ind_H^GA$ at the 
unit $e\in G$.
If $\rho\times V\in \Rep(A\rtimes H)$ then  
$(\rho\circ \eps_e) \times V$ is a representation 
of $\Ind_H^GA\rtimes H$. We then get

\begin{corollary}\label{cor-inducedalg}
The induced representation  $\ind_H^G\big((\rho\circ \eps_e)\times V\big)$
(induction from $H$ to $G$ for the system $(\Ind_H^GA,G,\Ind\alpha)$)
is unitarily equivalent to $\ind^{X_H^G(A)}(\rho\times V)$ (induction via Green's 
$\Ind_H^GA\rtimes G-A\rtimes H$ imprimitivity bimodule $X_H^G(A)$).
\end{corollary}
\begin{proof} By Proposition \ref{prop-blattind} and Corollary \ref{cor-indblatt},
both representations can be realized on the Hilbert space $H_{\ind(\rho\times V)}$
whose construction only depends on $G$ and the unitary representation $V$ of $H$.
Applying the formula for $\pi$ in (\ref{eq-ind2}) to the present situation, we 
see that the $\Ind_H^GA$-part  of $\ind_H^G\big((\rho\circ \eps_e)\times V\big)$ is given by the 
formula
$$(\pi(f)\xi)(s)=\rho\big(\ind\alpha_{s^{-1}}(f)(e)\big)\xi(s)=\rho(f(s))\xi(s)=(\sigma(f)\xi)(s)$$
with $\sigma$ as in (\ref{eq-ind1}).
\end{proof}

We now turn to some further properties of induced representations. To obtain those properties
we shall pass from
 Green's to Blattner's realizations of the induced representations and back 
whenever it seems  convenient.
We start the discussion with the
theorem of induction in steps. For this suppose that $L\subseteq H$ 
are closed subgroups of $G$. To avoid confusion, we write
$\Phi_H^G$ for the left action of $A\rtimes G$ on $X_H^G(A)$ (i.e.,
$\Phi_H^G=k_A\times k_G$ in the notation used above) and we write
$\Phi_L^G$ and $\Phi_L^H$ for the left actions 
of $A\rtimes G$ and $A\rtimes H$ on $X_L^G(A)$ and $X_L^H(A)$,
respectively. Then the theorem of induction in steps reads as

\begin{theorem}[Green]\label{thm-steps}
Let $(A,G,\alpha)$ and $L\subseteq H$  be as above. 
Then  
$$[X_H^G(A),\Phi_H^G]\circ [X_L^H(A),\Phi_L^H] 
=[X_L^G(A),\Phi_L^G]$$
 as morphisms from $A\rtimes G$ to 
$A\rtimes L$ in the correspondence category $\Corr$. As a consequence,
we have 
$$\ind_H^G\big(\ind_L^H(\rho\times V)\big)=\ind_L^G(\rho\times V)$$
for all $\rho\times V\in \Rep(A\rtimes L)$.
\end{theorem}
\begin{proof} For the proof one has to check that
$X_H^G\otimes_{A\rtimes H} X_L^H\cong X_L^G(A)$ as 
Hilbert $A\rtimes G-A\rtimes L$ bimodule. Indeed, one can check that such 
isomorphism is given on the level functions by the pairing
$C_c(G,A)\otimes C_c(H,A)\to C_c(G,A)$
as given by the second formula in (\ref{eq-products}).
We refer to \cite{Green1} and \cite{Wi-crossed}*{Theorem 5.9} for more details.
\end{proof}

By an {\em automorphism} $\gamma$ of a system $(A,G,\alpha)$ we
understand a pair $\gamma=(\gamma_A,\gamma_G)$, where $\gamma_A$
is a $*$-automorphism of $A$ and $\gamma_G:G\to G$ is an automorphism of $G$
such that $\alpha_{\gamma_G(t)}=\gamma_A\circ \alpha_t\circ \gamma_A^{-1}$
for all $t\in G$. 
An {\em inner automorphism} of $(A,G,\alpha)$ is an automorphism of the form
 $(\alpha_s, C_s)$, $s\in G$,  with $C_s(t)=sts^{-1}$.
If $\gamma=(\gamma_A,\gamma_G)$ is an automorphism of $(A,G,\alpha)$ 
and if $H$ is a closed subgroup of $G$, then $\gamma$ induces an isomorphism
$\gamma_{A\rtimes H}:A\rtimes H\to A\rtimes H_\gamma$ with $H_\gamma:=\gamma_G(H)$ via 
$$\gamma_{A\rtimes H}(f)(h):=\gamma_A\big(f(\gamma_G^{-1}(h))\big)\quad \text{for $h\in H_{\gamma}$
and $f\in C_c(H,A)$,}$$
where we adjust Haar measures on $H$ and $H_\gamma$ such that 
$\int_{H} f(\gamma_G(h))\,dh
=\int_{H_\gamma} f(h')\,dh'$ for $f\in C_c(H_\gamma)$. 
Note that if $(\rho,V)\in \Rep(A,H_{\gamma})$, then $(\rho\circ \gamma_A, V\circ \gamma_G)\in \Rep(A,H)$ 
and we have
$$(\rho\times V)\circ \gamma_{A\rtimes H}\cong (\rho\circ \gamma_A)\times(V\circ \gamma_G)$$
for their integrated forms.
\begin{remark}\label{rem-inner}
If $H=N$ is normal in $G$ and if $\gamma_s=(\alpha_s, C_s)$ is an inner 
automorphism of $(A,G,\alpha)$, then we will write $\alpha^N_s$ for the 
corresponding automorphism of $A\rtimes N$. Then $s\mapsto \alpha^N_s$ is an
action of $G$ on $A\rtimes N$. This action will serve as a starting point for the 
study of twisted actions in \S \ref{sec-twists} below.
\end{remark}

\begin{proposition}\label{prop-auto}
Suppose that $\gamma=(\gamma_A,\gamma_G)$ is an automorphism of $(A,G,\alpha)$
and let $H\subseteq L$ be two closed subgroups of $G$. Then 
$$\ind_{H}^{L}\big((\rho\times V)\circ \gamma_{A\rtimes H} \big)\cong 
\big(\ind_{H_\gamma}^{L_\gamma}(\rho\times V)\big)\circ \gamma_{A\rtimes L}$$
for all $\rho\times V\in \Rep(A\rtimes H_{\gamma})$, where ``$\cong$'' denotes unitary equivalence.
In particular, if $\rho\times V\in \Rep(A,H)$ and  $(\alpha_s, C_s)$ is an inner automorphism of 
$(A,G,\alpha)$ 
then 
$$\ind_H^G(\rho\times V)\cong  \ind_{sHs^{-1}}^G\big(s\cdot(\rho\times V)\big),$$
where we put $s\cdot(\rho\times V):=(\rho\circ\alpha_{s^{-1}})\times 
(V\circ C_{s^{-1}})\in \Rep(A,sHs^{-1})$.
\end{proposition}
\begin{proof} Simply check that the map $\gamma_{A\rtimes L}:C_c(L,A)\to C_c(L_{\gamma},A)$ as defined above
also extends to a bijection $\Phi_L:X_H^L(A)\to  X_{H_{\gamma}}^{L_\gamma}(A)$ 
which is compatible with the isomorphisms $\gamma_{A\rtimes L}:A\rtimes L\to A\rtimes L_\gamma$ 
and $\gamma_{A\rtimes H}:A\rtimes H\to A\rtimes H_\gamma$ on the left and right. This implies that
$\gamma_{A\rtimes L}^*\circ [X_{H_\gamma}^{L_\gamma}(A), k_A\times k_{L_\gamma}]=
[X_H^L(A), k_A\times k_L]\circ\gamma_{A\rtimes H}^*$ in $\Corr$ and the first statement follows.
The second statement follows from the first applied to $L=G$ and $\gamma=(\alpha_s,C_s)$
together with the fact that for any $\pi\times U\in\Rep(A,G)$ the unitary $U_s\in \U(H_{\pi})$
implements a unitary equivalence between 
$s\cdot(\pi\times U)=(\pi\circ \alpha_{s^{-1}})\times (U\circ C_{s^{-1}})$ and
 $\pi\times U$.
\end{proof}

As a direct consequence we get:

\begin{corollary}\label{cor-regular}
Let $(A,G,\alpha)$ be a system. For $J\in \I(A)$ let 
$$J^G:=\cap\{\alpha_s(J):s\in G\}.$$
 Then 
$\ind_{\{e\}}^G J^G=\ind_{\{e\}}^G J$ in $A\rtimes G$.
As a consequence, if $\rho\in \Rep(A)$ 
 such that
$\cap\{\ker(\rho\circ \alpha_s): s\in G\}=\{0\}$, then 
$\ind_{\{e\}}^G\rho$ factors through a faithful representation of 
the reduced crossed product $A\rtimes_rG$.
\end{corollary}
\begin{proof} Let $J=\ker\rho$ for some $\rho\in \Rep(A)$ and let
 $\rho^G:=\oplus_{s\in G}\rho\circ \alpha_s$. Then $J^G=\ker \rho^G$.
 It follows from  Proposition \ref{prop-auto}
 that $\ind_{\{e\}}^G\rho\circ \alpha_s\cong \ind_{\{e\}}^G\rho$ for all $s\in G$.
 Since induction preserves direct sums, it follows that 
 $$\ind_{\{e\}}^GJ=\ker(\ind_{\{e\}}^G\rho)=\ker(\ind_{\{e\}}^G\rho^G)=\ind_{\{e\}}^GJ^G.$$
 If $\cap\{\ker(\rho\circ \alpha_s): s\in G\}=\{0\}$, then $\rho^G$ is faithful
 and  it follows from 
Remark \ref{rem-universal} (3) and Remark \ref{rem-induce} (3) that 
$\ker\Lambda_A^G=\ker(\ind_{\{e\}}^G\rho^G)=\ker(\ind_{\{e\}}^G\rho)$.
\end{proof}

\begin{remark}\label{rem-reducedimp}
From the previous results it is now possible  to obtain a fairly easy proof 
of the fact that Green's $\Ind_H^GA\rtimes G-A\rtimes H$ imprimitivity bimodule $X_H^G(A)$ factors to give a $\Ind_H^GA\rtimes_r G-A\rtimes_r H$ imprimitivity 
bimodule for the reduced crossed products (compare with Remark \ref{rem-Green}).
Indeed, if $\rho$ is any faithful representation of $A$, and
if $\eps_e:\Ind A\to A$ denotes evaluation at the unit $e$,
it follows from Corollary \ref{cor-regular}, that 
$\ker(\Lambda_{\Ind A}^G)=\ker\big(\ind_{\{e\}}^G(\rho\circ\eps_e)\big)$.
The latter coincides with
$\ker(\ind_H^G(\ind_{\{e\}}^H(\rho\circ\eps_e)))$ by Theorem \ref{thm-steps}.
If $\sigma\times V$ denotes the representation $\ind_{\{e\}}^H\rho\in \Rep(A\rtimes H)$,
then $\ker(\sigma\times V)=\ker\Lambda_A^H$ since 
$\rho$ is faithful on $A$ and
a short computation shows that $\ind_{\{e\}}^H(\rho\circ\eps_e)=
(\sigma\circ \eps_e)\times V$, where on the left-hand side
we use induction in the system $(\Ind A,H,\Ind \alpha)$.
Putting all this together we get
 \begin{align*}
 \ker(\Lambda_{\Ind A}^G)&=\ker\big(\ind_{\{e\}}^G(\rho\circ\eps_e)\big)
 =\ker\big(\ind_H^G(\ind_{\{e\}}^H(\rho\circ\eps_e))\big)\\
& =\ker\big(\ind_H^G\big((\sigma\circ \eps_e)\times V)\big)
\stackrel{*}{=}\ker\big(\ind^{X_H^G(A)}(\sigma\times V)\big)\\
&=\ker\big(\ind^{X_H^G(A)}\Lambda_A^H\big)\stackrel{**}{=}\ind^{X_H^G(A)}\big(\ker\Lambda_A^H\big),
\end{align*}
where * follows from Corollary \ref{cor-inducedalg} and ** follows from Equation (\ref{eq-kerindpi}).
The desired result then follows from the Rieffel correspondence (Proposition \ref{prop-rief-cor}).
\end{remark}

We now come to some important results concerning the relation between induction 
and restriction of representations and ideals (see Definition \ref{def-resrep} for 
the definition of the restriction maps). We start with
%
%\begin{definition}\label{def-res}
%If $(A,G,\alpha)$ is a system and 
%$H$ is a closed subgroup of $G$, we may restrict a representation 
%$\pi\times U\in \Rep(A\rtimes G)$
%to the representation $\pi\times U|_H\in \Rep(A\rtimes H)$
%and we obtain a {\em restriction map}
%$$\res_H^G:\Rep(A\rtimes G)\to \Rep(A\rtimes H); \;\;\res_H^G(\pi\times U):=\pi\times U|_H.$$
%Similarly, we obtain a restriction map
%$$\res_H^G:\I(A\rtimes G)\to\I(A\rtimes H); I\mapsto \res_H^GI$$
%by defining $\res_H^GI=\ker(\pi\times U|_H)$ whenever $I=\ker(\pi\times U)$ for 
%$\pi\times U\in \Rep(A\rtimes G)$.
%\end{definition}

%If $(i_A,i_G):(A,G)\to M(A\rtimes G)$ denote the canonical maps, then 
%the restriction $(i_A, i_G|_H)$ to $(A,H)$ is covariant and integrates to a non-degenerate 
%$*$-homomorphism 
%$$i_A\times i_G|_H:A\rtimes H\to M(A\rtimes G)$$
%such that $\res_H^G(\pi\times U)=(i_A\times i_G|_H)^*(\pi\times U)$ for all 
%$\pi\times U\in \Rep(A\rtimes G)$. Hence, the restriction map is given via composition 
%by a morphism in the Morita category $\Corr$. 
%In particular, it follows from equation (\ref{eq-kerindpi}) that the restriction map
%$\res_H^G:\I(A\rtimes G)\to \I(A\rtimes H)$ is well defined.
%Moreover, Proposition \ref{prop-cont} implies 

%\begin{proposition}\label{prop-contres}
%The restriction  maps 
%$$\res_H^G:\Rep(A\rtimes G)\to \Rep(A\rtimes H)\quad\text{and}\quad
%\res_H^G:\I(A\rtimes G)\to \I(A\rtimes H)$$
%are continuous with respect to the Fell topologies.
%\end{proposition}

%The next result is basic for the theory:

\begin{proposition}\label{prop-ind-res}
Suppose that $(A,G,\alpha)$ is a system and let $N\subseteq H$ be closed subgroups of 
$G$ such that $N$ is normal in $G$. Let $\F_{\rho\times V}$ be the dense subspace 
of Blattner's induced Hilbert space $H_{\ind(\rho\times V)}$ as constructed above.
Then 
\begin{equation}\label{eq-indres}
\big(\res_N^G(\ind_H^G(\rho\times V))(f)\xi\big)(s)=\res_N^H(\rho\times V)(\alpha^N_{s^{-1}}(f))\xi(s)
\end{equation}
for all $f\in A\rtimes N$, $\xi\in \F_{\rho\times V}$ and $s\in G$,
where $\alpha^N:G\to \Aut(A\rtimes N)$ is the canonical action of $G$ on 
$A\rtimes N$ (see Remark \ref{rem-inner} above).  As a consequence, if $J\in \I(A\rtimes H)$, we get
\begin{equation}\label{eq-indresideal}\res_N^G\big(\ind_H^GJ\big)=
\cap\{\alpha^N_s(\res_N^H(J)): s\in G\}.
\end{equation}
\end{proposition}
\begin{proof} Define $\sigma:A\rtimes N\to \B(H_{\ind(\rho\times V)})$ by 
$(\sigma(f)\xi)(s)=\rho\times V|_N(\alpha^N_{s^{-1}}(f))\xi(s)$ for $f\in A\rtimes N$ and $\xi\in \F_{\rho\times V}$. Then $\sigma$ is a non-degenerate $*$-representation
and hence it suffices to check that  the left-hand side of (\ref{eq-indres}) coincides with $(\sigma(f)\xi)(s)$
for $f\in C_c(N,A)$. But using (\ref{eq-ind2}) together with the transformation $n\mapsto sns^{-1}$ and the equation $\xi(sn^{-1})=V_n\xi(s)$ for $s\in G, n\in N$, the left-hand side becomes
{\small
\begin{align*}
\big(\res_N^G\big(\ind_H^G(\rho\times V))(f)\xi\big)\big)(s)&=\int_N
\rho(\alpha_{s^{-1}}(f(n)))\xi(n^{-1}s)\,dn\\
&=
\delta(s^{-1})\int_N\rho(\alpha_{s^{-1}}(f(sns^{-1})))\xi(sn^{-1})\,dn\\
&=\int_N\rho(\alpha^N_{s^{-1}}(f)(n))V_n\xi(s)\,dn
=(\sigma(f)\xi)(s).
\end{align*} }
\end{proof}

\begin{remark}\label{compat-ind}
Suppose that $(A,G,\alpha)$ is a system,  $H$ is a closed subgroup of $G$,
and $J\subseteq A$ is a $G$-invariant ideal of $A$. If $\rho\times V$ is a representation 
of $A\rtimes H$ and if we put $\pi\times U:=\ind_H^G(\rho\times V)$, 
then it follows from the above proposition that
$$J\subseteq\ker\rho\quad \Longleftrightarrow\quad
J\subseteq\ker\pi.$$
Hence, we see that the induction map for the system $(A,G,\alpha)$
determines a map 
from $\Rep(A/J\rtimes H)$ to $\Rep(A/J\rtimes G)$ if we identify representations 
of $A/J\rtimes H$ with  the representations $\rho\times V$ of $A\rtimes H$ wich
satisfy $J\subseteq \ker\rho$ (and similarly for $G$).
It is easy to check
(e.g., by using Blattner's construction of the induced representations) that 
this map coincides with the induction map for the system $(A/J,G,\alpha)$.

Also, if $\rho\times V$ is a representation of $A\rtimes H$ such that $\rho$ restricts 
to a non-degenerate representation of $J$, then one can  check that
the restriction of $\ind_H^G(\rho\times V)$ to $J\rtimes G$ conicides
with the induced representation  $\Ind_H^G(\rho|_J\times V)$
where the latter representation is induced from $J\rtimes H$ to 
$J\rtimes G$ via $X_H^G(J)$.\footnote{Of course, these results are also consequences
of the naturality of the assignment $A\mapsto X_H^G(A)$ as stated in Remark \ref{rem-Green} (3).}
We shall use these facts quite frequently below.
\end{remark}

We close this section with some useful results on tensor products of representations.
If $(\pi,U)$ is a covariant representation of the system  $(A,G,\alpha)$ on $H_{\pi}$ and if 
$V$ is a unitary representation of $G$ on $H_V$, then $(\pi\otimes 1_{H_V},U\otimes V)$ 
is a covariant representation of $(A,G,\alpha)$ on $H_\pi\otimes H_V$
and we obtain a pairing
\begin{align*}
\otimes:& \Rep(A\rtimes G)\times \Rep(G)\to \Rep(A\rtimes G);\\
\quad\quad\quad\quad 
&\big( (\pi\times U), V\big)\mapsto (\pi\times U)\otimes V:=
(\pi\otimes 1_{H_V})\times (U\otimes V).
\end{align*}
Identifying $\Rep(G)\cong \Rep(C^*(G))$, this map can also be obtained via the composition
$$\Rep(A\rtimes G)\times \Rep(C^*(G))\to \Rep\big((A\rtimes G)\otimes C^*(G)\big)
\stackrel{\widehat{\alpha}^*}{\to}\Rep(A\rtimes G),$$
where the first map sends a pair $(\pi\times U, V)$ to the external tensor-product representation 
$(\pi\times U)\hat\otimes V$ of $(A\rtimes G)\otimes C^*(G)$ and
 $\widehat{\alpha}:A\rtimes G\to M\big((A\rtimes G)\otimes C^*(G)\big)$ denotes the integrated form
of the tensor product $(i_A\otimes 1_{C^*(G)}, i_G\otimes i_G)$ of the canonical inclusions $(i_A,i_G):(A,G)\to M(A\rtimes G)$ with the 
inclusion  $i_G:G\to M(C^*(G))$ ($\widehat{\alpha}$ is the dual {\em coaction} of $G$ on $A\rtimes G$). Thus,  from Propositions \ref{prop-cont}
and   \ref{prop-tensor} we get

\begin{proposition}\label{prop-diagtensor}
The map $\otimes: \Rep(A\rtimes G)\times \Rep(G)\to \Rep(A\rtimes G)$ 
preserves weak containment in 
both variables and is jointly continuous with respect to the Fell topologies.
\end{proposition}

\begin{proposition}\label{prop-indtensor}
Let $(A,G,\alpha)$ be a system and let $H$ be a closed subgroup of $G$. Then
\begin{enumerate}
\item $\ind_H^G\big((\rho\times V)\otimes U|_H\big)\cong \big(\ind_H^G(\rho\times V)\big)\otimes U$ for all
$\rho\times V\in \Rep(A\rtimes H)$ and $U\in \Rep(G)$;
\item $\ind_H^G\big((\pi\times U|_H)\otimes V\big)\cong (\pi\times U)\otimes \ind_H^GV$
for all $V\in \Rep(H)$ and $\pi\times U\in \Rep(A\rtimes G)$.
\end{enumerate}
In particular, if $\pi\times U\in \Rep(A\rtimes G)$ and $N$ is a normal subgroup of $G$, then
$$\ind_N^G(\pi\times U|_N)\cong (\pi\times U)\otimes \lambda_{G/N},$$
 where $\lambda_{G/N}$ 
denotes the regular representation of $G/N$, viewed as a representation of $G$.
\end{proposition}
\begin{proof} This result can be most easily shown using Blattner's realization 
of the induced representations: In the first case define
$$W:\F_{\rho\times V}\otimes H_U\to \F_{(\rho\times V)\otimes U|_H};\;W(\xi\otimes v)(s)=\xi(s)\otimes U_{s^{-1}}v.$$
Then a short computation shows that $W$ is a unitary intertwiner of
$\big(\ind_H^G(\rho\times V)\big)\otimes U$ and $\ind_H^G\big((\rho\times V)\otimes U|_H)\big)$.
A similar map works for the second equivalence. Since $\lambda_{G/N}=\ind_N^G1_N$,
the last assertion follows from (ii) for the case $V=1_N$.
\end{proof}

\begin{corollary}\label{cor-amenable-tensor}
Suppose that $(A,G,\alpha)$ is a system and that $N$ is a normal subgroup of $G$ such that
$G/N$ is amenable. Then $\pi\times U$ is weakly contained in 
$\ind_N^G(\res_N^G(\pi\times U))$ for all $\pi\times U\in \Rep(A\rtimes G)$.
As a consequence,
$\ind_N^G(\res_N^GI)\subseteq I$ for all $I\in \I(A\rtimes G)$.
\end{corollary}
\begin{proof} Since $G/N$ is amenable if and only if $1_{G/N}\prec \lambda_{G/N}$
we obtain from Proposition \ref{prop-indtensor}
$$\pi\times U=(\pi\times U)\otimes 1_{G/N}\prec (\pi\times U)\otimes \lambda_{G/N}
=\ind_N^G(\pi\times U|_N),$$
which proves the first statement. The second statement follows from the first by choosing
$\pi\times U\in \Rep(A\rtimes G)$ such that $I=\ker(\pi\times U)$.
\end{proof}

\subsection{The ideal structure of crossed products}\label{subsec-MRG}

In this section we come to the main results on the Mackey-Rieffel-Green machine, namely
the description of the spectrum $(A\rtimes G)\dach$ and the primitive ideal space 
$\Prim(A\rtimes G)$ in terms of induced representations (resp.\ ideals) under 
some favorable circumstances. We start with  some topological notations:

\begin{definition}\label{defn-top}
Let $Y$ be a topological space. 
\begin{enumerate}
\item We say that  $Y$ is {\em almost Hausdorff} if 
every nonempty closed subset $F$ of $Y$ contains a nonempty relatively  open Hausdorff subset
$U$ (which can then be chosen to be dense in $F$).
\item A subset $C\subseteq Y$ is called {\em locally closed} if $C$ is relatively open 
in its closure $\overline{C}$, i.e. $\overline{C}\setminus C$ is closed in $Y$.
\end{enumerate}
\end{definition}

It is important to note that if $A$ is a type I algebra, then the spectrum
$\widehat{A}$ (and then also $\Prim(A)\cong\widehat{A}$) is almost Hausdorff
with respect to the Jacobson topology.
This follows from the fact that every quotient of a type I algebra is type I and that 
every nonzero type I algebra contains a nonzero continuous-trace ideal, and hence
its spectrum contains a nonempty Hausdorff subset $U$ 
(see \cite{Dix}*{Chapter 4} and \S \ref{subsec-spectrum}).
Notice also that if $Y$ is almost Hausdoff, then the one-point sets $\{y\}$ are 
locally closed for all $y\in Y$. 
\\

If $A$ is a $C^*$-algebra and if $J\subseteq I$ are two closed two-sided ideals
of $A$, then we may view $\widehat{I/J}$ (resp. $\Prim(I/J)$) 
as a locally closed subset of $\widehat{A}$ (resp. $\Prim(A)$). Indeed, we 
first identify $\widehat{A/J}$ with the closed subset $\{\pi\in \widehat{A}:J\subseteq\ker\pi\}$
of $\widehat{A}$ and then we identify $\widehat{I/J}$ with the open subset $\{\pi\in \widehat{A/J}:
\pi(I)\neq\{0\}\}$ (and similarly for $\Prim(I/J)$ --- compare with \S \ref{subsec-spectrum}). 

Conversely, if $C$ is a locally closed subset of $\widehat{A}$ ,
then $C$ is canonically  homeomorphic to $\widehat{I_C/J_C}$ 
if we take  $J_C:=\ker(C)$ 
and $I_C:=\ker\big(\overline{C}\smallsetminus C\big)$ (we write 
$\ker(E):=\cap\{\ker\pi:\pi\in E\}$ if $E\subseteq \widehat{A}$ and similarly
$\ker(D):=\cap\{P:P\in D\}$ for $D\subseteq\Prim(A)$). 
If we apply this observation to commutative $C^*$-algebras, we recover the well known fact
that the locally closed 
subsets of a locally compact Hausdorff space $Y$ are precisely those subsets of
$Y$ which are locally compact in the relative topology.

\begin{definition}\label{defn-LC}
Let $A$ be a $C^*$-algebra and let $C$ be a locally closed subset of $\widehat{A}$ (resp. 
$\Prim(A)$). Then $A_C:=I_C/J_C$ with $I_C, J_C$ as above is called 
the {\em restriction of $A$ to $C$}. In the same way, we define the restriction $A_D$ of
$A$ to $D$ for a locally closed subset $D$ of $\Prim(A)$.
\end{definition}

In what follows, we shall use the following notations:

\begin{notation}\label{not-orbit}
If $(A,G,\alpha)$ is a system, we consider $\Prim(A)$ as a 
$G$-space via the continuous action 
$G\times\Prim(A)\to\Prim(A);\;(s,P)\mapsto s\cdot P:=\alpha_s(P)$.
We write 
$$G_P:=\{s\in G:s\cdot P=P\}\quad\text{and}\quad
G(P):=\{s\cdot P : s\in G\}$$ 
 for the {\em stabiliser}  and the {\em $G$-orbit} of $P\in\Prim(A)$,
respectively. Moreover, we put 
$$P^G:=\ker G(P)=\cap\{s\cdot P: s\in G\}.$$
Note that the stabilisers $G_P$ are closed subgroups of $G$
for all $P\in \Prim(A)$.\footnote{The fact that $G_P$ is closed in $G$ follows
from the fact that  $\Prim(A)$ is a T$_0$-space. Indeed, if $\{s_i\}$ is a net in ${G_P}$ 
which converges to some $s\in G$, then $P=s_i\cdot P\to s\cdot P$, so $s\cdot P$ is in the 
closure of $\{P\}$. Conversely, we have $s\cdot P=ss_i^{-1}\cdot P\to P$, and hence $\{P\}\in \overline{\{sP\}}$. Since $\Prim(A)$ is a T$_0$-space it follows that $P=s\cdot P$.}
\end{notation}

\begin{remark}\label{rem-GAhat}
Similarly, we may consider the $G$-space $\widehat{A}$ with $G$-action 
$(s,\pi)\mapsto s\cdot \pi:=\pi\circ \alpha_{s^{-1}}$ (identifying representations with their
equivalence classes) and we then write $G_\pi$ and $G(\pi)$ for 
the stabilisers and the $G$-orbits, respectively. However, the stabilisers
$G_{\pi}$ are not necessarily  closed in $G$ if $A$ is not a type I algebra.
If $A$ is type I, then $\pi\mapsto\ker\pi$ is a $G$-equivariant 
homeomorphism from $\widehat{A}$ to $\Prim(A)$.
\end{remark}

The following theorem is due to Glimm \cite{Gl}:

\begin{theorem}[Glimm's theorem]\label{thm-Glimm}
Suppose that $(A,G,\alpha)$ is a separable type I system (i.e., $A$ is a separable
type I algebra and $G$ is second countable). Then the following are equivalent:
\begin{enumerate}
\item The quotient space $G\backslash\!\Prim(A)$ is almost Hausdorff.
\item $G\backslash\!\Prim(A)$ is a T$_0$-space.
\item All points in $G\backslash\!\Prim(A)$ are locally closed.
\item For  all $P\in\Prim(A)$ the quotient $G/G_P$ is homeomorphic to $G(P)$ via
\mbox{$s\cdot G_P\mapsto s\cdot P$}.
\item There exists an ordinal number $\mu$ and an 
increasing sequence $\{I_{\nu} \}_{\nu\leq \mu}$ of $G$-invariant ideals of 
$A$ such that $I_0=\{0\}$, $I_\mu=A$ and $G\backslash\!\Prim({I_{\nu+1}/I_\nu})$
is Hausdorff for all $\nu<\mu$.
\end{enumerate}
\end{theorem}

Hence, if $(A,G,\alpha)$ is a separable type I system satisfying one of the equivalent conditions 
above, then 
$(A,G,\alpha)$ is smooth in the sense of:

\begin{definition}\label{defn-smooth}
A system $(A,G,\alpha)$ is called {\em smooth} if the following two conditions
are satisfied:
\begin{enumerate}
\item  The map $G/G_P\to G(P); \;s\cdot G_P\to s\cdot P$ is a homeomorphism for
all  $P\in \Prim(A)$.
\item The quotient $G\backslash\!\Prim(A)$ is almost Hausdorff, or
$A$ is separable and all orbits $G(P)$ are locally closed 
in $\Prim(A)$.
\end{enumerate}
\end{definition}

If $G(P)$ is a locally closed orbit of $\Prim(A)$, then we may identify 
$G(P)$ with $\Prim(A_{G(P)})$, where $A_{G(P)}=I_{G(P)}/J_{G(P)}$ denotes the restriction 
of $A$ to $G(P)$ as in Definition \ref{defn-LC} 
(note that  $J_{G(P)}=P^G=\cap\{\alpha_s(P): s\in G\}$). Since the ideals $I_{G(P)}$ and $J_{G(P)}$
are $G$-invariant, the action of $G$ on $A$ 
restricts to an action of $G$ on $A_{G(P)}$.
Using exactness of the full crossed-product functor, we get 
\begin{equation}\label{eq-subquot}
A_{G(P)}\rtimes G=(I_{G(P)}/J_{G(P)})\rtimes G\cong (I_{G(P)}\rtimes G)/(J_{G(P)}\rtimes G).
\end{equation}
If $G$ is exact, a similar statement holds for the reduced crossed products.

\begin{proposition}\label{prop-live}
Suppose that $(A,G,\alpha)$ is a system such that 
\begin{enumerate}
\item $G\backslash\!\Prim(A)$ is almost Hausdorff, or 
\item $A$ is separable.
\end{enumerate}
Then, for each $\pi\times U\in (A\rtimes G)\dach$, there exists 
an orbit $G(P)\subseteq \Prim(A)$ such that $\ker\pi=P^G$.
If, in addition, all orbits 
in $\Prim(A)$ are locally closed (which is automatic 
in case of (i)),  then $G(P)$ is uniquely determined by $\pi\times U$.
\end{proposition}
\begin{proof} (Following ideas from \cite{Rie2})
Let $J=\ker\pi$. By passing from $A$ to $A/J$ we may assume without loss 
of generality  that $\pi$ is faithful.
We then have to show that there exists a $P\in \Prim(A)$ such that
$G(P)$ is dense in $\Prim(A)$. 

We first show that under these assumptions
every open subset $W\subseteq G\backslash\!\Prim(A)$ is dense.
Indeed, since $\pi$ is faithful, it follows that 
$\pi\times U$ restricts to a non-zero, and hence irreducible representation 
of $I\rtimes G$, whenever $I$ is a nonzero $G$-invariant ideal of $A$.
In particular, $\pi(I)H_{\pi}=H_{\pi}$ for all such ideals $I$.
Assume now that there are two nonempty $G$-invariant open sets 
$U_1,U_2\subseteq \Prim(A)$ with
$U_1\cap U_2=\emptyset$. Put $I_i=\ker(\Prim(A)\smallsetminus U_i)$, $ i=1,2$.
Then $I_1,I_2$ would be nonzero $G$-invariant ideals such that $I_1\cdot I_2=I_1\cap I_2=\{0\}$,
and  then  
$$H_{\pi}=\pi(I_1)H_{\pi}=\pi(I_1)\big(\pi(I_2)H_{\pi}\big)=\pi(I_1\cdot I_2)H_{\pi}=\{0\},$$
which is a contradiction.

Assume that $G\backslash\!\Prim(A)$ is almost Hausdorff.
If there is no dense orbit $G(P)$ in $\Prim(A)$, then  $G\backslash\!\Prim(A)$ 
contains an open dense Hausdorff subset which contains at least two different points. 
But then there exist nonempty  $G$-invariant open subsets 
$U_1$, $U_2$ of $\Prim(A)$ with $U_1\cap U_2=
\emptyset$, which is impossible.

If $A$ is separable, then $G\setminus\!\Prim(A)$ is second countable (see \cite{Dix}*{Chapter 3})
and we find a countable basis $\{U_n:n\in \NN\}$ for  its topology with $U_n\neq \emptyset$ for all $n\in \NN$. By the first part of this proof we know that all $U_n$ are dense in $G\backslash\Prim(A)$.
Since $G\backslash\!\Prim(A)$ is a Baire space by \cite{Dix}*{Chapter 3}, it follows that 
$D:=\cap_{n\in \NN}U_n$ is also dense in $G\backslash\Prim(A)$. Note that every open subset of 
$G\backslash\!\Prim(A)$ contains $D$.
Hence, if we pick any orbit $G(P)\in D$ then
$G(P)$ is dense in $\Prim(A)$, since otherwise
$D$ would be a subset
of the non-empty open set $G\backslash\!\big(\Prim(A)\smallsetminus \overline{G(P)}\big)$, 
which is impossible.

If the dense orbit $G(P)$ is locally closed then $G(P)$ is open in its closure $\Prim(A)$,
which implies that $G(P)$ is the unique dense orbit in $\Prim(A)$. This gives 
the uniqueness assertion of the proposition.
\end{proof}

Recall that if $G(P)$ is a locally closed orbit in $\Prim(A)$, we get
$A_{G(P)}\rtimes G\cong (I_{G(P)}\rtimes G)/(J_{G(P)}\rtimes G)$  and similarly
$A_{G(P)}\rtimes_r G\cong (I_{G(P)}\rtimes_r G)/(J_{G(P)}\rtimes_r G)$ if
$G$ is exact. Using this we get

\begin{corollary}\label{cor-smooth}
Suppose that $(A,G,\alpha)$ is smooth. Then we obtain a decomposition of
$(A\rtimes G)\dach$ 
(resp. $\Prim(A\rtimes G)$) as the disjoint
union of the locally closed subsets $(A_{G(P)}\rtimes G)\dach$ (resp. 
$\Prim(A_{G(P)}\rtimes G)$), where $G(P)$ runs through all $G$-orbits
in $\Prim(A)$.
If $G$ is exact, similar statements hold for the reduced crossed products.
\end{corollary}
\begin{proof} It follows from Proposition \ref{prop-live}
that for each $\pi\times U\in (A\rtimes G)\dach$, there exists a unique 
orbit $G(P)$ such that $\ker\pi=P^G=J_{G(P)}$ and then $\pi\times U$
restricts to an irreducible representation of $A_{G(P)}\rtimes G$.
Hence
 $$(A\rtimes G)\dach=\cup\{(A_{G(P)}\rtimes G)\dach :{G(P)\in G\backslash \Prim(A)}\}.$$
To see that this union is disjoint, assume that there exists an element
 $\rho\times V\in (A_{G(P)}\rtimes G)\dach$ (viewed as a representation of $A\rtimes G$)
 such that $\ker\rho\neq P^G$. Since $\rho$ is a representation 
 of $A_{G(P)}=I_{G(P)}/P^G$ we have $\ker\rho\supseteq P^G$.
  By Proposition \ref{prop-live} there exists a $Q\in \Prim(A)$ such that 
 $\ker\rho= Q^G$. Then  $Q^G\supseteq P^G$, which implies that 
 $G(Q)\subseteq \big(\overline{G(P)} \smallsetminus G(P)\big)$.
But then
$$\ker\rho=Q^G=\ker G(Q)\supseteq\ker\big(\overline{G(P)}\smallsetminus G(P)\big)=I_{G(P)},$$
which contradicts the assumption that $\rho\times V\in (A_{G(P)}\rtimes G)\dach$.
\end{proof}

It is now easy to give a proof of the Mackey-Green-Rieffel theorem, which
is the main result of this section. If $(A,G,\alpha)$ is smooth, one can easily check 
that points in $\Prim(A)$ are automatically locally closed (since they are closed in 
their orbits). 
Hence, for each $P\in \Prim(A)$ the restriction $A_P:=I_P/P$ of $A$ to $\{P\}$
is a simple subquotient of $A$. 
Since $I_P$ and $P$ are invariant under the action of the stabiliser $G_P$,
the action of $G_P$ on $A$ factors through an action of $G_P$ on $A_P$.
It is then straightforward to check (using the same arguments as given in the 
proof of Corollary \ref{cor-smooth})  that
there is a canonical one-to-one correspondence between the 
irreducible  representations of $A_P\rtimes G_P$ and the set of all irreducible 
representations $\rho\times V$ of $A\rtimes G_P$ satisfying $\ker\rho=P$.
In case where $A=C_0(X)$ is commutative, we will study this problem in \S \ref{subsec-trans} below, and 
the case where $A$ is type I will be studied in  \S \ref{subsec-littlegroup}.

\begin{remark}\label{rem-typeImackey}
If $A$ is type I, then $A_P\cong \K(H_\pi)$, the compact operators on the Hilbert space 
$H_\pi$, where $\pi:A\to \B(H_\pi)$ is the unique (up to equivalence) irreducible 
representation of $A$ with $\ker\pi=P$. To see this we first pass to $A/P\cong \pi(A)$.
Since $A$ is type I we know that $\K(H_\pi)\subseteq \pi(A)$.  Hence, if we identify 
$\K(H_\pi)$ with an ideal of $A/P$, we see (since $\pi$ does not vanish on this ideal)
that this ideal must correspond to 
the open set $\{\pi\}$ (resp. $\{P\}$) in its closure $\widehat{A/P}$ (resp. $\Prim(A/P)$).
\end{remark}

\begin{theorem}[Mackey-Rieffel-Green]\label{thm-MRG}
Suppose that $(A,G,\alpha)$ is  smooth. Let $\mathcal S\subseteq \Prim(A)$ be 
a cross-section for the orbit space $G\backslash \Prim(A)$, i.e. $\mathcal S$ intersects each orbit $G(P)$ 
in exactly one point.
Then induction of representations and ideals induces bijections
$$
\Ind:\cup_{P\in \mathcal S} (A_P\rtimes G_P)\dach \to (A\rtimes G)\dach;\rho\times V\mapsto \ind_{G_P}^G(\rho\times V)\quad \text{and}
$$
$$
\Ind:\cup_{P\in \mathcal S} \Prim(A_P\rtimes G_P) \to \Prim(A\rtimes G); \; Q\mapsto \ind_{G_P}^GQ.
$$
 If $G$ is exact, these maps restrict to similar bijections 
$$\cup_{P\in \mathcal S} (A_P\rtimes_r G_P)\dach \stackrel{\Ind}{\to} (A\rtimes_r G)\dach\quad\text{and}
\quad
\cup_{P\in \mathcal S} \Prim(A_P\rtimes_rG_P) \stackrel{\Ind}{\to} \Prim(A\rtimes_rG)$$
for the reduced crossed products.
\end{theorem}
\begin{proof}
We show that
$\Ind:\cup_{P\in \mathcal S}(A_P\rtimes G_P)\dach \to(A\rtimes G)\dach$ is a bijection. 
Bijectivity of the other maps follows similarly.

By Corollary \ref{cor-smooth} it suffices to show that 
$\Ind:(A_P\rtimes G_P)\dach \to (A_{G(P)}\rtimes G)\dach$ is a bijection for all $P\in \mathcal S$. By definition of $A_{G(P)}$ we have $\Prim(A_{G(P)})\cong G(P)$ 
and by the smoothness of the action we have $G(P) \cong G/G_P$ as $G$-spaces.
Hence, it follows from 
 Theorem \ref{thm-ind}  that $A_{G(P)}\cong\Ind_{G_P}^G A_P$.
Hence induction via Green's $A_{G(P)}\rtimes G-A_P\rtimes G_P$ imprimitivity
bimodule $X_P:=X_{G_P}^G(A_P)$ gives the desired bijection 
$\ind^{X_P}:(A_P\rtimes G_P)\dach \to (A_{G(P)}\rtimes G)\dach$.
By Corollary \ref{cor-inducedalg}, induction via $X_P$ coincides with the usual
induction for the system $(A_{G(P)},G,\alpha)$, which by Remark \ref{compat-ind}
is  compatible
with inducing the corresponding representations for  $(A,G,\alpha)$.
\end{proof}

The above result shows that for smooth systems, all representations are induced 
from the stabilisers for the corresponding action of $G$ on $\Prim(A)$. In fact the above result
is much stronger, since it shows that $A\rtimes G$ has a ``fibration'' over $G\backslash \Prim(A)$
such that the fiber $A_{G(P)}\rtimes G$ over an orbit $G(P)$ is Morita equivalent
to $A_P\rtimes G_P$, hence, up to the global structure of the fibration, the study 
of $A\rtimes G$ reduces to the study of the fibers $A_P\rtimes G_P$. Note that 
under the assumptions of Theorem \ref{thm-MRG} the algebra $A_P$ is always simple.
We shall give a more detailed study of the crossed products $A_P\rtimes G_P$ in 
the important special case where $A$ is type I  in
\S\ref{subsec-littlegroup} below. The easier situation where $A=C_0(X)$ is treated in 
\S\ref{subsec-trans} below.

Note that the study of the global structure
of $A\rtimes G$, i.e., of the global structure of the fibration over $G\backslash \Prim(A)$
is in general quite complicated, even in the situation where $G\backslash\Prim(A)$ is Hausdorff. 
In general, it is also very difficult (if not impossible) to describe the global topology
of $\Prim(A\rtimes G)$ in terms of the bijection of Theorem \ref{thm-MRG}.
Some progress has been made in the case where $A$ is a continuous-trace $C^*$-algebra
and/or where the stabilisers are assumed to vary continuously, 
and we refer to \cites{CKRW, E-ct,  EE, ER, EN, EW2, EW-proper, RW} and the references
given in those papers and books for more information on this problem. 
If $A=C_0(X)$ is commutative and $G$ is abelian, a very satisfying description  
of the topology of $\Prim(C_0(X)\rtimes G)$ has been obtained by Dana Williams 
in \cite{Wi0}. We shall discuss this situation in  \S\ref{subsec-trans} below.

Even worse, the assumption of having a smooth action is a very strong one
and for arbitrary systems one cannot expect that one can compute all irreducible
representations via induction from stabilisers. Indeed, in general it is not possible
to classify all irreducible representations of a non-type I $C^*$-algebra, and a 
similar problem occurs for crossed products $A\rtimes G$ if the action of 
$G$ on $\Prim(A)$ fails to be smooth.
However, at least if $(A,G,\alpha)$ is separable and $G$ is amenable,
there is a positive result towards the description 
of $\Prim(A\rtimes G)$ which was obtained by work of Sauvageot and Gootman-Rosenberg,
thus giving a positive answer to an earlier formulated  conjecture by Effros and Hahn
(see \cite{EH}). To give precise statements, we need 

\begin{definition}\label{def-homogeneous-rep}
A non-degenerate representation $\rho$ of a $C^*$-algebra $A$ is called 
{\em homogeneous} if all non-trivial subrepresentations of $\rho$
have the same kernel as $\rho$. 
\end{definition}

It is clear that every irreducible representation is homogeneous 
and one can show that the kernel of any  homogeneous representation
is a prime ideal, and hence it is primitive if $A$ is second countable.
We refer to \cite{Wi-crossed} for a discussion on this and for 
very detailed proofs of Theorems \ref{thm-sauv} and \ref{thm-GR} stated below:

\begin{theorem}[Sauvageot (\cite{S2})]\label{thm-sauv}
Suppose that $(A,G,\alpha)$ is a separable system (i.e., $A$ is separable
and $G$ is second countable). Let 
$P\in \Prim(A)$ and let $G_P$ denote the stabiliser of $P$ in $G$.
Suppose that $\rho\times V$ is a homogeneous representation
of $A\rtimes G_P$ such that $\rho$ is a homogeneous representation of 
$A$ with $\ker \rho=P$. Then $\ind_{G_P}^G(\rho\times V)$ 
is a homogeneous representation of $A\rtimes G$ and
$\ker\big(\ind_{G_P}^G(\rho\times V)\big)$ is a primitive ideal of $A\rtimes G$.
\end{theorem}

We say that a primitive ideal of $A\rtimes G$ is {\em induced}
if it is obtained as in the above theorem. Note that Sauvageot already showed
in \cite{S2} that in case where $G$ is amenable, every primitive ideal of $A\rtimes G$
contains an induced primitive ideal and in case where $G$ is discrete
every primitive ideal is contained in an induced primitive ideal.
Together, this shows that for actions of discrete amenable groups all primitive
ideals of $A\rtimes G$ are induced from the stabilisers.
Sauvageot's  result was generalized by Gootman and Rosenberg in \cite{GR}*{Theorem 3.1}:

\begin{theorem}[Gootman-Rosenberg]\label{thm-GR}
Suppose that $(A,G,\alpha)$ is a separable system. Then every primitive ideal of 
$A\rtimes G$ is contained in an induced ideal. As a consequence, if $G$ is amenable,
then every primitive ideal of $A\rtimes G$ is induced.
\end{theorem}
%
%Note that similar statements also hold for twisted systems, which follows 
%easily from Theorem \ref{thm-twistedequiv}. 

The condition in Theorem \ref{thm-sauv} that the representations 
$\rho\times V$ and $\rho$ are homogeneous is a little bit unfortunate.
In fact, a somehow more natural formulation of Sauvageot's theorem
(using the notion of induced ideals) would be to state that whenever 
$Q\in \Prim(A\rtimes G_P)$ such that
 $\res_{\{e\}}^{G_P}(Q)=P$, then $\ind_{G_P}^G(Q)$ is a primitive ideal
 of $A\rtimes G$. 
 Note that if $\rho\times V$ is as in Theorem \ref{thm-sauv},
then $Q=\ker(\rho\times V)$ is an element of $\Prim(A\rtimes G_P)$
which satisfies the above conditions.
At present time, we do not know whether this more general statement 
is true, and we want to take this opportunity to 
point out that the statement of  \cite{E-ct}*{Theorem 1.4.14} is not correct 
(or at least not known) as it stands.
We are very grateful to Dana Williams 
for pointing out this error
and we refer to the paper \cite{EW-induced} for a more elaborate discussion of this problem. 
But let us indicate here that the problem vanishes if all points in 
$\Prim(A)$ are locally closed (which is in particular true if $A$ is type I).

\begin{proposition}\label{prop-sauv}
Suppose that $(A,G,\alpha)$ is a separable system such
 that one of the following conditions is satisfied:
\begin{enumerate}
\item All points in $\Prim(A)$ are locally closed (which is automatic if 
$A$ is type I).
\item All stabilisers $G_P$ for $P\in \Prim(A)$ are normal subgroups of $G$
(which is automatic if $G$ is abelian).
\end{enumerate}
Then $\ind_{G_P}^GQ\in \Prim(A\rtimes G)$ for all $P\in \Prim(A)$ and 
$Q\in \Prim(A\rtimes G_P)$
such that $\res_{\{e\}}^{G_P}Q=P$. 
If, in addition, $G$ is amenable, then all primitive ideals of $A\rtimes G$ are 
induced in this way.
\end{proposition}
\begin{proof}
Let us first assume condition (i). Choose $\rho\times V\in (A\rtimes G_P)\dach$
such that $\ker(\rho\times V)=Q$ and $\ker\rho=P$. Then we may
regard $\rho$ as a representation of $A_P$, the simple subquotient of
$A$ corresponding to the locally closed subset $\{P\}$ of $\Prim(A)$.
Since $A_P$ is simple, all nontrivial subrepresentations of $\rho$
have kernel $\{0\}$ in $A_P$ (and hence they have kernel $P$ in $A$).
Hence $\rho$ is homogeneous and 
the result follows from Theorems \ref{thm-sauv} and \ref{thm-GR}.

Let us now assume (ii). If $N:=G_P$ is normal, we may use the theory of twisted actions,
which we shall present in \S \ref{sec-twists} below,
to pass to the system $((A\rtimes N)\otimes\K, G/N, \beta)$. If $\rho\times V\in (A\rtimes N)\dach$ 
with $\ker(\rho\times V)=P$, then the corresponding representation  of 
$(A\rtimes N)\otimes \K$ has trivial stabiliser in $G/N$, and therefore
the induced representation has primitive kernel in 
$A\rtimes G\sim_M \big((A\rtimes N)\otimes\K\big)\rtimes G/N$
 by Theorem \ref{thm-sauv}.
\end{proof}

Recall that if $M$ is a topological $G$-space, then 
two elements $m_1,m_2\in M$ are said to be in the same 
{\em quasi-orbit} if $m_1\in \overline{G(m_2)}$ and $m_2\in \overline{G(m_1)}$.
Being in the same quasi-orbit is clearly an equivalence relation on $M$
and we denote by $G_q(m)$ the quasi-orbit (i.e., the equivalence class)
of $m$ and by $\Q_G(M)$ the set of all quasi-orbits in $M$ equipped 
with the quotient topology. 
Note that 
$\Q_G(M)$ is  always a T$_0$-space.
If $G\backslash M$ is a T$_0$-space, then 
$\Q_G(M)$ coincides with $G\backslash M$.

If $(A,G,\alpha)$ is a system, it follows from the definition of the 
Jacobson topology that two elements $P,Q\in \Prim(A)$ are in the same 
quasi-orbit if and only if $P^G=Q^G$. 
If the action of $G$ on $A$ is smooth, then all points in $G\backslash\Prim(A)$ are
locally closed, which implies in particular that $G\backslash\Prim(A)$ is
a T$_0$-space. Hence in this case we have $\Q_G(\Prim(A))=G\backslash\Prim(A)$.
In what follows, we let
$$\Prim^G(A):=\{P^G: P\in \Prim(A)\}\subseteq \I(A)$$ 
 equipped  with the relative Fell topology.
Then  \cite{Green1}*{Lemma on p. 221}  gives

\begin{lemma}\label{lem-quasi}
Let $(A,G,\alpha)$ be a system. Then the map 
$$q:\Prim(A)\to \Prim^G(A):P\mapsto P^G$$ 
is a continuous and open surjection
and therefore factors through a homeomorphism 
between $\Q_G(\Prim(A))$ and $\Prim^G(A)$.
\end{lemma}

As a consequence of the previous results we get

\begin{corollary}\label{cor-simple}
Suppose that $(A,G,\alpha)$ is smooth or that $(A,G,\alpha)$ is separable and $G$ is 
amenable. Suppose further that the action of $G$ on $\Prim(A)$ is free (i.e., all stabilisers 
are trivial). 
Then the map 
$$\Ind:\Prim^G(A)\cong \Q_G(\Prim(A))\to \Prim(A\rtimes G); P\mapsto \ind_{\{e\}}^GP^G$$ 
is a homeomorphism. 
In particular, $A\rtimes G$ is simple if and only if every $G$-orbit is dense
in $\Prim(A)$, and $A\rtimes G$ is primitive (i.e., $\{0\}$ is a primitive ideal of $A\rtimes G$)
if and only if there exists a dense $G$-orbit in $\Prim(A)$.
\end{corollary}
\begin{proof} 
It follows from  Theorem \ref{thm-MRG} and Theorem \ref{thm-GR} that the map
$\ind_{\{e\}}^G:\Prim(A)\to \Prim(A\rtimes G); P\mapsto \ind_{\{e\}}^GP$ is well defined and 
surjective.
By Corollary \ref{cor-regular} we know that
 $\ind_{\{e\}}^GP=\ind_{\{e\}}^GP^G$, so the induction map 
 $\Ind:\Prim^G(A)\to \Prim(A\rtimes G)$ is also well defined and surjective.
Equation (\ref{eq-indresideal}) applied to $H=\{e\}$ gives
 $\res_{\{e\}}^G(\ind_{\{e\}}^GP)=P^G$, which shows that
 $\res_{\{e\}}^G:\Prim(A\rtimes G)\to \Prim^G(A)$
 is the inverse of $\Ind$. Since induction and restriction are continuous
by Proposition \ref{prop-contHG} the result follows.
\end{proof}

A quite recent result of Sierakowski (see \cite{Sier}*{Proposition 1.3 and Theorem 1.20} and \cite{EL}*{Corollary 2.7}) shows that 
for countable discrete groups $G$, the assumptions for the action of $G$ on $\Prim(A)$ and of amenability of $G$ 
can be weakened considerably. We need

\begin{definition}
An action of a group $G$ on a topological space $X$ is called {\em essentially free} if every $G$-invariant closed 
subset $C\subseteq X$ contains a dense $G$-invariant subset $D$ such that $G$ acts freely on $D$.
\end{definition}

\begin{theorem}[Sierakowski] Suppose that $(A,G,\alpha)$ is a $C^*$-dynamical system with $A$ separable and $G$ countable (hence discrete) and exact.
Suppose further that the action of $G$ on $\widehat{A}$ is essentially free (which is true, if the action of $G$ on $\Prim(A)$ is essentially free).
Then the map 
$$\Ind: \Prim^G(A)\cong \Q_G(\Prim(A))\to \Prim(A\rtimes_rG); P\mapsto \Ind_e^GP$$
is a well-defined homeomorphism.
\end{theorem}

Since the induced ideals in the above theorem clearly contain the kernel of the regular representation of $A\rtimes G$, it is clear that
a similar statement cannot hold for the full crossed product $A\rtimes G$ if it differs from $A\rtimes_rG$.

We should note that Sierakowski's original results \cite{Sier}*{Proposition 1.3 and Theorem 1.20} show, that under the assumptions 
of the theorem the map $\res: \I(A\rtimes_rG)\to \I^G(A); J\mapsto J\cap A$ is a bijection between the 
set of closed two-sided ideals in $A\rtimes_rG$ and the set of $G$-invariant closed  two-sided ideals in $A$, with inverse given by
$I\mapsto I\rtimes_rG$. The straight-forward translation of this into the statement of the above theorem has been 
given in \cite{EL}. We should also mention that Sierakowski's result still holds under some slightly weaker assumptions, which he calls the
{\em residual Rokhlin* property}. We refer the interested reader to \cite{Sier} for more details on this property.

\begin{remark}\label{rem-conststab} 
If $(A,G,\alpha)$ is a system with constant stabiliser $N$ for the
action of $G$ on $\Prim(A)$, then $N$ is normal in $G$ and one can pass to the iterated twisted 
system \mbox{$(A\rtimes N, G, N, \alpha^N, \tau^N)$} (see \S\ref{sec-twists} below), and then 
to an equivariantly  Morita equivalent system $(B,G/N,\beta)$ (see Proposition \ref{prop-twistequiv}) to see that
induction of primitive ideals gives a homeomorphism between 
$\Q_{G/N}(\Prim(A\rtimes N))$ and $\Prim(A\rtimes G)$ if one of the following conditions 
are satisfied:
\begin{enumerate}
\item $(A,G,\alpha)$ is smooth.
\item $(A\rtimes N, G, N,\alpha^N,\tau^N)$ is smooth (i.e., the action of $G/N$ on $\Prim(A\rtimes N)$
via $\alpha^N$ satisfies the conditions of Definition \ref{defn-smooth}.
\item $(A,G,\alpha)$ is separable and $G/N$ is amenable.
\end{enumerate}
A similar result can be obtained for systems with continuously 
varying stabilisers (see \cite{E-prim}). In the case of constant stabilisers,
the problem of describing the topology of $\Prim(A\rtimes G)$ now reduces 
to the description the topology of $\Prim(A\rtimes N)$ and the action of $G/N$ on 
$\Prim(A\rtimes N)$. In general, both parts can be quite difficult to perform, but in some
interesting special cases, e.g. if $A$ has continuous trace,
some good progress has been made for the description of $\Prim(A\rtimes N)$
(e.g. see \cites{EW1, EW2, EN} and the references given there). Of course, if 
$A=C_0(X)$ is abelian, and $N$ is the constant stabiliser of the elements
of $\Prim(A)=X$, then $N$ acts trivially on $X$ and $\Prim(C_0(X)\rtimes N)=
\Prim(C_0(X)\otimes C^*(N))=X\times\Prim(C^*(N))$.
\end{remark}

\begin{example}\label{ex-Atheta}
As an  easy application of Corollary \ref{cor-simple} we get the simplicity of the irrational rotation 
algebra $A_{\theta}$, for $\theta$ an irrational number 
in $(0,1)$. Recall  that $A_{\theta}= C(\TT)\rtimes_{\theta} \ZZ$
where $n\in \ZZ$ acts on $z\in \TT$ via $n\cdot z:= e^{2\pi i \theta n}z$.
Since $\theta$ is irrational, the action of $\ZZ$ on $\Prim(C(\TT))=\TT$
is free and all $\ZZ$-orbits are dense in 
$\TT$. Hence, there exists only one quasi-orbit in $\TT$ and the crossed product 
is simple. Of course, there are other more  elementary proofs for the simplicity
of $A_{\theta}$ which do not use such heavy machinery, but this example 
illustrates quite well how one can use the above results.
\end{example}

\subsection{The Mackey-machine  for transformation groups}\label{subsec-trans}
Suppose that $X$ is a locally compact $G$-space and consider the corresponding
action of $G$ on $A=C_0(X)$ given by $(s\cdot\varphi)(x)=\varphi(s^{-1}x)$ for 
$s\in G$, $\varphi\in C_0(X)$. 
Then $\Prim(A)=X$ and $A_x\cong \CC$ for all $x\in X$, so that $A_x\rtimes G_x\cong C^*(G_x)$
for all $x\in X$, where $G_x$ denotes the stabiliser of $x$. 
Hence, if the action of $G$ on $X$
is smooth in the sense of Definition \ref{defn-smooth}, then it follows from 
Theorem \ref{thm-MRG} that $C_0(X)\rtimes G$ is ``fibered'' over 
$G\backslash X$ with fibres $C_0(G(x))\rtimes G\sim_M C^*(G_x)$ (compare the discussion
following Theorem \ref{thm-MRG}).

If $V\in \widehat{G_x}$ and if $\eps_x:C_0(X)\to \CC$ denotes evaluation at $x$,
then $\eps_x\times V$ is  the representation of $C_0(X)\rtimes G_x$ which corresponds 
to $V$ by regarding  $\widehat{G_x}\cong(A_x\rtimes G_x)\dach$ as a subset of 
$(A\rtimes G_x)\dach$ as described in the discussion preceeding Theorem \ref{thm-MRG}.
In this situation, the result of Theorem \ref{thm-sauv} can be improved by showing:

\begin{proposition}[{cf.  \cite{Wi0}*{Proposition 4.2}}]\label{prop-irred}
Let 
$\eps_x\times V\in (C_0(X)\rtimes G_x)\dach$ be as above.
Then $\ind_{G_x}^G(\eps_x\rtimes V)$ is irreducible. Moreover, if $V,W\in \widehat{G_x}$,
then $$\ind_{G_x}^G(\eps_x\times V)\cong \ind_{G_x}^G(\eps_x\times W)\quad\Longleftrightarrow\quad V\cong W.$$
\end{proposition}

Combining this with Theorem \ref{thm-MRG} and Theorem \ref{thm-GR} gives:

\begin{theorem}\label{thm-Mtrans}
Suppose that $X$ is a locally compact $G$-space. 
\begin{enumerate}
\item If $G$ acts smoothly on $X$,
and if $\mathcal S\subseteq X$ is a section for $G\backslash X$, then we get a bijection
$$\Ind: \cup_{x\in {\mathcal S}}\widehat{G_x}\to \big(C_0(X)\rtimes G\big)\dach;
V\mapsto \ind_{G_x}^G(\eps_x\times V).$$

\item If $X$ and $G$ are second countable and if $G$ is amenable, then 
 every primitive ideal of $C_0(X)\rtimes G$ is the kernel of 
some induced irreducible representation $\ind_{G_x}^G(\eps_x\times V)$.
\end{enumerate}
\end{theorem}

If  $G$ is abelian, then so are the stabilisers $G_x$ for all $x\in X$. Then $\widehat{G}_x$ 
is the Pontrjagin dual group of $G_x$ and we get a short exact sequence
$$0\to \widehat{G/G_x}\to \widehat{G}\stackrel{\res}{\to}\widehat{G}_x\to 0$$
for all $x\in X$. Moreover, since $G_x$ is normal in $G$, it follows that the stabilisers are 
constant on quasi-orbits $G_q(x)$ in $X$. 
We can then consider an equivalence relation on $X\times \widehat{G}$ by 
$$(x,\chi)\sim (y, \mu)\Leftrightarrow G_q(x)=G_q(y)\;\text{and}\; \chi|_{G_x}=\mu|_{G_y}.$$
Then the following result is \cite{Wi0}*{Theorem 5.3}:

\begin{theorem}[Williams]\label{thm-abelian} Suppose that $G$ is abelian and the action of $G$ on $X$ is smooth or $G$ and $X$ are 
second countable. Then the map
$$\Ind: (X\times \widehat{G})/\sim \to\Prim(C_0(X)\rtimes G); [(x,\chi)]\mapsto \ker(\Ind_{G_x}^G (\eps_x\rtimes \chi|_{G_x}))$$
is a homeomorphism.
\end{theorem}

Note that in case where the action of $G$ on $X$ is smooth and $G$ is abelian, the crossed product $C_0(X)\rtimes G$ 
is type I, since $C^*(G_x)$ is type I. Hence in this case we get a homeomorphism between 
$(X\times \widehat{G})/\sim $ and $(C_0(X)\rtimes G)\dach$.
We now want to present some applications to group representation theory:
 
 \begin{example}\label{ex-trans}
 Suppose that $G=N\rtimes H$ is the semi-direct product of the 
 abelian group $N$ by the group $H$. Then, as seen in Example \ref{ex-deco}, we have 
 $$C^*(N\rtimes H)\cong C^*(N)\rtimes H\cong
 C_0(\widehat{N})\rtimes H,$$
 where the last isomorphism is given via the Gelfand-transform 
 $C^*(N)\cong C_0(\widehat{N})$. The corresponding action of $H$ on 
 $C_0(\widehat{N})$ is induced by the action of $H$ on $\widehat{N}$ 
 given by $\big(h\cdot\chi\big)(n):=\chi(h^{-1}\cdot n)$ if $h\in H$, $\chi\in \widehat{N}$ 
 and $n\in N$.
 Thus, if the action of $H$ on $\widehat{N}$ is smooth, we obtain every irreducible 
 representation of $C^*(N\rtimes H)\cong C_0(\widehat{N})\rtimes H$ as an induced 
 representation $\ind_{H_{\chi}}^H(\eps_\chi\times V)$ 
 for some $\chi\in \widehat{N}$ and  $V\in \widehat{H}_\chi$. The
 isomorphism $C_0(\widehat{N})\rtimes H_{\chi}\cong C^*(N\rtimes H_\chi)$,
 transforms the representation $\eps_{\chi}\times V$ to the representation 
 $\chi\times V$ of $N\rtimes H_{\chi}$ defined by $\chi\times V(n,h)=\chi(n)V(h)$
 and one can show that 
 $\ind_{N\rtimes H_{\chi}}^{N\rtimes H}(\chi\times V)$ corresponds to the representation
 $\ind_{H_\chi}^H(\eps_\chi\times V)$ under the isomorphism
 $C^*(N\rtimes H)\cong C_0(\widehat{N})\rtimes H$.
 Thus, choosing a cross-section $\mathcal S\subseteq \widehat{N}$ for $H\backslash\widehat{N}$,
 it follows from Theorem \ref{thm-MRG} that
 $$\Ind:\cup\{\widehat{H}_{\chi}:\chi\in \mathcal S\}\to \widehat{N\rtimes H}; V\mapsto \ind_{N\rtimes H_{\chi}}^{N\rtimes H} (\chi\times V)$$
 is a bijection. 
 
 If the action of $H$ on $\widehat{N}$ is not smooth, but $N\rtimes H$ is second countable and amenable, then we get at least all 
 primitive ideals of $C^*(N\rtimes H)$ as kernels of the induced representations
 $\ind_{N\rtimes H_\chi}^{N\rtimes H}(\chi\times V)$.
 
Let us now discuss some explicit examples:

{\bf (1)} Let $G=\RR\rtimes \RR^*$ denote the $ax+b$-group, i.e., $G$ is the 
semi-direct product for the action of the multiplicative group $\RR^*$ on $\RR$
via dilation. Identifying $\RR$ with $\widehat{\RR}$ via $t\mapsto\chi_t$ with $\chi_t(s)=e^{2\pi i ts}$,
we see easily that the 
action of $\RR^*$ on $\widehat{\RR}$ is also given by dilation. Hence there are precisely
two orbits in $ \widehat{\RR}$: $\{\chi_0\}$ and $\widehat{\RR}\smallsetminus\{\chi_0\}$.
Let $\mathcal S=\{\chi_0,\chi_1\}\subseteq \widehat{\RR}$. Then $\mathcal S$ is a cross-section 
for $\RR^*\backslash\widehat{ \RR}$,  the stabiliser of $\chi_1$ in $\RR^*$ is $\{1\}$ and the stabiliser
of $\chi_0$ is all of $\RR^*$. Thus we see that 
$$\widehat{G}=\{\chi_0\times\mu :\mu\in \widehat{\RR^*}\}\cup\{\ind_{\RR}^{\RR\rtimes\RR^*}\chi_1\}.$$
It follows from Theorem \ref{thm-abelian} that 
the single representation $\pi:=\ind_{\RR}^{\RR\rtimes\RR^*}\chi_1$ is dense in $\widehat{G}$ and that 
the set $\{\chi_0\times\mu :\mu\in \widehat{\RR^*}\}\subseteq \widehat{G}$ is homeomorphic to $\widehat{\RR^*}\cong \RR^*$.

Notice that we could also consider the $C^*$-algebra $C^*(G)$ as ``fibered'' over 
$\RR^*\backslash\widehat{\RR}$: The open orbit $\widehat{\RR}\smallsetminus\{\chi_0\}\cong\RR^*$
corresponds to the ideal $C_0(\RR^*)\rtimes\RR^*\cong\K(L^2(\RR^*))$ and the 
closed orbit $\{\chi_0\}$ corresponds to the quotient $C_0(\widehat{\RR^*})$ of 
$C^*(G)$, so that this picture yields the short exact sequence
$$0\to \K(L^2(\RR^*))\to C^*(G)\to C_0(\widehat{\RR^*})\to 0$$
(compare also with Example \ref{ex-Green}). 

{\bf (2)} A more complicated example is given by the  Mautner group. 
This group is the semi-direct product
$G=\CC^2\rtimes \RR$ with action given by 
$$t\cdot (z, w)= (e^{-2\pi i t}z, e^{-2\pi i\theta t} w),$$
where $\theta\in (0,1)$ is a fixed irrational number. Identifying $\CC^2$ with the dual group
$\widehat{\CC^2}$ via $(u,v)\mapsto \chi_{(u,v)}$ such that
$$\chi_{(u,v)}(z,w)=\exp(2\pi i \text{Re}(z\bar{u}+w\bar{v})),$$
we get $t\cdot\chi_{(u,v)}=\chi_{(e^{2\pi i t}u, e^{2\pi i \theta t} z)}$.
The quasi-orbit space for the action of $\RR$ on $\widehat{\CC^2}$ can then be parametrized by
the set $[0,\infty)\times [0,\infty)$: If $(r,s)\in [0,\infty)^2$,
then the corresponding quasi-orbit $\mathcal O_{(r,s)}$ 
consists of all $(u,v)\in \CC^2$ such that $|u|=r$ and $|v|=s$.
Hence, if $r,s>0$, then $\mathcal O_{(r,s)}$ is homeomorphic to $\TT^2$ 
and this homeomorphism  carries the action of $\RR$ on $\mathcal O_{(r,s)}$ 
to the irrational flow 
of $\RR$ on $\TT^2$ corresponding to $\theta$ as considered in 
part (4) of Example \ref{ex-Green}. 
In particular, $\RR$ acts freely but not smoothly
on those quasi-orbits.
If $r\neq 0$ and $s= 0$, the quasi-orbit $\mathcal O_{(r,s)} $ is
homeomorphic to $\TT$ with action $t\cdot u:=e^{2\pi i t} u$ and constant stabiliser
$\ZZ$. In particular, all those
quasi-orbits are orbits. Similarly, if $r=0$ and $s\neq 0$,
the quasi-orbit $\mathcal O_{(r,s)}$ is homeomorphic to $\TT$ with action 
$t\cdot v= e^{2\pi i \theta t}v$ and stabiliser $\frac{1}{\theta}\ZZ$. 
Finally, the quasi-orbit corresponding to $(0,0)$ is 
the point-set $\{(0,0)\}$ with stabiliser $\RR$. 

Since $G$ is second countable and amenable, we can therefore parametrize $\Prim(C^*(G))$
by the set
$$\{(r,s): r,s >0\}\cup  \big((0,\infty)\times\widehat{\ZZ}\big)\cup \big((0,\infty)\times{\widehat{\mbox{\small{$\frac{1}{\theta}$}}\ZZ}}\big)\cup \widehat{\RR}.$$
In fact, we can also view $C^*(G)$ as ``fibered'' over $[0,\infty)^2$ with 
fibers 
$$C^*(G)_{(r,s)}\cong C(\TT^2)\rtimes_{\theta}\RR\sim_M A_{\theta}\quad\quad\text{for 
$r,s>0$,}$$
 where $A_\theta$ denotes the irrational rotation algebra,
 $$C^*(G)_{(r,0)}\cong C(\TT)\rtimes \RR\sim_M C(\widehat{\ZZ})\cong C(\TT)\quad \quad \text{
for $r>0$,}$$ 
$$C^*(G)_{(0,s)}\cong C(\TT)\rtimes_{\theta}\RR\sim_MC(\widehat{\mbox{\small{$\frac{1}{\theta}$}}\ZZ})
\cong C(\TT)\quad\quad\text{for $s>0$},$$ 
and $C^*(G)_{(0,0)}\cong C_0(\RR)$. Using Theorem \ref{thm-abelian}, it is 
also possible to describe the topology of $\Prim(G)$, 
but we do not go into the details here.
We should mention that the Mautner group is the lowest dimensional example of a 
connected Lie-group $G$ with a non-type I group algebra $C^*(G)$.
\end{example}

 \begin{remark}\label{rem-littlegroup}
 It follows from Theorems \ref{thm-MRG} and \ref{thm-GR} that for understanding 
the ideal structure of $A\rtimes G$, it is necessary to understand the structure
of $A_P\rtimes G_P$ for $P\in \Prim(A)$. We saw in this section that 
this is the same as understanding the group algebras $C^*(G_x)$ for the stabilisers
$G_x$ if $A=C_0(X)$ is abelian. In general, the problem becomes much more
difficult.
 However, at least in the important special case where $A$ is type I
 one can still give a quite satisfactory description of $A_P\rtimes G_P$ in terms of the 
 stabilisers. 
 Since an elegant treatment of that case uses the theory of twisted
 actions and crossed products, we postpone the discussion of this case 
 to \S \ref{subsec-littlegroup} below.
 \end{remark}

\section{The Mackey-Rieffel-Green machine for twisted crossed products}\label{sec-twists}
 
\subsection{Twisted actions and twisted crossed products}\label{subsec-twists}

One draw-back of the theory of crossed products by ordinary actions
is the fact that crossed products $A\rtimes G$ (and their reduced
analogues) cannot be written as iterated crossed products
$(A\rtimes N)\rtimes G/N$ if $N$ is a normal subgroup such that
the extension
$$1\to N\to G\to G/N\to 0$$
is not topologically split (compare with 
Example \ref{ex-deco}). In order to close this gap, we now 
introduce 
twisted actions and twisted crossed products following 
Phil Green's 
approach of \cite{Green1}. Note that there is an alternative 
approach due to Leptin and 
Busby--Smith (see \cites{Lep1, BuS, PR89} 
for the construction 
of twisted crossed products within this theory), but Green's 
theory
seems to be better suited for our purposes.

As a motivation, 
consider a closed normal subgroup $N$ of the locally compact 
group 
$G$, and assume that $\alpha:G\to\Aut(A)$ is an action. Let 
$A\rtimes N$ be the crossed product of $A$ by $N$. 
Let $\delta:G\to 
\RR^+$ be the module for the conjugation action of $G$ on $N$,
i.e., 
$\delta(s) \int_N f(s^{-1}ns)\,dn=\int_Nf(n)\,dn$ for all $f\in 
C_c(N)$.
A short computation using the formula 
\begin{equation}\label{eq-Weil}
\int_G g(s)\,ds=\int_{G/N}\left( 
\int_N g(sn)\,dn\right)\,dsN
\end{equation}
(with respect to suitable 
choices of Haar measures)
shows that 
$\delta(s)=\Delta_G(s)\Delta_{G/N}(s^{-1})$ for all $s\in G$.
Similar 
to 
Example \ref{ex-deco} we define an action 
$\alpha^N:G\to\Aut(A\rtimes 
N)$
by
\begin{equation}\label{eq-decom}
\big(\alpha^N_s(f)\big)(n)=\delta(s) 
\alpha_s\big(f(s^{-1}ns)\big)
\end{equation}
for $f$ in the dense 
subalgebra $C_c(N,A)\subseteq A\rtimes N$.
If we denote by 
$\tau^N:=i_N:N\to UM(A\rtimes N)$ the canonical embedding 
as defined 
in part (1) of Remark \ref{rem-universal}, then the pair 
$(\alpha^N,\tau^N)$ satisfies the 
equations
\begin{equation}\label{eq-twist}
\tau^N_n 
x\tau^N_{n^{-1}}=\alpha^N_n(x)\quad\text{and}\quad
\alpha^N_s(\tau^N_n)=\tau^N_{sns^{-1}}
\end{equation}
for 
all $x\in A\rtimes N$, $n\in N$ and $s\in G$, where in the second formula 
 we extended the 
automorphism
$\alpha^N_s$ of $A\rtimes N$ to $M(A\rtimes N)$. Suppose 
now that 
$(\pi,U)$ is a covariant homomorphism of $(A,G,\alpha)$ 
into some $M(D)$. 
Let $(\pi, U|_N)$ denote its restriction to 
$(A,N,\alpha)$ and let 
$\pi\times U|_N:A\rtimes N\to M(D)$ be its 
integrated form. Then $(\pi\times U|_N, 
U)$ is a non-degenerate covariant 
homomorphism of $(A\rtimes N, G,\alpha^N)$ which 
satisfies
$$\pi\times U|_N(\tau^N_n)=U_n$$
for all $n\in N$ (see Remark \ref{rem-universal}).
The 
pair $(\alpha^N,\tau^N)$ is the prototype for a twisted action 
(which we denote the {\em decomposition twisted action}) and 
$(\pi\times U|_N, U)$ is the prototype of a twisted covariant 
homomorphism 
as in

\begin{definition}[Green]\label{def-twist}
Let $N$ be 
a closed normal subgroup of $G$. A {\em twisted action} of $(G,N)$ 
on a $C^*$-algebra $A$ is a pair $(\alpha,\tau)$ such 
that
$\alpha:G\to \Aut(A)$ is an action and 
 $\tau:N\to UM(A)$ is a 
strictly continuous homomorphism such that
$$\tau_n 
a\tau_{n^{-1}}=\alpha_n(a)\quad\text{and}\quad
\alpha_s(\tau_n)=\tau_{sns^{-1}}
$$
for 
all $a\in A$, $n\in N$ and $s\in G$. We then say that 
$(A,G,N,\alpha,\tau)$ 
is a twisted system.
A {\em (twisted) 
covariant homomorphism}  of $(A,G,N,\alpha,\tau)$ into some $M(D)$ 
is a covariant homomorphism $(\rho,V)$ of $(A,G,\alpha)$ into $M(D)$ 
which {\em preserves $\tau$} in the sense that $\rho(\tau_na)=V_n\rho(a)$ for 
all $n\in N, a\in A$. \footnote{The latter condition becomes $\rho(\tau_n)=V_n$
if $(\rho,V)$ is non-degenerate.}
\end{definition}

\begin{remark}\label{rem-twist}
Note that the kernel of the regular representation $\Lambda_A^N:A\rtimes N\to 
A\rtimes_rN$ is easily seen to be invariant under the decomposition 
twisted action $(\alpha^N,\tau^N)$ (which just means that it is invariant 
under $\alpha^N$), so that  $(\alpha^N,\tau^N)$ induces a twisted action 
on the quotient $A\rtimes_rN$. In what follows, we shall make no notational difference
between the decomposition twisted actions on the full or the 
reduced crossed products.
\end{remark}

Let $C_c(G,A,\tau)$ denote the set 
of all 
continuous $A$-valued functions on $G$ with compact support mod 
$N$
and which satisfy
$$f(ns)=f(s)\tau_{n^{-1}}\quad\quad\text{for 
all $n\in N$, $s\in G$.}$$
Then $C_c(G,A,\tau)$ becomes a $*$-algebra 
with convolution and 
involution defined by
$$f*g(s)=\int_{G/N} 
f(t)\alpha_t(g(t^{-1}s))\,dtN\quad\text{and}
\quad 
f^*(s)=\Delta_{G/N}(s^{-1})\alpha_s\big(f(s^{-1})^*\big).$$
If 
$(\rho,V)$ is a covariant representation of $(A,G,N,\alpha,\tau)$, 
then the equation
$$\rho\times V(f)=\int_{G/N} \rho(f(s))V_s\, 
dsN$$
defines a $*$-homomorphism $\rho\times V: C_c(G,A,\tau)\to 
M(D)$, and the 
{\em full twisted crossed product} 
$A\rtimes_{\alpha,\tau}(G,N)$ (or just 
$A\rtimes (G,N)$ if 
$(\alpha,\tau)$ is understood) is defined as the 
completion of 
$C_c(G,A,\tau)$ with respect to 
$$\|f\|_{\max}:=\sup\{\|\rho\times 
V(f)\|: \text{$(\rho,V)$ is a covariant homomorphism 
of 
$(A,G,N,\alpha,\tau)$}\}.$$
Note that the same 
formulas as given in Remark \ref{rem-universal}
define a twisted covariant homomorphism
$(i_A, i_G)$ of 
$(A,G,N,\alpha,\tau)$ into $M(A\rtimes (G,N))$ such that
any 
non-degenerate homomorphism $\Phi:A\rtimes (G,N)\to M(D)$ is the 
integrated form $\rho\times V$ with $\rho=\Phi\circ i_A$ and 
$V=\Phi\circ i_G$.

\begin{remark}\label{rem-quotient}
It is important to notice that for any twisted action $(\alpha,\tau)$
of $(G,N)$ the map
$$\Phi: C_c(G,A)\to C_c(G,A,\tau); \Phi(f)(s)=\int_N f(sn)\tau_{sns^{-1}}\,dn$$
extends to a quotient map $\Phi:A\rtimes G\to A\rtimes (G,N)$
of the full crossed products, such that
$\ker\Phi=\cap\{\ker(\pi\times U): (\pi, U)\;\text{preserves $\tau$}\}.$
The ideal $I_{\tau}:=\ker\Phi$ is called the {\em twisting ideal}
of $A\rtimes G$. Note that if $G=N$, then 
$A\rtimes (N,N)\cong A$ via $f\mapsto f(e); C_c(N,A,\tau)\to A$.
\end{remark}

For the definition of the  {\em reduced twisted 
crossed products} 
$A\rtimes_{\alpha,\tau,r}(G,N)$ (or just 
$A\rtimes_r(G,N)$)
we define a Hilbert $A$-module 
$L^2(G,A,\tau)$ by 
taking the completion of $C_c(G,A,\tau)$ with respect to 
the 
$A$-valued inner product
$$\lk 
\xi,\eta\rk_A:=\xi^**\eta(e)=\int_{G/N}\alpha_{s^{-1}}\big(\xi(s)^*\eta(s)\big)\,dsN.$$
The 
{\em regular representation} 
$$\Lambda_A^{G,N}:C_c(G,A,\tau)\to 
\L_A(L^2(G,A,\tau)); \;\;\Lambda_A^{G,N}(f)\xi=f*\xi$$
embeds 
$C_c(G,A,\tau)$ into the algebra of adjointable operators on 
$L^2(G,A,\tau)$
and then 
$A\rtimes_r(G,N):=\overline{\Lambda_A^{G,N}\big(C_c(G,A,\tau)\big)}
\subseteq 
\L_A(L^2(G,A,\tau))$. 
If $N=\{e\}$ is trivial, then $\L_A(L^2(G,A))$ 
identifies naturally
with $M(A\otimes \K(L^2(G)))$, and we recover the 
original definition of the 
regular representation $\Lambda_A^G$ of 
$(A,G,\alpha)$ and 
of the reduced crossed product $A\rtimes_rG$ of 
$A$ by $G$.

\begin{remark}\label{rem-quotred}
The analogue of Remark \ref{rem-quotient} {\bf does not hold} 
in general for the reduced crossed products, i.e. $A\rtimes_r(G,N)$
is in general not a quotient of $A\rtimes_rG$. 
For example, if $N$ is not amenable, the 
algebra $C_r^*(G/N)=\CC\rtimes_{\id,1,r}(G,N)$ is not a quotient
of $C_r^*(G)=\CC\rtimes_{\id,r}G$ -- at least not in a canonical way. 
\end{remark}

We are now coming back to the decomposition 
problem

\begin{proposition}[Green]\label{prop-decom}
Let $\alpha:G\to 
\Aut(A)$ be an action, let $N$ be a closed normal subgroup 
of $G$,
and let $(\alpha^N,\tau^N)$ be the decomposition twisted action of 
$(G,N)$ on 
$A\rtimes N$. Then the 
map
\begin{equation}\label{eq-isodecom}
\Psi:C_c(G,A)\to 
C_c(G,C_c(N,A), 
\tau^N);\;\;\;\Psi(f)(s)(n)=\delta(s)f(ns)
\end{equation}
extends to 
isomorphisms $A\rtimes G\cong (A\rtimes N)\rtimes (G,N)$
and 
$A\rtimes_rG\cong (A\rtimes_r N)\rtimes_r (G,N)$. In particular, if 
$A=\CC$ we obtain isomorphisms  $C^*(G)\cong C^*(N)\rtimes (G,N)$
and 
$C_r^*(G)\cong C_r^*(N)\rtimes_r(G,N)$.
Under the isomorphism 
of the 
full crossed products, a representation $\pi\times U$ of $A\rtimes 
G$
corresponds to the representation $(\pi\times U|_N)\times U$ of 
$(A\rtimes N)\rtimes (G,N)$.
\end{proposition}

A similar result 
holds if we start with a twisted action of $(G,M)$ on $A$ with 
$M\subseteq N$ (see \cite{Green1}*{Proposition 1} and \cite{CE1}).
We should note that Green only considered  full crossed products in \cite{Green1}.
The above decomposition of reduced crossed products was first shown 
by Kirchberg and Wassermann in \cite{KW2}.
Note that all results 
stated in \S \ref{sec-actions} for ordinary crossed products 
have their complete analogues in the twisted case, where $G/N$ plays 
the r\^ole of $G$. In particular, the 
full and reduced crossed 
products coincide if $G/N$ is amenable.
Indeed, we shall see in \S \ref{sec-twistmor} below that there is a convenient way to 
extend results known for ordinary actions 
to the twisted case via Morita equivalence (see Theorem \ref{thm-twistedequiv}
below).

\subsection{The twisted equivariant corespondence category and the stabilisation trick}\label{sec-twistmor}

As done for ordinary actions in \S\ref{sec-mor}  we may consider the 
{ \em  twisted equivariant correspondence category} $\Corr(G,N)$ (resp. the {\em compact  twisted equivariant correspondence category}
$ \Corr_c(G,N)$) in which the objects are
twisted systems $(A,G,N,\alpha,\tau)$ and in which the morphism from 
$(A,G,N, \alpha,\tau)$ and $(B,G, N, \beta,\sigma)$ are given by morphisms
$[E,\Phi,u]$ from $(A,G,\alpha)$ to $(B,G,\beta)$ in $\Corr(G)$ (resp. $\Corr_c(G)$) which preserve 
the twists in the sense that
\begin{equation}\label{eq-mortwist}
\Phi(\tau_n)\xi= u_n(\xi)\sigma_n\quad\text{for all $n\in N$.}
\end{equation}
As for ordinary actions, the  crossed product construction
$(A,G,N,\alpha,\tau)\mapsto A\rtimes_{(r)}(G,N)$
extend to (full and reduced) descent functors
$$\rtimes_{(r)}:\Corr(G,N)\to \Corr.$$
If $[E,\Phi, u]\in \Mor(G,N)$ is a morphism from 
$(A,G,N, \alpha,\tau)$ to $(B,G, N, \beta,\sigma)$, then the descent
$[E\rtimes_{(r)} (G,N),\Phi\rtimes_{(r)}(G,N)]$ can be defined  
by setting $E\rtimes_{(r)} (G,N):=(E\rtimes G)/\big((E\rtimes G)\cdot I_{(r)}\big)$ 
with $I_{(r)}:=\ker\big(B\rtimes G\to B\rtimes_{(r)}(G,N)\big)$. 
Alternatively, one can construct $E\rtimes_{(r)}G$ as the closure of 
$C_c(G,E,\sigma)$, the continuous $E$-valued functions $\xi$ on $G$ with compact support
modulo $N$ satisfying $\xi(ns)=\xi(s)\sigma_{n^{-1}}$ for $s\in G$, $n\in N$, with respect 
to the $B\rtimes_{(r)}(G,N)$-valued inner product given by
$$\lk \xi, \eta\rk_{B\rtimes_{(r)}(G,N)}(t)=\int_{G/N}\beta_{s^{-1}}(\lk\xi(s), \eta(ts)\rk_B)\,dsN$$
(compare with the formulas given in \S\ref{subsec-equivmor}).

There is a natural inclusion functor $\inf: \Corr(G/N)\to \Corr(G,N)$ given as follows:
If $(A, G/N,\alpha)$ is an action of $G/N$, we let $\inf\alpha:G\to \Aut(A)$ 
denote the inflation of $\alpha$ from $G/N$ to $G$ and we let
$1_N:N\to U(A)$ denote the trivial homomorphism $1_N(s)=1$.
Then $(\inf\alpha, 1_N)$ is a twisted action of $(G,N)$ on $A$
and we set
$$\inf\big((A,G/N,\alpha)\big):=(A,G,N,\inf\alpha, 1_N).$$
Similarly, on morphisms we set $\inf\big([E,\Phi,u]\big):=[E,\Phi,\inf u]$,
where $\inf u$ denotes the inflation of $u$ from $G/N$ to $G$.
The dense subalgebra $C_c(G,A,1_N)$ of the crossed product 
$A\rtimes_{(r)}(G,N)$ for $(\inf\alpha, 1_N)$ consists of functions 
which are constant on $N$-cosets and which have compact supports in $G/N$, 
hence it coincides with $C_c(G/N,A)$ (even as a $*$-algebra). The 
identification $C_c(G,A,1_N)\cong C_c(G/N,A)$ extends to the crossed 
products, and we 
obtain canonical isomorphisms $A\rtimes_{(r)} G/N\cong 
A\rtimes_{(r)}(G,N)$. A similar observation can be made for the crossed products
of morphism and we see that the inclusion $\inf:\Corr(G/N)\to \Corr(G,N)$ 
is compatible with the crossed product functor in the sense that the diagram
$$
\begin{CD}
\Corr(G/N)   @>\inf>> \Corr(G,N)\\
@V\rtimes VV       @VV\rtimes V\\
\Corr  @= \Corr
\end{CD}
$$ 
commutes. 

In what follows next we want to see that every twisted action is Morita equivalent
(and hence isomorphic in $\Corr(G,N)$) to some inflated twisted action as above.
This will allow us to pass to an untwisted system whenever a theory 
(like the theory of induced representations, or $K$-theory of crossed products, etc.)
only  depends on the Morita equivalence class of a given twisted action.

To do this, we first note that 
Green's imprimitivity theorem (see Theorem \ref{thm-Green}) extends easily
to crossed products by twisted actions:
If $N$ is a closed normal subgroup of $G$ such that $N\subseteq H$
for some closed subgroup $H$ of $G$,
and if $(\alpha,\tau)$ is a twisted action of $(H,N)$ 
on $A$, then we obtain a twisted action 
$(\Ind\alpha,\Ind\tau)$ of $(G,N)$ on $\Ind_H^G(A,\alpha)$
by defining 
$$(\Ind\tau_n f)(s)=\tau_{s^{-1}ns}f(s)\quad\text{ for}\;f\in \Ind A, s\in G\;\text{and}\; n\in N.$$
One can check that the twisting ideals 
$I_{\tau}\subseteq A\rtimes H$ and $I_{\Ind\tau}\subseteq 
\Ind A\rtimes G$ (see Remark \ref{rem-quotient})
are linked via the Rieffel correspondence  of the 
$\Ind A\rtimes G - A\rtimes H$ imprimitivity bimodule $X_H^G(A)$.
Similarly,  the kernels 
$I_{\tau,r}:=\ker\big(A\rtimes H\to A\rtimes_r(H,N)\big)$ and 
$I_{\Ind\tau,r}:=\ker\big(\Ind A\rtimes G\to \Ind A\rtimes_r(G,N)\big)$ are linked 
via the Rieffel correspondence
(we refer to \cite{Green1} and  \cite{KW2} for the details). 
Thus,  from  Proposition \ref{prop-rief-cor} it follows:

\begin{theorem}\label{thm-greentwist}
 The quotient $Y_H^G(A):=X_H^G(A)/(X_H^G(A)\cdot I_{\tau}) $ (resp. 
 $Y_H^G(A)_r:=X_H^G(A)/(X_H^G(A)\cdot I_{\tau,r})$) 
becomes an $\Ind_H^G(A,\alpha)\rtimes(G,N)$-$A\rtimes (H,N)$ (resp. $\Ind_H^G(A,\alpha)\rtimes_r(G,N)$ - $A\rtimes_r(H,N)$)
imprimitivity bimodule.
\end{theorem}

\begin{remark}\label{rem-resttwistgreen}
{\bf (1)} Alternatively, one can construct the modules $Y_H^G(A)$ and $Y_H^G(A)_r$
by taking completions of $Y_0(A):=C_c(G,A,\tau)$
with respect to suitable
$C_c(G,\Ind A,\Ind\tau)$- and $C_c(N,A,\tau)$-valued 
inner products. The formulas are precisely those of (\ref{eq-products})
if we integrate over $G/N$ and $H/N$, respectively (compare with the formula 
for convolution in $C_c(G,A,\tau)$ as given in \S\ref{sec-twists}).

{\bf (2)} If we start with a twisted action $(\alpha,\tau)$ of $(G,N)$ on $A$ and restrict this to $(H,N)$,
then the induced algebra $\Ind_H^G(A,\alpha)$ is isomorphic to $C_0(G/H,A)\cong C_0(G/H)\otimes A$
 as in Remark \ref{rem-Ind}. The isomorphism transforms the action $\Ind\alpha$ to the action 
$l\otimes\alpha:G\to \Aut(C_0(G/H,A))$, with $l:G\to \Aut(C_0(G/H))$ being left-translation action,
and the twist $\Ind\tau$ is transformed to the twist $1\otimes\tau:N\to U(C_0(G/H)\otimes A)$.
Hence, in this setting, the above theorem provides Morita equivalences 
$$A\rtimes_{(r)}(H,N)\sim_M C_0(G/H,A)\rtimes_{(r)}(G,N)$$
for the above described twisted action $(l\otimes\alpha, 1\otimes \tau)$ of $(G,N)$.
\end{remark}

We want to use Theorem \ref{thm-greentwist} to construct a functor 
$$\F:\Corr(G,N)\to \Corr(G/N)$$ which, up to a
natural equivalence, inverts the inflation functor $\inf:\Corr(G/N)\to \Corr(G,N)$.
We start with the special case of the decomposition twisted actions 
$(\alpha^N,\tau^N)$ of $(G,N)$ on $A\rtimes N$ with respect to 
a given system $(A,G,\alpha)$ and a normal subgroup $N$ of $G$
(see \S\ref{sec-twists} for the construction).
Since $A$ is a $G$-algebra, it follows from 
Remark \ref{rem-Ind} that $\Ind_N^G(A,\alpha)$ is isomorphic to 
$C_0(G/N,A)$ as a $G$-algebra. 
Let $X_N^G(A)$ be Green's $C_0(G/N,A)\rtimes G-A\rtimes N$
imprimitivity bimodule. Since right translation of $G/N$ 
on $C_0(G/N,A)$ commutes with $\Ind\alpha$, 
it induces an action 
$$\beta^N:G/N\to \Aut\big(C_0(G/N,A)\rtimes G\big)$$
on the crossed product. 
For $s\in G$ and $\xi\in C_c(G,A)\subseteq X_N^G(A)$ let
$$u_s^N (\xi)(t):=\sqrt{\delta(s)}\alpha_s(\xi(ts)),\quad \xi\in C_c(G,A)$$
where $\delta(s)=\Delta_G(s)\Delta_{G/N}(s^{-1})$. 
This formula determines an action 
$u^N:G\to  \Aut(X_N^G(A))$ such that 
$(X_N^G(A),u^N)$ becomes a $(G,N)$-equivariant 
$C_0(G/N,A)\rtimes G$-$A\rtimes N$ Morita equivalence
with respect to the twisted actions $(\inf\beta^N, 1_N)$ and $(\alpha^N,\tau^N)$, respectively. 
All these twisted actions pass to the quotients to 
give also a $(G,N)$-equivariant
equivalence $(X_N^G(A)_r, u^N)$ 
for the reduced crossed products.
Thus we get

\begin{proposition}[{cf. \cite{E-mor}*{Theorem 1}}]\label{prop-twistequiv}
The decomposition action $(\alpha^N,\tau^N)$ of $(G,N)$ on 
$A\rtimes_{(r)}N$ is canonically Morita equivalent to the 
(untwisted) action
$\beta^N$ of $G/N$ on $C_0(G/N,A)\rtimes_{(r)}G$ as described
above.
\end{proposition}

If one starts with an arbitrary twisted action $(\alpha,\tau)$ 
of $(G,N)$ on $A$, one checks that the twisting ideals 
$I_{\tau}\subseteq A\rtimes N$ and $I_{\Ind\tau}\subseteq 
C_0(G/N,A)\rtimes G$ are $(G,N)$-invariant and  that the twisted
 action on $A\cong (A\rtimes N)/I_{\tau}$
(cf. Remark \ref{rem-quotient}) induced 
from $(\alpha^N,\tau^N)$ is equal to 
$(\alpha,\tau)$. Hence, if $\beta$ denotes the 
action of $G/N$ on $C_0(G/N,A)\rtimes (G,N)\cong \big(C_0(G/N,A)\rtimes G\big)/I_{\Ind\tau}$ induced from $\beta^N$, then 
$u^N$ factors through an action $u$ of $G$ on 
$Y_N^G(A)=X_N^G(A)/(X_N^G(A)\cdot I_{\tau})$ such that
$(Y_N^G(A), u)$ becomes a $(G,N)$-equivariant 
$C_0(G/N,A)\rtimes (G,N)$-$A$ Morita equivalence with respect to the twisted
actions $(\inf \beta, 1_N)$ and $(\alpha,\tau)$, respectively.
Following the arguments given in 
 \cite{EKQR1} one can show  that there is a functor
$\F: \Corr(G,N)\to \Corr(G/N)$ given on objects by the assignment
$$(A,G,N,\alpha,\tau)\stackrel{\F}{\mapsto} 
(C_0(G/N,A)\rtimes(G,N),\, G/N,\, \beta)$$
(and a similar crossed-product construction on the morphisms)
such that

\begin{theorem}[{cf. \cite{E-mor}*{Theorem 1} and \cite{EKQR1}*{Theorem 4.1}}]
\label{thm-twistedequiv}
The assignment 
$$(A,G,N,\alpha,\tau)\mapsto (Y_N^G(A),u)$$
 is a natural equivalence between the identity functor
on $\Corr(G,N)$ and the functor $\inf\circ \F:\Corr(G,N)\to \Corr(G,N)$,
where $\inf:\Corr(G/N)\to \Corr(G,N)$ denotes the inflation functor.
In particular, every twisted action of $(G,N)$ is Morita equivalent
to an ordinary action of $G/N$ (viewed as a twisted action via inflation).
\end{theorem}

Note that a first 
version of the above Theorem was obtained by 
Packer and Raeburn in the setting
of Busby-Smith twisted actions  (\cite{PR89}). We therefore call it 
the {\em Packer-Raeburn stabilisation trick}.
As mentioned before, it allows to extend results known for ordinary actions 
to the twisted case as soon as they are invariant under 
Morita equivalence.
If $A$ is separable and $G$ is second countable,
the algebra $B=C_0(G/N,A)\rtimes (G,N)$ is separable, too. 
Thus, it follows from a theorem of
Brown, Green, and Rieffel (see \cite{BGR}) that $A$ and $B$ are stably isomorphic
(a direct isomorphism $B\cong A\otimes \K(L^2(G/N))$ is obtained 
in \cite{Green2} but  see also  Proposition \ref{prop-Green-structure}).
Hence, as a consequence of  Theorem \ref{thm-twistedequiv} we get

\begin{corollary}\label{cor-twistedequiv}
If $G$ is second countable and  $A$ is separable, then
every twisted action of $(G,N)$ on $A$ is Morita equivalent 
to some action $\beta$ of $G/N$ on $A\otimes\K$.
\end{corollary}

We want to discuss some further consequences of Theorem \ref{thm-twistedequiv}:

\subsection{Twisted Takesaki-Takai duality}
If $(A,G,N,\alpha,\tau)$ is a twisted system
with $G/N$ abelian, then we define the dual action 
$$\widehat{(\alpha,\tau)}:\widehat{G/N}\to \Aut\big(A\rtimes(G,N)\big)$$
as in the previous section by pointwise multiplying characters
of $G/N$ with functions 
in the dense subalgebra $C_c(G,A,\tau)$. 
Similarly, we can define actions of $\widehat{G/N}$ on 
(twisted) crossed products of Hilbert bimodules, so that taking dual actions
gives a descent functor
$\rtimes :\Corr(G,N)\to \Corr(\widehat{G/N})$.
The Takesaki-Takai duality theorem shows that on 
$\Corr(G/N)\subseteq \Corr(G,N)$ this functor
is inverted, up to a natural equivalence, by the functor 
$\rtimes:\Corr(\widehat{G/N})\to \Corr(G/N)$. 
Using Theorem \ref{thm-twistedequiv},
this directly extends to the twisted case. 

\subsection{Stability of exactness under group extensions}
Recall from \S\ref{subsec-permanence} that a group is called exact if for every short exact sequence 
$0\to I\to A\to A/I\to 0$ of $G$-algebras the resulting sequence
$$0\to I\rtimes_rG\to A\rtimes_rG\to (A/I)\rtimes_rG\to 0$$
of reduced crossed products is exact. We want to use Theorem \ref{thm-twistedequiv}
to give a proof of the following result of Kirchberg and Wassermann:

\begin{theorem}[{Kirchberg and  Wassermann \cite{KW2}}]\label{thm-extexact}
Suppose that $N$ is a closed normal subgroup of the locally compact group
$G$ such that $N$ and $G/N$ are exact. Then $G$ is exact.
\end{theorem}

The result will follow  from 

\begin{lemma}\label{lem-exactmor} Suppose that $N$ is a closed normal subgroup of $G$
and  that $(X,u)$ is a $(G,N)$-equivariant 
Morita equivalence for the twisted actions
 $(\beta,\sigma)$ and $(\alpha,\tau)$ of $G$ on $B$ and $A$, respectively.
Let $I\subseteq A$ be a $(G,N)$-invariant ideal of $A$, and 
let  $J:=\Ind^XI\subseteq B$ denote the ideal of $B$ induced from $I$ via $X$ (which is 
a $(G,N)$-equivariant ideal of $B$).

Then $J\rtimes_{(r)}(G,N)$ (resp. $(B/J)\rtimes_{(r)}(G,N)$) 
corresponds to $I\rtimes_{(r)}(G,N)$ (resp. $(A/I)\rtimes_{(r)}(G,N)$) under the 
Rieffel correspondence for $X\rtimes_{(r)}(G,N)$.
\end{lemma}
\begin{proof} Let $Y:=X\cdot I\subseteq X$. Then  the closure
$C_c(G,Y,\tau)\subseteq C_c(G,X,\tau)$ is a \mbox{$B\rtimes_{(r)}(G,N)-A\rtimes_{(r)}(G,N)$}
submodule of $X\rtimes_{(r)}(G,N)$ which corresponds to the ideals
$J\rtimes_{(r)}(G,N)$ and $I\rtimes_{(r)}(G,N)$ under the Rieffel correspondence.
For the quotients observe that the obvious quotient map
$C_c(G,X,\tau)\to C_c(G,X/X\cdot I, \tau)$ extends to an imprimitivity bimodule quotient map
\mbox{$X\rtimes_{(r)}(G,N)\to (X/X\cdot I)\rtimes_{(r)}(G,N)$,} whose kernel 
corresponds to the ideals $K_B:=\ker\big(B\rtimes_{(r)}(G,N)\to (B/J)\rtimes_{(r)}(G,N)\big)$
and $K_A:=\ker\big(A\rtimes_{(r)}(G,N)\to (A/I)\rtimes_{(r)}(G,N)\big)$ under the 
Rieffel correspondence (see Remark \ref{rem-imphom}).
\end{proof}

As a consequence  we get

\begin{lemma}\label{lem-exacttwist}
Suppose that $N$ is a closed normal subgroup of $G$ such that $G/N$ is exact.
Suppose further that $0\to I\to A\to A/I\to 0$ is a short exact sequence of 
$(G,N)$-algebras. Then the sequence
$$0\to I\rtimes_r(G,N)\to A\rtimes_r(G,N)\to (A/I)\rtimes_r(G,N)\to 0$$
is exact.
\end{lemma}
\begin{proof} By Theorem \ref{thm-twistedequiv} there exists a system $(B,G/N, \beta)$
such that\linebreak
 $(B, G,N,\inf\beta,1_N)$ is Morita equivalent to the given twisted system
$(A,G,N,\alpha,\tau)$ via some equivalence $(X,u)$. 
If $I$ is a $(G,N)$-invariant ideal of $A$, let $J:=\Ind^XI\subseteq B$. It follows then from
Lemma \ref{lem-exactmor} and the Rieffel correspondence, that
$$0\to I\rtimes_r(G,N)\to A\rtimes_r(G,N)\to (A/I)\rtimes_r(G,N)\to 0$$
is exact if and only if
$$0\to J\rtimes_r(G,N)\to B\rtimes_r(G,N)\to (B/J)\rtimes_r(G,N)\to 0$$
is exact. But the latter sequence is equal to the sequence
$$0\to J\rtimes_rG/N\to B\rtimes_rG/N\to (B/J)\rtimes_rG/N\to 0,$$ 
which is exact  since $G/N$ is exact.
\end{proof}

\begin{proof}[Proof of Theorem \ref{thm-extexact}]
Suppose that $0\to I\to A\to A/I\to 0$ is an exact sequence of $G$-algebras and consider
the decomposition twisted action $(\alpha^N,\tau^N)$ of $(G,N)$ on $A\rtimes_rN$.
Since $N$ is exact, the sequence
$$0\to I\rtimes_rN\to A\rtimes_rN\to (A/I)\rtimes_rN\to 0$$
is a short exact sequence of $(G,N)$-algebras. Since $G/N$ is exact, it follows therefore
from Lemma \ref{lem-exacttwist} that
$$0\to (I\rtimes_rN)\rtimes_r(G,N)\to (A\rtimes_rN)\rtimes_r(G,N)\to \big((A/I)\rtimes_rN\big)\rtimes_r(G,N)
\to 0$$
is exact. But it follows from Proposition \ref{prop-decom} that this sequence equals
$$0\to I\rtimes_rG\to A\rtimes_rG\to (A/I)\rtimes_rG\to 0.$$
\end{proof}

\subsection{Induced representations of twisted crossed products}\label{subsec-indtwist}

Using Green's imprimitivity theorem for twisted systems, we can define induced representations
and ideals for twisted crossed products $A\rtimes (G,N)$ as in the untwisted case, using 
the spaces $C_c(G,A,\tau)$ and $C_c(G, \Ind A, \Ind\tau)$ etc. (e.g. see 
\cite{E-ct}*{Chapter 1} for this approach). An alternative but equivalent way, as followed 
in Green's original paper \cite{Green1} is to define induced representations via
the untwisted crossed products: Suppose that $(\alpha,\tau)$ is a twisted action 
of $(G,N)$ on $A$ and let $H\subseteq G$ be a closed subgroup of $G$ such that
$N\subseteq H$. Since $A\rtimes (H,N)$ is a quotient of $A\rtimes H$ we can 
regard every representation of $A\rtimes (H,N)$ as a representation of $A\rtimes H$.
We can use the untwisted theory to induce the representation to $A\rtimes G$.
But then we have to check that this representation factors through the quotient 
$A\rtimes (G,N)$ to have a satisfying theory. This has been done in 
\cite{Green1}*{Corollary 5}, but one can also obtain it as an easy consequence of 
Proposition \ref{prop-ind-res}: Let $I_\tau^N\subset A\rtimes N$ denote the twisting 
ideal for $(\alpha|_N,\tau)$. It is then clear from the definition of representations $\pi\times U$
of $A\rtimes H$ (resp. $A\rtimes G$) which preserve $\tau$, that $\pi\times U$ preserves
$\tau$ iff $\pi\times U|_N$ preserves $\tau$ as a representation of $A\rtimes N$.
Hence, $\pi \times U$ is a representation of $A\rtimes (H,N)$ (resp. $A\rtimes (G,N)$)
iff $I_\tau^N\subseteq \ker(\pi\times U|_N)$. Since $I_\tau^N$ is easily seen to 
be a $G$-invariant ideal of $A\rtimes N$, this property is preserved under induction by 
Proposition \ref{prop-ind-res}.

The procedure of inducing representations is compatible with
passing to Morita equivalent systems. To be more precise:
Suppose that $(X,u)$ is a Morita equivalence for the systems 
$(A,G,\alpha)$ and $(B,G,\beta)$. If $H$ is a closed subgroup of $G$ and 
$\pi\times U$ is a representation of $B\rtimes H$, then we get an equivalence
$$\Ind_H^G\big(\Ind^{X\rtimes H}(\pi\times U)\big)\cong \Ind^{X\rtimes G}\big(\Ind_H^G(\pi\times U)\big).
$$
This result follows from an isomorphism of $A\rtimes G$-$B\rtimes H$ bimodules
$$X_H^G(A)\otimes_{A\times H}(X\rtimes H)\cong (X\rtimes G)\otimes_{B\rtimes G}X_H^G(B),$$
which just means that the respective compositions in the correspondence categories coincide.
A similar result can be shown for the reduction of representations to subgroups.
Both results will follow from a linking algebra trick as introduced in \cite{ER}*{\S4}.
Similar statements holds for twisted systems.

  \subsection{Twisted group algebras, actions on $\K$ and 
  Mackey's little group method}\label{subsec-littlegroup}

In this section we want to study crossed products of the form $\K\rtimes_{(r)}G$, 
where $\K=\K(H)$  is the algebra of compact operators on some 
Hilbert space $H$. As we shall see below, such actions are strongly related to twisted 
actions on the algebra $\CC$ of complex numbers.
While there are only trivial actions of groups on 
$\CC$, there are usually many nontrivial twisted
actions of pairs $(G,N)$ on $\CC$. However, in a certain sense they 
are all equivalent to twisted actions of the following type:

\begin{example}\label{ex-om} Assume that $1\to \TT\to \tilde{G}\to G\to 1$ 
is a central extension of the locally compact group $G$ by the 
circle group $\TT$. Let $\iota:\TT\to \TT; \iota(z)=z$ denote the identity character 
on $\TT$. Then $(\id,\iota)$ is a twisted action of $(\tilde{G},\TT)$ on $\CC$.
A (covariant) representation of the twisted system $(\CC, \tilde{G},\TT,\id,\iota)$
on a Hilbert space $H$ 
consists of  the  representation $\lambda\mapsto \lambda 1_H$ of $\CC$ together with a unitary representation $U:\tilde{G}\to \U(H)$ satisfying $U_z=z\cdot 1_H$ for
all $z\in \TT\subseteq \tilde{G}$, i.e., of those representations of $\tilde{G}$ 
which restrict to a multiple of $\iota$ on the central subgroup $\TT$ of $\tilde{G}$.
Hence, the twisted crossed product $\CC\rtimes(\tilde{G},\TT)$ is the 
quotient of $C^*(\tilde{G})$ by the ideal 
$I_{\iota}=\cap\{\ker U: U\in \Rep(\tilde{G})\;\text{and $U|_\TT=\iota\cdot 1_H$}\}$.
Note that the isomorphism class of $\CC\rtimes(\tilde{G},\TT)$ only depends on the 
isomorphism class of the extension $1\to \TT\to \tilde{G}\to G\to 1$.

If $G$ is second countable\footnote{This assumptions is made to avoid measurability problems. 
With some extra care, much of the following discussion also works in the non-separable case (e.g. see \cite{Kl})}, we can choose a Borel section $c:G\to\tilde{G}$ in the above extension,
and we then obtain a Borel map $\om:G\times G\to \TT$ by
$$\om(s,t):=c(s)c(t)c(st)^{-1}\in \TT.$$
A short computation then shows that $\om$ satisfies the cocycle conditions
$\om(s,e)=\om(e,s)=1$ 
and $\om(s,t)\om(st,r)=\om(s,tr)\om(t,r)$ for all $s,t,r\in G$. Hence it is 
a 2-cocycle in $Z^2(G,\TT)$ of Moore's group cohomology with Borel cochains
 (see \cites{Mo1,Mo2,Mo3,Mo4}). 
 The cohomology
class $[\om]\in H^2(G,\TT)$ then only depends on the isomorphism class of the 
given extension $1\to\TT\to\tilde{G}\to G\to 1$.\footnote{Two cocycles $\om$ and $\om'$ are
in the same class in $H^2(G,\TT)$
 iff they differ by a boundary $\partial f(s,t):=f(s)f(t)\overline{f(st)}$ of some 
Borel function $f:G\to \TT$.}
Conversely, if  $\om:G\times G\to \TT$ is any Borel $2$-cocycle on $G$, 
let $G_{\om}$ denote the cartesian product
$G\times \TT$ with multiplication given by
$$(s,z)\cdot (t,w)=(st, \om(s,t)zw).$$
By \cite{Ma4} there exists a unique 
locally compact topology on $G_{\om}$ whose
 Borel structure coincides with the product Borel structure.
Then $G_{\om}$ is a central extension of $G$ by $\TT$ corresponding to $\om$
(just consider the section $c:G\to G_\om; c(s)=(s,1)$) and  we obtain a complete
classification of the (isomorphism classes of) 
central extensions of $G$ by $\TT$ in terms of $H^2(G,\TT)$.
We then write $C^*_{(r)}(G,\om):=\CC\rtimes_{(r)}(G_\om,\TT)$ for the corresponding
full (resp. reduced) twisted crossed products, which we now call the
(full or reduced)  {\em twisted group algebra of $G$ corrsponding to $\om$}.

There is a canonical one-to-one correspondence
between the (non-degenerate) covariant representations of the twisted system $(\CC,G_{\om},\TT,\id,\iota)$ on a Hilbert space $H$ and the {\em projective
$\om$-representations} of $G$ on $H$, which are defined as Borel maps
$$V:G\to\U(H)\quad\text{satisfying}\quad V_sV_t=\om(s,t)V_{st}\quad s,t\in G.$$
Indeed, if $\tilde{V}:G_\om\to U(H)$ is a unitary representation of $G_\om$ which 
restricts to a multiple of $\iota$ on $\TT$, then $V_s:=\tilde{V}(s, 1)$ is the corresponding
$\om$-representation of $G$.

A convenient alternative realization of the twisted group algebra $C^*(G,\om)$ is obtained 
by taking a completion of the convolution algebra 
$L^1(G,\om)$, where $L^1(G,\om)$ denotes the algebra of all $L^1$-functions on $G$ with 
convolution and involution given by 
$$f*g(s)=\int_Gf(t)g(t^{-1}s)\om(t,t^{-1}s)\,dt \quad\text{and}\quad
f^*(s)=\Delta_G(s^{-1})\overline{\om(s,s^{-1})}\overline{f(s^{-1})}.$$
One checks that the $*$-representations of $L^1(G,\om)$ are given by integrating 
projective $\om$-representations and hence the corresponding $C^*$-norm for
completing  $L^1(G,\om)$ to $C^*(G,\om)$ is given by 
$$\|f\|_{\max}=\sup\{\|V(f)\|: \text{$V$ is an $\om$-representation of $G$}\}.$$
The map 
$$\Phi: C_c(G_{\om},\CC,\iota)\to L^1(G,\om);\Phi(f)(s):=f(s,1)$$
then extends to an isomorphism between the two pictures of $C^*(G,\om)$.\footnote{Use the identity $\overline{\om(t,t^{-1})}\om(t^{-1},s)\om(t,t^{-1}s)=1$
in order to check that $\Phi$ preserves multiplication.}
Similarly, we can define a {\em left  $\om$-regular representation} 
$L_\om$ of $G$ on $L^2(G)$ by setting
$$\big(L_\om(s)\xi\big)(t)=\om(s, s^{-1}t)\xi(s^{-1}t),\quad \xi\in L^2(G),$$ and then realize 
$C_r^*(G,\om)$ as $\overline{L_{\om}\big(L^1(G,\om)\big)}\subseteq \B(L^2(G))$.

%
%We shall see  below, that the theory of twisted group algebras for $G$
%is \mbox{(Morita-)} equivalent to the theory of crossed products $\K(H)\rtimes G$,
%where $\K(H)$ denotes the algebra
%of compact operators on a Hilbert space $H$.
\end{example}

\begin{example}\label{ex-tori}
Twisted group algebras appear quite often in $C^*$-algebra 
theory. For instance the rational and irrational rotation algebras 
$A_{\theta}$ for $\theta \in [0,1)$ are 
isomorphic to the twisted group algebras 
$C^*(\ZZ^2,\om_{\theta})$ with $\om_{\theta}\big((n,m), (k,l)\big)=e^{i2\pi\theta mk}$.
Note that every cocycle on $\ZZ^2$ is equivalent to $\om_\theta$ for some $\theta\in [0,1)$.
If $\theta=0$ we simply get $C^*(\ZZ^2)\cong C(\TT^2)$, the classical commutative $2$-torus.
For this reason the $A_\theta$ are often denoted as noncommutative $2$-tori. 

More generally, a 
 {\em noncommutative $n$-torus} is
a twisted group algebra $C^*(\ZZ^n,\om)$ for some cohomology class
 $[\om]\in H^2(\ZZ^n,\TT)$.
 
An extensive study of $2$-cocycles on abelian groups is given by Kleppner in \cite{Kl}.
In particular, for $G=\RR^n$, every cocycle is similar to a cocycle of the form
$\om(x,y)=e^{\pi i\lk Ax, y\rk}$, where $A$ is a skew-symmetric real $n\times n$-matrix,
and every cocycle of $\ZZ^n$ is similar to a restiction to $\ZZ^n$ of some
cocycle on $\RR^n$. 
The general structure of the twisted group algebras $C^*(G,\om)$ for abelian $G$ is studied
extensively in  \cite{ERo} in the type I case and in \cite{Po} in the general case.
If $G$ is abelian, then the {\em symmetry group} $S_{\om}$ of $\om$ is defined by
$$S_{\om}:=\{s\in G: \om(s,t)=\om(t,s)\;\text{for all $t\in G$}\}.$$
 Poguntke shows in \cite{Po} (in case $G$ satisfies some mild extra conditions, which are always satisfied
if $G$ is compactly generated) that $C^*(G,\om)$ is stably isomorphic to an algebra of the form
$C_0(\widehat{S}_\om)\otimes  C^*(\ZZ^n,\mu)$, where
$C^*(\ZZ^n,\mu)$ is some simple noncommutative $n$-torus (here we allow $n=0$ in which case we put   $C^*(\ZZ^n,\mu):=\CC$).\footnote{Two $C^*$-algebras $A$ and $B$ are called {\em stably isomorphic} if $A\otimes \K\cong B\otimes \K$, where $\K=\K(l^2(\NN))$.}
 \end{example}

 It follows from Theorems \ref{thm-MRG} and \ref{thm-GR} that for understanding 
the ideal structure of $A\rtimes G$, it is necessary to understand the structure
of $A_P\rtimes G_P$ for $P\in \Prim(A)$. In the special case $A=C_0(X)$, we saw
in the previous section that 
this is the same as understanding the group algebras $C^*(G_x)$ for the stabilisers
$G_x$, $x\in X$. In general, the problem becomes much more
difficult.
 However, at least in the important special case where $A$ is type I
 one can still give a quite satisfactory description of $A_P\rtimes G_P$ in terms of the 
 stabilisers. 
If $A$ is type I, we have $\widehat{A}\cong \Prim(A)$ via $\sigma\mapsto\ker\sigma$
and if $P=\ker\sigma$ for some $\sigma\in \widehat{A}$,
then the simple subquotient $A_P$ of $A$ corresponding to $P$
is isomorphic to $\K(H_{\sigma})$ (see Remark \ref{rem-typeImackey}). Thus, we have to understand the 
structure of the crossed products $\K(H_\sigma)\rtimes G_{\sigma}$, where 
$G_{\sigma}$ denotes the stabiliser of
$\sigma\in \widehat{A}$.

Hence, in what follows we shall always assume that $G$ is a locally compact group acting 
on the algeba $\K(H)$ of compact operators on some Hilbert space $H$.
In order to avoid measerability problems, we shall  always assume
that $G$ is second countable
 and  that $H$ is 
separable 
(see Remark \ref{rem-littlegroup1} for a brief discussion of the general case).
Since every automorphism of $\K(H)$ is given by conjugation 
with some unitary $U\in \B(H)$, it follows that  the automorphism group 
of $\K(H)$ is isomorphic (as topological groups) to the group $\PU:=\U/\TT\cdot 1$, 
where $\U=\U(H)$ denotes the group of
unitary operators on $H$ equipped with the strong operator topology.

Choose a Borel section $c:\PU\to\U$. 
If  $\alpha:G\to \PU$ is a continuous homomorphism, let $V_{\alpha}:=c\circ \alpha:G\to\U$.
Since $V_{\alpha}(s)V_{\alpha}(t)$ and $V_{\alpha}({st})$ both implement the automorphism 
$\alpha_{st}$, there exists a number $\om_{\alpha}(s,t)\in \TT$ with
$$\om_{\alpha}(s,t)\cdot 1=V_{\alpha}(st)V_{\alpha}(t)^*V_{\alpha}(s)^*.$$
A short computation using the identity $\Ad(V_\alpha(s)V_\alpha(tr))=\Ad(V_\alpha(s)V_\alpha(t)V_\alpha(r))=\Ad (V_\alpha(st) V_\alpha(r))$ for $s,t,r\in G$ shows that $\om_{\alpha}$ is a Borel $2$-cocycle on $G$ as 
in Example \ref{ex-om} and that $V_{\alpha}$ is a projective $\bar{\om}_{\alpha}$-representation 
of $G$ on $H$.

The class $[\om_\alpha]\in H^2(G,\TT)$
only depends on $\alpha$ and it vanishes if and only if 
 $\alpha$ is {\em unitary} in the sense that 
$\alpha$ is implemented by a strongly continuous homomorphism $V:G\to \U$. \footnote{To see
 this one should
  use the fact that any measurable homomorphism between polish groups is automatically 
continuous by \cite{Mo3}*{Proposition 5}.}
Therefore, the class $[\om_{\alpha}]\in H^2(G,\TT)$ is called the 
{\em Mackey obstruction for $\alpha$ being unitary}. An easy computation
gives:

\begin{lemma}\label{lem-cocycle}
Let $\alpha:G\to\Aut(\K(H))$, $V_{\alpha}:G\to \U(H)$ and $\om_{\alpha}$ be as above.
Let $G_{\om_{\alpha}}$ denote the central extension of 
$G$ by $\TT$ corresponding to $\om_\alpha$ as described in 
Example \ref{ex-om} and let $\iota:\TT\to\CC$ denote 
the inclusion. Let
$$\tilde{V}_{\alpha}:G_{\om_{\alpha}}\to\U(H);\quad \tilde{V}_{\alpha}(s,z)=\bar{z}V_{\alpha}(s).$$
Then $(H,\tilde V_{\alpha})$ is a $(G_{\om_{\alpha}},\TT)$-equivariant
Morita equivalence between the action $\alpha$ of 
$G\cong G_{\om_{\alpha}}/\TT$ on $\K(H)$ and the twisted action $(\id,\iota)$
of $(G_{\om_{\alpha}},\TT)$ on $\CC$.
\end{lemma}

We refer to \S \ref{sec-twistmor} for the definition of twisted
equivariant Morita equivalences.
Since Morita equivalent twisted systems have Morita equivalent full and reduced
crossed products, it follows that $\K(H)\rtimes_\alpha G$ is Morita equivalent
to the twisted group algebra $C^*(G,\om_{\alpha})$ (and similarly for
$\K(H)\rtimes_rG$ and $C_r^*(G,\om_\alpha)$). 
Recall from Example \ref{ex-om} that
there is a one-to-one correspondence between the representations of 
$C^*(G,\om_{\alpha})$ (or the covariant representations of $(\CC,G_{\om_{\alpha}},\TT,\id,\iota)$)
and the projective $\om_{\alpha}$-representations
of $G$. 
Using the above lemma and induction of covariant representations 
via the bimodule $(H,\tilde V_{\alpha})$ then gives:

\begin{theorem}\label{thm-cocycle}
Let $\alpha:G\to \Aut(\K(H))$ be an action and let $\om_{\alpha}$ and $V_\alpha:G\to \U(H)$ 
be as above. Then the assignment
$$L\mapsto (\id\otimes 1, V_{\alpha}\otimes L)$$
gives a homeomorphic bijection between 
the (irreducible) $\om_{\alpha}$-projective representations of $G$ 
and the (irreducible) non-degenerate covariant representations of $(\K(H),G,\alpha)$.
\end{theorem}

\begin{remark}\label{rem-littlegroup1}
{\bf (1)} It is actually quite easy to give a direct isomorphism between $C^*(G,\om_{\alpha})\otimes \K$
and the crossed product $\K\rtimes_{\alpha}G$, where we write $\K=\K(H)$. If $V_{\alpha}:G\to \U(H)$ is as above, then
one easily checks that 
$$\Phi: L^1(G,\om_{\alpha})\odot \K\to L^1(G,\K); \Phi(f\otimes k)(s)=f(s)kV_s^*.$$
is a $*$-homomorphism with dense range such that
$$(\id\otimes 1)\times (V_{\alpha}\otimes L)\big(\Phi(f\otimes k)\big)=L(f)\otimes k$$
for all $f\in L^1(G,\om)$ and $k\in \K$,
and hence the above theorem implies that $\Phi$ is isometric 
with respect to the $C^*$-norms. A similar argument also shows that
$\K\rtimes_rG\cong C_r^*(G,\om_{\alpha})\otimes \K$.

{\bf (2)} The separability assumptions made above are not really necessary: Indeed, 
if $\alpha:G\to\Aut(\K(H))\cong \PU(H)$ is an action of any locally compact group
on the algebra of compact operators on any Hilbert space $H$, 
then 
$$\tilde{G}:=\{(s, U)\in G\times \U(H): \alpha_s=\Ad(U)\}$$
fits into the central extension
$$
\begin{CD}
1 @>>>  \TT @>{z\mapsto (e,z\cdot 1)}>>\tilde{G} @>{(s,U)\mapsto s}>>G @>>>1.
\end{CD}
$$
If we define  $u:\tilde{G}\to \U(H); u(s,U)=U$,
then it is easy to check $(H,u)$ implements a Morita equivalence between 
$(\K(H), G,\alpha)$ and the twisted system $(\CC, \tilde{G}, \TT, \id, \iota)$.
Thus we obtain a one-to-one correspondence between the representations 
of $\K(H)\rtimes_\alpha G$ and the representations of $\CC\rtimes_{\id,\iota}(\tilde{G},\TT)$.
We refer to \cite{Green1}*{Theorem 18} for more details.
\end{remark}

Combining the previous results (and using the identification 
$\widehat{A}\cong\Prim(A)$ if $A$ is type I) with Theorem \ref{thm-MRG}
now gives:

\begin{theorem}[Mackey's little group method]\label{thm-Mackey3}
Suppose that $(A,G,\alpha)$ is a smooth separable system such that 
$A$ is type $I$. Let $S\subseteq \widehat{A}$ be a section for the 
quotient space $G\backslash\widehat{A}$ and for each $\pi\in S$
let $V_{\pi}:G_\pi\to \U(H_{\pi})$ be a measurable map such that 
$\pi(\alpha_s(a))=V_\pi(s)\pi(a)V_\pi(s)^*$ for all $a\in A$ and $s\in G_\pi$
(such map always exists). Let 
$\om_{\pi}\in Z^2(G_{\pi},\TT)$ be the $2$-cocycle satisfying
$$\om_{\pi}(s,t)\cdot 1_{H_\pi}:=V_{\pi}(st)V_{\pi}(t)^*V_{\pi}(s)^*.$$
Then
$$\IND:\cup_{\pi\in S}C^*(G_{\pi},\om_{\pi})\dach\to (A\rtimes G)\dach; 
\IND(L)=\ind_{G_{\pi}}^G(\pi\otimes 1)\times (V_{\pi}\otimes L). $$
is a bijection, which restricts to homeomorphisms between 
$C^*(G_{\pi},\om_\pi)\dach$ and it's image $(A_{G_{\pi}}\rtimes G)\dach$
for each $\pi\in  S$.
\end{theorem}

\begin{remark}\label{rem-Mackey3}
{\bf (1)}  If $G$ is exact, then a similar result holds for the reduced crossed product
$A\rtimes_rG$, if we also use the reduced twisted group algebras $C_r^*(G_\pi,\om_\pi)$
of the stabilisers.

{\bf (b)} If $(A,G,\alpha)$ is a  type I smooth system which is not separable, then a similar result 
can be formulated  using the approach described in part (2) of Remark \ref{rem-littlegroup1}.

 \end{remark}

Notice that the above result in particular applies to all systems 
$(A,G,\alpha)$ with $A$ type I and $G$ compact, since actions 
of compact groups on type I algebras are always smooth in the sense of 
Definition \ref{defn-smooth}. Since the central extensions $G_\om$
 of a compact group $G$ by $\TT$ are compact, and since 
 $C^*(G,\om)$ is a quotient of $C^*(G_\om)$ (see Example \ref{ex-om}),
 it follows that 
the twisted group algebras $C^*(G,\om)$ are direct sums of matrix algebras
if $G$ is compact. Using this,
we easily get from Theorem \ref{thm-Mackey3}:

 \begin{corollary}\label{cor-compact}
 Suppose that $(A,G,\alpha)$ is a system with $A$ type I and $G$
 compact. Then $A\rtimes G$ is type I. If, moreover, $A$ is CCR, then
 $A\rtimes G$ is CCR, too.
 \end{corollary}
 \begin{proof}
Since the locally closed subset 
$(A_{G_{\pi}} \rtimes G)\dach$ corresponding 
to some orbit $G(\pi)\subseteq\widehat{A}$ is homeomorphic (via Morita equivalence)
to $(\K(H_{\pi})\rtimes G)\dach \cong C^*(G_{\pi},\om_{\pi})\dach$,
it follows that $(A_{G(\pi)}\rtimes G)\dach$ is a discrete set in the induced
topology. This implies that all points in $(A\rtimes G)\dach$ 
are locally closed. Moreover, if $A$ is CCR, then the points in 
$\widehat{A}$ are closed. Since $G$ is compact, it follows then that
the $G$-orbits in $\widehat{A}$ are closed, too.
But then the discrete set 
$(A_{G_{\pi}} \rtimes G)\dach$ is closed in $(A\rtimes G)\dach$,
which implies that the points in  $(A\rtimes G)\dach$ are closed.
\end{proof}

\bibliography{references}

\begin{bibdiv}
\begin{biblist}

\bib{ADR}{book}{
      author={Anantharaman-Delaroche, C.},
      author={Renault, J.},
       title={Amenable groupoids},
      series={Monographie no. 36},
   publisher={L'Enseignement Math\'ematique},
     address={Geneva},
        date={2000},
}

\bib{AW}{article}{
      author={Akemann, C.},
      author={Weaver, N.},
       title={Consistency of a counterexample to {N}aimark's problem},
        date={2004},
     journal={Proc. Natl. Acad. Sci. USA},
      volume={101},
       pages={7522\ndash 7525},
}

\bib{BCL}{article}{
      author={Brodzki, J.},
      author={Cave, C.},
      author={Li, K.},
       title={Exactness of locally compact groups},
        date={2016},
     journal={preprint, arXiv:1603.01829},
         url={https://arxiv.org/abs/1603.01829},
}

\bib{BEW}{article}{
      author={Buss, A.},
      author={Echterhoff, S.},
      author={Willett, R.},
       title={{Exotic crossed products and the Baum-Connes conjecture}},
        date={2014},
     journal={preprint, arXiv:1409.4332},
         url={https://arxiv.org/abs/1409.4332},
}

\bib{BGR}{article}{
      author={Brown, L.},
      author={Green, P.},
      author={Rieffel, M.A.},
       title={Stable isomorphism and strong {M}orita equivalence of
  {$C^*$}-algebras},
        date={1977},
     journal={Pacific J. Math.},
      volume={71},
       pages={349\ndash 363},
}

\bib{Bla06}{book}{
      author={Blackadar, B.},
       title={{Operator algebras: Theory of C*-algebras and von Neumann
  algebras}},
      series={Encyclopaedia of Mathematical Sciences Vol. 122},
   publisher={Springer},
     address={Berlin},
        date={2006},
}

\bib{Blat1}{article}{
      author={Blattner, R.J.},
       title={On induced representations},
        date={1961},
     journal={Amer. J. Math.},
      volume={83},
       pages={79\ndash 98},
}

\bib{Blat2}{article}{
      author={Blattner, R.J.},
       title={{On a theorem of G. W. Mackey}},
        date={1962},
     journal={Bull. Amer. Math. Soc.},
      volume={68},
       pages={585\ndash 587},
}

\bib{BourT}{book}{
      author={Bourbaki, N.},
       title={{Topologie g\'en\'erale. Chapitres 1 \`a 4}},
      series={\'El\'ements de math\'ematique},
   publisher={Hermann},
     address={Paris},
        date={1971},
}

\bib{BuS}{article}{
      author={Busby, R.C.},
      author={Smith, H.A.},
       title={Representations of twisted group algebras},
        date={1970},
     journal={Trans. Amer. Math. Soc.},
      volume={149},
       pages={503\ndash 537},
}

\bib{BS2}{article}{
      author={Baaj, S.},
      author={Skandalis, G.},
       title={{Unitaires multiplicatifs et dualit\'e pour les produits
  crois\'es de $C^*$-alg\`ebres}},
        date={1993},
     journal={Ann. Sci. \'Ecole Norm. Sup.},
      volume={26},
       pages={425\ndash 488},
}

\bib{CE1}{article}{
      author={Chabert, J.},
      author={Echterhoff, S.},
       title={{Twisted equivariant $KK$-theory and the Baum-Connes conjecture
  for group extensions}},
        date={2001},
     journal={$K$-Theory},
      volume={23},
       pages={157\ndash 200},
}

\bib{CKRW}{article}{
      author={Crocker, D.},
      author={Kumjian, A.},
      author={Raeburn, I.},
      author={Williams, D.},
       title={An equivariant {B}rauer group and actions of groups on
  ${C}^*$-algebras},
        date={1997},
     journal={J. Funct. Anal.},
      volume={146},
       pages={151\ndash 184},
}

\bib{Combes}{article}{
      author={Combes, F.},
       title={Crossed products and {M}orita equivalence},
        date={1984},
     journal={Proc. London Math. Soc.},
      volume={49},
       pages={289\ndash 306},
}

\bib{Dav}{book}{
      author={Davidson, K.},
       title={${C}^*$-algebras by example},
      series={Fields Institute Monographs Vol. 6},
   publisher={American Mathematical Society},
     address={Providence, RI},
        date={1996},
}

\bib{DE}{book}{
      author={Deitmar, A.},
      author={Echterhoff, S.},
       title={Principles of harmonic analysis},
     edition={Second edition},
      series={Universitext},
   publisher={Springer},
     address={Cham},
        date={2014},
}

\bib{Dix}{book}{
      author={Dixmier, J.},
       title={${C}^*$-algebras},
      series={North-Holland Mathematical Library, Vol. 15},
   publisher={North-Holland Publishing Co.},
     address={Amsterdam-New York-Oxford},
        date={1977},
}

\bib{KK}{incollection}{
      author={Echterhoff, S.},
       title={Bivariant $kk$-theory and the baum-connes conjecure},
        date={2017},
   booktitle={{to appear in: $K$-theory for group C*-algebras and semigroup
  C*-algebras, Oberwolfach Seminars, Birkhäuser}},
      series={Oberwolfach Seminars},
   publisher={Birkh\"auser},
     address={Basel},
}

\bib{E-ind}{article}{
      author={Echterhoff, S.},
       title={On induced covariant systems},
        date={1990},
     journal={Proc. Amer. Math. Soc.},
      volume={108},
       pages={703\ndash 706},
}

\bib{E-prim}{article}{
      author={Echterhoff, S.},
       title={The primitive ideal space of twisted covariant systems with
  continuously varying stabilizers},
        date={1992},
     journal={Math. Ann.},
      volume={292},
       pages={59\ndash 84},
}

\bib{E-mor}{article}{
      author={Echterhoff, S.},
       title={{Morita equivalent twisted actions and a new version of the
  Packer-Raeburn stabilization trick}},
        date={1994},
     journal={J. London Math. Soc.},
      volume={50},
       pages={170\ndash 186},
}

\bib{E-ct}{article}{
      author={Echterhoff, S.},
       title={Crossed products with continuous trace},
        date={1996},
     journal={Mem. Amer. Math. Soc.},
      volume={123},
       pages={1\ndash 134},
}

\bib{EE}{article}{
      author={Echterhoff, S.},
      author={Emerson, H.},
       title={Structure and {K}-theory of crossed products by proper actions},
        date={2011},
     journal={Expositiones Mathematicae},
      volume={29},
       pages={300\ndash 344},
}

\bib{EH}{article}{
      author={Effros, E.G},
      author={Hahn, F.},
       title={Locally compact transformation groups and ${C}^*$-algebras},
        date={1967},
     journal={Mem. Amer. Math. Soc.},
      number={75},
       pages={1\ndash 92},
}

\bib{EKQR1}{article}{
      author={Echterhoff, S.},
      author={Kaliszewski, S.},
      author={Quigg, J.},
      author={Raeburn, I.},
       title={Naturality and induced representations},
        date={2000},
     journal={Bull. Austral. Math. Soc.},
      volume={61},
       pages={415\ndash 438},
}

\bib{EKQR2}{article}{
      author={Echterhoff, S.},
      author={Kaliszewski, S.},
      author={Quigg, J.},
      author={Raeburn, I.},
       title={{A Categorical Approach to Imprimitivity Theorems for
  $C^*$-Dynamical Systems}},
        date={2006},
     journal={Mem. Amer. Math. Soc.},
      volume={180},
       pages={viii+169},
}

\bib{EL}{article}{
      author={Echterhoff, S.},
      author={Laca, M.},
       title={The primitive ideal space of the ${C}^*$-algebra of the affine
  semigroup of algebraic integers},
        date={2013},
     journal={Math. Proc. Cambridge Philos. Soc.},
      volume={154},
       pages={119\ndash 126},
}

\bib{EN}{article}{
      author={Echterhoff, S.},
      author={Nest, R.},
       title={{The structure of the Brauer group and crossed products by
  $C_0(X)$-linear actions on $C_0(X,\mathcal{K})$}},
        date={2001},
     journal={Trans. Amer. Math. Soc.},
      volume={353},
       pages={3685\ndash 3712},
}

\bib{ERo}{article}{
      author={Echterhoff, S.},
      author={Rosenberg, J.},
       title={{Fine structure of the Mackey machine for actions of abelian
  groups with constant Mackey obstruction}},
        date={1995},
     journal={Pacific J. Math.},
      volume={170},
       pages={17\ndash 52},
}

\bib{ER}{article}{
      author={Echterhoff, S.},
      author={Raeburn, I.},
       title={The stabilisation trick for coactions},
        date={1996},
     journal={J. Reine Angew. Math.},
       pages={181\ndash 215},
}

\bib{EW1}{article}{
      author={Echterhoff, S.},
      author={Williams, D.P.},
       title={{Locally inner actions on $C_0(X)$-algebras}},
        date={2001},
     journal={J. Operator Th.},
      volume={45},
       pages={131\ndash 160},
}

\bib{EW-induced}{article}{
      author={Echterhoff, S.},
      author={Williams, D.P.},
       title={Inducing primitive ideals},
        date={2008},
     journal={Trans. Amer. Math. Soc.},
      volume={360},
       pages={6113\ndash 6129},
}

\bib{EW-proper}{article}{
      author={Echterhoff, S.},
      author={Williams, D.P.},
       title={Structure of crossed products by strictly proper actions on
  continuous-trace algebras},
        date={2014},
     journal={Trans. Amer. Math. Soc.},
      volume={366},
       pages={3649\ndash 3673},
}

\bib{EW2}{article}{
      author={Echterhoff, S.},
      author={Williams, D.P.},
       title={{Crossed Products by $C_0(X)$-actions}},
        date={1998},
     journal={J. Funct. Anal.},
      volume={158},
       pages={113\ndash 151},
}

\bib{Fell1}{article}{
      author={Fell, J.M.G.},
       title={The dual spaces of ${C}^*$-algebras},
        date={1960},
     journal={Trans. Amer. Math. Soc.},
      volume={94},
       pages={365\ndash 403},
}

\bib{Fell2}{article}{
      author={Fell, J.M.G.},
       title={Weak containment and induced representations of groups},
        date={1962},
     journal={Canad. J. Math},
      volume={14},
       pages={237\ndash 268},
}

\bib{Fol}{book}{
      author={Folland, G.B.},
       title={A course in abstract harmonic analysis},
      series={Studies in Advanced Mathematics},
   publisher={CRC Press},
     address={Boca Raton, FL},
        date={1995},
}

\bib{Ghys}{article}{
      author={Ghys, E.},
       title={{Groupes al\'eatoires (d'apr\'es Misha Gromov, ...)}},
        date={2004},
     journal={Ast\'erisque},
      volume={294},
       pages={173\ndash 204},
}

\bib{Gl}{article}{
      author={Glimm, J.},
       title={Locally compact transformation groups},
        date={1961},
     journal={Trans. Amer. Math. Soc.},
      volume={101},
       pages={124\ndash 138},
}

\bib{GPS}{inproceedings}{
      author={Giordano, T.},
      author={Putnam, I.F.},
      author={Skau, C.F.},
       title={{The orbit structure of Cantor minimal $\mathbb{Z}^2$-systems}},
        date={2006},
   booktitle={{Operator Algebras: The Abel Symposium 2004, Abel Symp., 1}},
   publisher={Springer},
     address={Berlin},
       pages={145\ndash 160},
}

\bib{GR}{article}{
      author={Gootman, E.C.},
      author={Rosenberg, J.},
       title={{The structure of crossed product $C^*$-algebras: a proof of the
  generalized Effros-Hahn conjecture}},
        date={1979},
     journal={Invent. Math.},
      volume={52},
       pages={283\ndash 298},
}

\bib{Green1}{article}{
      author={Green, P.},
       title={The local structure of twisted covariance algebras},
        date={1978},
     journal={Acta Math.},
      volume={140},
       pages={191\ndash 250},
}

\bib{Green2}{article}{
      author={Green, P.},
       title={The structure of imprimitivity algebras},
        date={1980},
     journal={J. Funct. Anal.},
      volume={36},
       pages={88\ndash 104},
}

\bib{Gr1}{article}{
      author={Gromov, M.},
       title={Spaces and questions},
        date={2000},
     journal={Geom. Funct. Anal.},
       pages={118\ndash 161},
}

\bib{IW}{article}{
      author={Ionescu, M.},
      author={Williams, D.P.},
       title={{The generalized Effros-Hahn conjecture for groupoids}},
        date={2009},
     journal={Indiana Univ. Math. J.},
      volume={58},
       pages={2489\ndash 2508},
}

\bib{K3}{article}{
      author={Kasparov, G.G.},
       title={{Equivariant $KK$-theory and the Novikov conjecture}},
        date={1988},
     journal={Invent. Math.},
      volume={91},
       pages={147\ndash 201},
}

\bib{K2}{inproceedings}{
      author={Kasparov, G.G.},
       title={{K-theory, group $C^*$-algebras, and higher signatures
  (conspectus)}},
        date={1995},
   booktitle={{Novikov conjectures, index theorems and rigidity, Vol. 1
  (Oberwolfach, 1993)}},
      volume={226},
   publisher={Cambridge Univ. Press},
     address={Cambridge},
       pages={101\ndash 146},
}

\bib{Kl}{article}{
      author={Kleppner, A.},
       title={Multipliers on abelian groups},
        date={1965},
     journal={Math. Ann.},
      volume={158},
       pages={11\ndash 34},
}

\bib{KV}{article}{
      author={Kustermans, J.},
      author={Vaes, S.},
       title={Locally compact quantum groups},
        date={2000},
     journal={Ann. Sci. {\'E}cole Norm. Sup.},
      volume={33},
       pages={837\ndash 934},
}

\bib{KW2}{article}{
      author={Kirchberg, E.},
      author={Wassermann, S.},
       title={Permanence properties of ${C}^*$-exact groups},
        date={2000},
     journal={Doc. Math.},
      volume={5},
       pages={513\ndash 558},
}

\bib{KW1}{article}{
      author={Kirchberg, E.},
      author={Wassermann, S.},
       title={Exact groups and continuous bundles of ${C}^*$-algebras},
        date={1999},
     journal={Math. Ann.},
      volume={315},
       pages={169\ndash 203},
}

\bib{Lan}{incollection}{
      author={Landsman, N.P.},
       title={Quantized reduction as a tensor product},
        date={2001},
   booktitle={Quantization of singular symplectic quotients},
      series={Progr. Math.},
      volume={198},
   publisher={Birkh\"auser},
     address={Basel},
       pages={137\ndash 180},
}

\bib{Lance}{book}{
      author={Lance, E.C.},
       title={Hilbert ${C}^*$-modules. {A} toolkit for operator algebraists},
      series={London Mathematical Society Lecture Note Series},
   publisher={Cambridge University Press},
     address={Cambridge},
        date={1995},
      volume={210},
}

\bib{Lep1}{article}{
      author={Leptin, H.},
       title={{Verallgemeinerte $L^1$-Algebren}},
        date={1965},
     journal={Math. Ann.},
      volume={159},
       pages={51\ndash 76},
}

\bib{LeG}{article}{
      author={Le~Gall, P.-Y.},
       title={{Th{\'e}orie de Kasparov {\'e}quivariante et groupo{\"\i}des I}},
        date={1999},
     journal={$K$-Theory},
      volume={16},
}

\bib{Ma1}{article}{
      author={Mackey, G.W.},
       title={On induced representations of groups},
        date={1951},
     journal={Amer. J. Math.},
      volume={73},
       pages={576\ndash 592},
}

\bib{Ma2}{article}{
      author={Mackey, G.W.},
       title={Induced representations of locally compact groups {I}},
        date={1952},
     journal={Annals of Math. (2)},
      volume={55},
       pages={101\ndash 139},
}

\bib{Ma3}{article}{
      author={Mackey, G.W.},
       title={{Induced representations of locally compact groups II. The
  Frobenius reciprocity theorem}},
        date={1953},
     journal={Annals of Math. (2)},
      volume={58},
       pages={193\ndash 221},
}

\bib{Ma4}{article}{
      author={Mackey, G.W.},
       title={Borel structure in groups and their duals},
        date={1957},
     journal={Trans. Amer. Math. Soc.},
      volume={85},
       pages={134\ndash 165},
}

\bib{Ma5}{article}{
      author={Mackey, G.W.},
       title={Unitary representations of group extensions. {I}},
        date={1958},
     journal={Acta Math.},
      volume={99},
       pages={265\ndash 311},
}

\bib{Mo1}{article}{
      author={Moore, C.C.},
       title={{Extensions and Low Dimensional Cohomology Theory of Locally
  Compact Groups. I}},
        date={1964},
     journal={Trans. Amer. Math. Soc.},
      volume={113},
       pages={40\ndash 63},
}

\bib{Mo2}{article}{
      author={Moore, C.C.},
       title={{Extensions and Low Dimensional Cohomology Theory of Locally
  Compact Groups. II}},
        date={1964},
     journal={Trans. Amer. Math. Soc.},
      volume={113},
       pages={64\ndash 86},
}

\bib{Mo3}{article}{
      author={Moore, C.C.},
       title={{Group extensions and cohomology for locally compact groups.
  III}},
        date={1976},
     journal={Trans. Amer. Math. Soc.},
      volume={221},
       pages={1\ndash 33},
}

\bib{Mo4}{article}{
      author={Moore, C.C.},
       title={{Group extensions and cohomology for locally compact groups.
  IV}},
        date={1976},
     journal={Trans. Amer. Math. Soc.},
      volume={221},
       pages={35\ndash 58},
}

\bib{Mur}{book}{
      author={Murphy, G.J.},
       title={${C}^*$-algebras and operator theory},
   publisher={Academic Press Inc.},
     address={Boston, MA},
        date={1990},
}

\bib{Osajda}{article}{
      author={Osajda, D.},
       title={Small cancellation labellings of some infinite graphs and
  applications},
        date={2014},
     journal={preprint, arXiv:1406.5015},
         url={https://arxiv.org/abs/1406.5015},
}

\bib{Oz}{inproceedings}{
      author={Ozawa, N.},
       title={Amenable actions and applications},
        date={2006},
   booktitle={{International Congress of Mathematicians}},
      volume={Vol. II},
   publisher={Eur. Math. Soc.},
     address={Z{\"u}rich},
       pages={1563\ndash 1580},
}

\bib{Pat}{book}{
      author={Paterson, A.L.T.},
       title={Amenability},
      series={Mathematical Surveys and Monographs},
   publisher={Amer. Math. Soc.},
     address={Providence, RI},
        date={1988},
      number={29},
}

\bib{Ped}{book}{
      author={Pedersen, G.K.},
       title={{$C^*$-algebras and their Automorphism Groups}},
   publisher={Academic Press},
     address={London},
        date={1979},
}

\bib{Po}{article}{
      author={Poguntke, Detlev},
       title={The structure of twisted convolution {$C^*$}-algebras on abelian
  groups},
        date={1997},
        ISSN={0379-4024},
     journal={J. Operator Theory},
      volume={38},
      number={1},
       pages={3\ndash 18},
      review={\MR{1462012}},
}

\bib{PR89}{article}{
      author={Packer, J.A.},
      author={Raeburn, I.},
       title={Twisted crossed products of ${C}^*$-algebras},
        date={1989},
     journal={Math. Proc. Cambridge Philos. Soc.},
      volume={106},
       pages={293\ndash 311},
}

\bib{QS}{article}{
      author={Quigg, J.},
      author={Spielberg, J.},
       title={Regularity and hyporegularity in ${C}^*$-dynamical systems},
        date={1992},
     journal={Houston J. Math},
      volume={18},
       pages={139\ndash 152},
}

\bib{Ren1}{book}{
      author={Renault, J.},
       title={A groupoid approach to ${C}^*$-algebras},
      series={Lecture Notes in Math.},
   publisher={Springer},
     address={Berlin},
        date={1980},
      volume={793},
}

\bib{Ren2}{article}{
      author={Renault, J.},
       title={Repr{\'e}sentations des produits crois{\'e}s d'alg{\`e}bres de
  groupo{\"i}des},
        date={1987},
     journal={J. Operator Th.},
       pages={67\ndash 97},
}

\bib{Rie1}{article}{
      author={Rieffel, M.A.},
       title={Induced representations of ${C}^*$-algebras},
        date={1974},
     journal={Adv. Math.},
      volume={13},
       pages={176\ndash 257},
}

\bib{Rie2}{incollection}{
      author={Rieffel, M.A.},
       title={{Unitary representations of group extensions; an algebraic
  approach to the theory of Mackey and Blattner}},
        date={1979},
   booktitle={Studies in analysis},
      series={Advances in Math. Supplementary Studies},
      volume={4},
   publisher={Academic Press},
     address={New York-London},
       pages={43\ndash 82},
}

\bib{RudRC}{book}{
      author={Rudin, W.},
       title={Real and complex analysis},
     edition={Third edition},
   publisher={Tata McGraw-Hill Book Co.},
     address={New York},
        date={1987},
}

\bib{RW}{book}{
      author={Raeburn, I.},
      author={Williams, D.P.},
       title={Morita equivalence and continuous-trace ${C}^*$-algebras},
      series={Mathematical Surveys and Monographs},
   publisher={Amer. Math. Soc.},
        date={1998},
      number={60},
}

\bib{S1}{article}{
      author={Sauvageot, J.-L.},
       title={Id{\'e}aux primitifs de certains produits crois{\'e}s},
        date={1977},
     journal={Math. Ann.},
      volume={231},
       pages={61\ndash 76},
}

\bib{S2}{article}{
      author={Sauvageot, J.-L.},
       title={Id{\'e}aux primitifs induits dans les produits crois{\'e}s},
        date={1979},
     journal={J. Funct. Anal.},
      volume={32},
       pages={381\ndash 392},
}

\bib{Sier}{article}{
      author={Sierakowski, A.},
       title={The ideal structure of reduced crossed products},
        date={2010},
     journal={M\"unster J. of Math.},
      volume={3},
       pages={237\ndash 262},
}

\bib{T1}{article}{
      author={Takesaki, M.},
       title={Covariant representations of ${C}^*$-algebras and their locally
  compact automorphism groups},
        date={1967},
     journal={Acta Math.},
      volume={119},
       pages={273\ndash 303},
}

\bib{Wi-crossed}{book}{
      author={Williams, D.P.},
       title={Crossed products of ${C}^*$-algebras},
      series={Mathematical Surveys and Monographs},
   publisher={Amer. Math. Soc.},
        date={2007},
      number={134},
}

\bib{Wi0}{article}{
      author={Williams, D.P.},
       title={{The topology on the primitive ideal space of transformation
  group $C^*$-algebras and CCR transformation group $C^*$-algebras}},
        date={1981},
     journal={Trans. Amer. Math. Soc.},
      volume={266},
       pages={335\ndash 359},
}

\end{biblist}
\end{bibdiv}
\end{document}